\documentclass[11pt,a4paper]{article}
\pdfoutput=1
\usepackage[margin=2.5cm]{geometry} 
\usepackage{amssymb}
\usepackage{amsmath}
\usepackage{mathtools}
\usepackage{amsthm}
\usepackage{amsfonts}
\usepackage{kpfonts} 
\usepackage{bbm} 
\usepackage{dutchcal} 
\usepackage[mathscr]{euscript} 
\usepackage{tikz-cd} 
\usepackage{todonotes} 
\usepackage{graphicx} 
\usepackage[toc,page]{appendix} 
\usepackage{xcolor} 
\usepackage{float} 
\usepackage[inline]{enumitem} 
\usepackage{comment} 
\usepackage{authblk} 
\usepackage[normalem]{ulem} 

\usepackage{url}
\usepackage[style = alphabetic,
            backref = true,
            backrefstyle = two]{biblatex}
\renewcommand\mathfrak[1]{\mbox{\usefont{U}{euf}{m}{n}#1}}

\definecolor{winered}{rgb}{0.5,0,0}
\usepackage[citecolor = winered, 
            urlcolor = winered, 
            linkcolor = winered,
            menucolor = winered, 
            colorlinks = true]{hyperref}

\DefineBibliographyStrings{english}{%
  backrefpage = {p.}, 
  backrefpages = {pp.}, 
}

\theoremstyle{definition}
\newtheorem{Th}{Theorem}[subsection]
\newtheorem{Lemma}[Th]{Lemma}
\newtheorem{Cor}[Th]{Corollary}
\newtheorem{Prop}[Th]{Proposition}
\newtheorem{Def}[Th]{Definition}

\newtheorem{Rem}[Th]{Remark}
\newtheorem{Ex}[Th]{Example}

\newcommand{\N}{\mathbb{N}}
\newcommand{\ncat}{\mathbf} 
\newcommand{\cat}{\mathcal} 
\renewcommand{\u}{\underline}
\newcommand{\colim}{\text{colim}}

\newcommand{\op}{\text{op}} 

\newcommand{\sd}{\mathfrak{D}}
\newcommand{\const}{\text{cst}}

\newcommand{\tw}{\textbf{tw}}
\newcommand{\Gr}{\ncat{Gr}}
\newcommand{\dGr}{\ncat{dGr}}
\newcommand{\mlt}{\text{mlt}}
\newcommand{\G}{\mathcal{G}}
\newcommand{\Lyr}{\text{Lyr}}
\newcommand{\Gaif}{\text{Gaif}}
\newcommand{\Hyp}{\ncat{Hyp}}

\definecolor{diagramBlue}{RGB}{94, 164 ,255}
\definecolor{diagramGold}{RGB}{255,162,23}

\usetikzlibrary{calc}
\usetikzlibrary{decorations.pathmorphing}

\tikzset{curve/.style={settings={#1},to path={(\tikztostart)
    .. controls ($(\tikztostart)!\pv{pos}!(\tikztotarget)!\pv{height}!270:(\tikztotarget)$)
    and ($(\tikztostart)!1-\pv{pos}!(\tikztotarget)!\pv{height}!270:(\tikztotarget)$)
    .. (\tikztotarget)\tikztonodes}},
    settings/.code={\tikzset{quiver/.cd,#1}
        \def\pv##1{\pgfkeysvalueof{/tikz/quiver/##1}}},
    quiver/.cd,pos/.initial=0.35,height/.initial=0}

\tikzset{tail reversed/.code={\pgfsetarrowsstart{tikzcd to}}}
\tikzset{2tail/.code={\pgfsetarrowsstart{Implies[reversed]}}}
\tikzset{2tail reversed/.code={\pgfsetarrowsstart{Implies}}}

\definecolor{emilioeditcolor}{rgb}{0.94, 0.97, 1.0}

\definecolor{beneditcolor}{rgb}{204, 104, 0}

\definecolor{zoltaneditcolor}{rgb}{250, 128, 114}

\title{Structured Decompositions: Structural and Algorithmic Compositionality.}

\author{Benjamin Merlin Bumpus\thanks{Instituto de Matemática e Estatística, Universidade de São Paulo. Rua do Matão, 1010 — 05508–090, São
Paulo, SP, Brasil. This author acknowledges ERC funding under the European Union's Horizon 2020 research and innovation program (grant \# 803421, ``ReduceSearch'') which was received while at the Department of Mathematics and Computer Science of Eindhoven University of Technology. benjamin.merlin.bumpus(at)gmail.com}, \, Zoltan A. Kocsis\thanks{School of Computer Science and Engineering, University of New South Wales, z.kocsis(at)unsw.edu.au },\\ Jade Edenstar Master\thanks{Independent researcher, jadeedenstarmaster(at)gmail.com} ~\& Emilio Minichiello\thanks{CUNY CityTech, eminichiello(at)gmail.com}}

\date{Last compilation: \today}

\begin{document}

\maketitle

\begin{abstract}
We introduce \textit{structured decompositions}, category-theoretic structures which simultaneously generalize notions from \textit{graph theory} (including treewidth, layered treewidth, co-treewidth, graph decomposition width, tree independence number, hypergraph treewidth and $H$-treewidth), \textit{geometric group theory} (specifically Bass-Serre theory), and \textit{dynamical systems} (e.g. hybrid dynamical systems). We define sd-functors, which provide a compositional way to analyze and relate different structural complexity measures, and establish a general duality between decompositions and completions of objects. 
\end{abstract}

\setcounter{tocdepth}{2} 
\tableofcontents

\section{Introduction}
Compositionality can be understood as the perspective that the semantics, structure or function of the whole should be given by that of its constituent parts~\cite{compositionality-SEP}.  This principle has long shaped thinking in mathematics and computer science and indeed basic tools like recursion and divide-and-conquer algorithms themselves rely on compositional reasoning. 

More recently, a strong research effort has concentrated around the systematic mathematical study of compositional systems and their occurrences in wider scientific domains~\cite{fong2016algebra, pollard2017open, cicala2019rewriting, courser2020open, master2021composing}. Thanks to these efforts, we can now make sense of how one might build graphs, Petri nets~\cite{open-petri-nets-master-baez}, chemical reaction networks~\cite{reaction-networks-pollard-baez}, stock and flow diagrams~\cite{open-stock-and-flow} or epidemiological models~\cite{libkind2022algebraic} from smaller constituent parts. Ideally, we want versatile algorithms which consider the categorial and compositional structure of their inputs, and work generically across all these instances. Compositional algorithms at this level of generality don't yet exist. This article presents a starting point for their development, based upon a rich, well-established theory for computing on graphs.

The field of parameterized complexity supplies a number of established techniques for building compositional algorithms well-suited for particular combinatorial objects that exhibit certain compositional structure. The classes of relevant objects are often identified via so-called ``width measures'': numerical quantities originating from structural and algorithmic graph theory which roughly can be thought as measuring compositional complexity. Treewidth, the most well-known among these, measures roughly how much the global connectivity of a graph differs from that of a tree. The focus of the present paper is to develop a general theory of those width measures which arise from colimits.

Our article continues the work of \cite{bumpus2021generalizing} and \cite{bumpus2023spined}, of studying category-theoretic generalizations of treewidth. We do this by introducing \textit{structured decomposition categories} (which we also refer to as \textit{sd-categories} for convenience, c.f. Definition \ref{def structured decomposition category}). These are categories $\cat{C}$ equipped with a given family of diagrams all of whose colimits exist in $\cat{C}$. The diagrams in question, which we call structured decompositions\footnote{We invite the computationally-minded reader to furthermore experiment with structured decompositions via their implementation \href{https://github.com/AlgebraicJulia/StructuredDecompositions.jl}{StructuredDecompositions.jl} in Algebraic Julia.}, are always of a specific type and roughly they can be understood as functorial ways of attaching data to combinatorial objects. This paper thus investigates which combinatorial invariants can be expressed as colimits: in a structured decomposition category one sees an object as being decomposed if it arises as a colimit of a structured decomposition. We equip sd-categories with additional structure that interacts nicely with structured decomposition diagrams, resulting in \textit{width categories} (c.f. Definition \ref{def width category}).

Slightly more concretely, structured decompositions are special kinds of diagrams in some fixed category. The following is a generic example of what these diagrams look like. 
\[\begin{tikzcd}
	{{\color{diagramGold}\bullet}} & {{\color{diagramGold}\bullet}} & {{\color{diagramGold}\bullet}} & {{\color{diagramGold}\bullet}} \\
	\\
	{{\color{diagramBlue}\bullet}} & {{\color{diagramBlue}\bullet}} & {{\color{diagramBlue}\bullet}} & {{\color{diagramBlue}\bullet}} & {{\color{diagramBlue}\bullet}}
	\arrow[from=1-1, to=3-1]
	\arrow[from=1-1, to=3-2]
	\arrow[from=1-2, to=3-1]
	\arrow[from=1-2, to=3-3]
	\arrow[from=1-3, to=3-2]
	\arrow[from=1-3, to=3-5]
	\arrow[from=1-4, to=3-3]
	\arrow[from=1-4, to=3-4]
\end{tikzcd}\]
Intuitively, a structured decomposition is a way of assigning objects of a category to the vertices and edges of a graph. One starts with a graph \(G\) and constructs a category \(\smallint G\) as follows: \(\smallint G\) has an object for each vertex of \(G\) (blue dots in the diagram above), an object for each edge of \(G\) (gold dots) and a span joining each edge to its source and target vertices\footnote{This is an instance of the more general notion of a Grothendieck construction.}. Then structured decompositions are functors of the form \(\smallint G \to \cat{C}\) for some category $\cat{C}$. 

As we already mentioned, we view structured decompositions as generalizations of graph decompositions. In graph theory, decompositions often come equipped with a notion of width -- a numerical quantity measuring the structural complexity of the pieces involved in the decomposition -- which itself then induces an invariant on the graphs being decomposed. For instance, consider treewidth, a very important invariant in both graph structure theory (e.g. Robertson and Seymour's celebrated graph minor theorem~\cite{RobertsonXX}) and in parameterized complexity (e.g. Courcelle's famous algorithmic meta-theorem for bounded tree-width graphs~\cite{courcelle1990monadic}): one says that a graph has treewidth at most $k$ if it admits a tree decomposition of width at most $k$. In the same spirit, if $\Gamma$ is a width category (c.f. Definition \ref{def width category}), we associate a corresponding notion of width, which we call $\Gamma$-width $\mathbf{w}_\Gamma(X)$, to each object $X \in \Gamma$.

Let us now summarize the results of this paper. In Section \ref{section sd-categories}, we introduce the main categorical machinery we will use, namely sd-categories, width categories and sd-functors. In Proposition \ref{prop Gamma width is monotone with respect to monomorphisms}, we prove that $\Gamma$-width is monotone with respect to monomorphisms. In Proposition \ref{prop sd-functors give inequality on width}, we prove that if $F : \Gamma \to \Delta$ is a functor satisfying certain easy-to-check conditions (Definition \ref{def sd functors}), then $F$ cannot increase width (i.e. for every $X \in \Gamma$, $\mathbf{w}_{\Delta}(F(X)) \leq \mathbf{w}_\Gamma(X)$). Inspired by chordal completions of graphs, we also introduce a second width measure one can associate to a width category (chordal width), and in Theorem \ref{th completions equiv to decomps in stable sd-categories}, we prove that chordal width and $\Gamma$-width agree for so-called complete width categories.

In Section \ref{section examples}, we prove that a large number of examples from the literature appear as particular instances of our framework. These include: 
\begin{enumerate*}
    \item \textit{treewidth}~\cite{bertele1972, halin1976s, robertsonII} (Proposition \ref{prop sd-category for treewidth}), 
    \item \textit{complemented treewidth}~\cite{sousa2021} (Proposition \ref{prop sd-category for co-treewidth}), 
    \item \textit{tree independence number}~\cite{dallard2024treeindependence} (Proposition \ref{prop sd-category for tree independence number}), 
    \item \textit{hypergraph treewidth}~\cite{heinz2013tree} (Proposition \ref{prop sd-category for hypergraph treewidth}), 
    \item \textit{layered treewidth}~\cite{shahrokhi2015layered, dujmovic2017layered} (Proposition \ref{prop sd-category for layered treewidth}) and 
    \item \textit{$\mathcal{H}$-treewidth}~\cite{jansen2022vertex} (Proposition \ref{prop sd-category for H-treewidth}).
\end{enumerate*} We explore relationships between these width categories and obtain bounds between their $\Gamma$-widths (Corollaries \ref{cor tree independence number and co-treewidth}, \ref{cor tree independence number and vertex cover number}, \ref{cor hypergraph treewidth equal to treewidth of Gaifman graph}).

Moreover, instances of our notion also turn up outside of combinatorial settings, When instantiated in the category of groups, structured decompositions coincide with `graphs of groups', a concept developed in the 1970s and central to Bass-Serre theory~\cite{serre1970groupes, serre2002trees, BASS19933, Higgins-fund-groupoid}.
When dealing with manifolds and hybrid dynamical systems, structured decompositions appear independently under the name of ``hybrid objects'' in Ames' PhD thesis~\cite{ames2006categorical}. We treat these notions in Section~\ref{section examples}. The fact that the same formalism shows up in combinatorics, geometric group theory, and hybrid systems supports the need for a general theory of structured decompositions and width measures built from colimits: that is precisely the present contribution.

One should not expect all combinatorial decomposition methods to arise as colimits. Colimits have a strong topological flavor, and there are combinatorial decomposition methods such as clique-width decomposition trees~\cite{courcelle1993, COURCELLE199687}, and rank decompositions~\cite{oum2006approximating} which do not display this topological kind of compositionality. For instance clique-width decomposition trees are defined via a grammar which allows joining two graphs together by adding edges between them: it is unclear how such an operation could be naturally described as a colimit. Section \ref{section future work} discusses the research problem of capturing these methods, along with other open questions. The appendix (Section \ref{section categories of graphs}) treats the properties of several different categories of graphs and hypergraphs in detail.

\subsection{Related Work}
Recent efforts in graph theory have sought to generalize tree decompositions in two different ways. The first considers more general decomposition ``shapes'', such as cycle and planar decompositions; this is  best exemplified by Carmesin's work on graph decompositions~\cite{carmesin-local2, diestel-homotopy-groups}. The second approach is to work with tree-shaped decompositions, but to allow for more complex  ``bags''; examples include $\mathcal{H}$-treewidth~\cite{H-tw-conference} and layered treewidth~\cite{shahrokhi2015layered, dujmovic2017layered}. Our notion of structured decompositions bridges and unifies both approaches, promising exciting new avenues of research at the intersection of categorial and combinatorial ideas.  

There have been a few categorifications of treewidth in the literature so far; one of which is the starting point for the present paper. This is the notion of spined categories~\cite{bumpus2021generalizing, bumpus2023spined}, categories $\cat{C}$ equipped with a so called proxy-pushout operation and a sequence $\Omega: \mathbb{N}_= \to \cat{C}$ of objects called the spine. This paper originated as a generalization of spined categories to the case where the spine is filtration on the objects of $\cat{C}$, rather than being a sequence. This change allows us to recover many more examples of combinatorial invariants. 

Treewidth has also received categorical treatment through the notion of monoidal width. This was introduced by Di Lavore and Soboci{\'n}ski~\cite{monoidal-width, monoidal-width-rkw} and it measures the width of a morphism in a monoidal category $(\cat{C},\otimes)$. Instead of a spine $s : \cat{N} \to \cat{C}$, monoidal width is determined by a width function $w : \cat{A} \subset \mathsf{Mor} \cat{C} \to \mathbb{N}$ sending a subclass of atomic morphisms to natural numbers. There are two main differences between our approach and that of Di Lavore and Soboci{\'n}ski. Firstly, for monoidal width, the morphisms are being decomposed whereas in the the present structured decomposition approach we decompose the objects. Setting that aside, the only other difference is that the spine and width function are going in opposite directions.  Indeed, a spine $s$ may be obtained from a width function $w$ by choosing $s(n) \in w^{-1}(n)$. 

The idea that the treewidth -- specifically -- can be encoded by taking pushouts has already been noted by Blume et. al.~\cite{blume2011treewidth}. Those authors define a cospan decomposition of a graph $H$ to be a sequence of connected cospans whose colimit is $H$. This notion is perhaps conceptually simpler than structured decompositions (indeed it is an instantiation, or a special case of our notion) but it is deployed only in the specific case of graphs and treewidth. In contrast, as we've already mentioned, our focus is one of developing a general theory that encapsulates many notions at once. Furthermore, a benefit of viewing combinatorial decompositions as diagrams, as we do, is that one can speak of morphisms between decompositions and functorial relationships between various notions of width. 

More broadly, on the topic of diagrammatic reasoning, note that although mathematicians have considered categories of diagrams (in the sense of Definition \ref{def diagram category}) since the beginnings of category theory \cite{eilenberg1945general}, apart from a handful of examples~\cite{kock1967limit, guitart1973foncteur, guitart1974remarques, guitart1977decompositions}, interest in their study has waned over time, but has picked up recently ~\cite{peschke2020diagrams, perrone2022kan, diagrammaticequations}. Structured decompositions, being particular kinds of diagrams, assemble into a subcategory of the category of diagrams. This further justifies, independently of combinatorial consideration, a systematic study of the kinds of objects that can be built via colimits.

Finally we should point out that structured decompositions bear significant similarity to undirected wiring diagrams~\cite{spivak2013operad}. These provide an operadic view on a construction which is very similar to what we would call $\ncat{FinSet}$-valued structured decompositions. Although it is beyond the scope of the article, investigating the connections between these two notions is a promising direction for further study.

\subsection{Notation}
We use the notation $\mathcal{C}, \mathcal{D}$ for generic categories and the bold font $\mathbf{Set}$ for named categories. We assume all categories are locally small, and we shall refer to those that have class-sized hom-sets as huge categories. We let $\ncat{Cat}$ denote the category of small categories and $\ncat{CAT}$ denote the huge category of categories. We let $\ncat{Set}$ denote the category of sets and $\ncat{FinSet}$ the full subcategory of finite sets. We fix a small skeleton for $\ncat{FinSet}$, and abuse notation by also calling that category $\ncat{FinSet}$. By a finite set, we then mean an object in this skeleton. This allows us to say that the class of all finite sets is a set. If $S$ is a set, then we let $P(S)$ denote the power set of $S$, i.e. the set of subsets of $S$, and we let $P_{\neq \varnothing}(S)$ denote the set of nonempty subsets of $S$.

\subsection{Acknowledgements}
We thank Karl Heuer, Reinhard Diestel, Bart Jansen, James Fairbanks, Elena Di Lavore, and Pawel Soboci{\'n}ski for their helpful conversations which aided in the development of this work. We are grateful to Will Turner for spotting a slight inaccuracy in a previous version of this article and for pointing out the connection to Bass-Serre theory. Furthermore we would like to thank Evan Patterson for suggesting to use the Grothendieck construction in the definition of structured decompositions; this has considerably simplified many of the technical arguments as well as the presentation of the paper. Finally we thank the anonymous reviewers for their constructive and detailed feedback. 

\section{Structured Decomposition Categories} \label{section sd-categories}

In this section we first define structured decompositions, the main categorical tool we will use in this paper to study various notions of width in the literature. After this we define the central notions of this paper, sd-categories, $\Gamma$-width and sd-functors. We then review the classical notion of treewidth from graph theory.

\subsection{The Category of Graphs}

Graphs notions vary throughout the literature, depending on the context and application. In what follows, we consider what the \textit{nLab}~\cite{nlab:graph} refers to as simple graphs, which have no loops and no multiple edges between the same pair of vertices. See Section \ref{section categories of graphs} for a discussion of different categories of graphs.

\begin{Def} \label{def graphs}
A \textbf{graph} $G$ consists of the following data:
\begin{itemize}
    \item a finite set $V(G)$ of \textbf{vertices},
    \item a symmetric, irreflexive binary \textbf{edge relation} $E(G) \subseteq V(G)^2$. 
\end{itemize}
A \textbf{graph homomorphism} $f: G \rightarrow H$ between the graphs $G$ and $H$ is a function $f: V(G) \rightarrow V(H)$ that preserves the edge relation: if $(x,y) \in E(G)$, then $(f(x), f(y)) \in E(H)$. Let $\ncat{Gr}$ denote the category whose objects are graphs, and whose morphisms are graph homomorphisms.
\end{Def}

We will frequently abuse notation and write $x \in G$ to mean that $x \in V(G)$, and say that there is an edge $xy$ when $(x,y) \in E(G)$. Since the edge relation is symmetric, we will often identify an edge $xy$ with its symmetry class $\{xy,yx\}$.

The following classes of graphs will be used throughout.

\begin{Def} \label{def complete graphs, cycles, trees and paths}
The \textbf{complete graph} on $n$ vertices, also known as the \textbf{$n$-clique}, denoted $K^n$ is the graph with $V(K^n) = \{1, 2, \dots, n \}$ and where every pair of vertices are connected by an edge. By convention we take $K^0 = \varnothing$.
\begin{figure}[H]
    \centering
    \begin{equation*}
    \begin{tikzcd}
	&&&& \bullet && \bullet & \bullet \\
	\bullet & \bullet & \bullet & \bullet && \bullet & \bullet & \bullet
	\arrow[no head, from=1-5, to=2-6]
	\arrow[no head, from=1-7, to=1-8]
	\arrow[no head, from=1-7, to=2-7]
	\arrow[no head, from=1-7, to=2-8]
	\arrow[no head, from=1-8, to=2-7]
	\arrow[no head, from=1-8, to=2-8]
	\arrow[no head, from=2-2, to=2-3]
	\arrow[no head, from=2-4, to=1-5]
	\arrow[no head, from=2-4, to=2-6]
	\arrow[no head, from=2-7, to=2-8]
\end{tikzcd}
\end{equation*}
    \caption{The first four nonempty complete graphs.}
\end{figure}
The $n$-\textbf{cycle} $C^n$ is a graph with $n$-many vertices $\{1, \dots, n \}$ where there is an edge between $i$ and $j$ if and only if $j = i + 1$ for $1 \leq i \leq n-1$ or $i = n$ and $j = 1$. A \textbf{circuit of length $n$} in a graph $G$ is a morphism $f : C^n \to G$, which we commonly denote by the image of the vertices $\{f(1), f(2), \dots, f(n) \}$. An $n$-cycle in a graph $G$ is a graph homomorphism $f : C^n \to G$, injective on vertices and edges, for some $n \geq 3$. 

Let $P^n$ denote the graph with $V(P^n) = \{1, \dots, n \}$ and such that $i$ and $j$ are adjacent if and only if $|j - i| = 1$. We call $P^n$ an \textbf{$n$-path}. Let $\text{Paths}$ denote the set of all $n$-paths for $n \geq 0$. Given $x,y \in V(G)$, a path connecting $x$ to $y$ in $G$ is a graph homomorphism $f : P^n \to G$, injective on vertices and edges, such that $f(1) = x$ and $f(n) = y$ for some $n \geq 1$. We say that a graph is \textbf{connected} if any two of its vertices have a path connecting them. We say a graph is a \textbf{tree} if it is connected and has no cycles. We let $\text{Trees}$ denote the set of all trees.
\end{Def}

\subsection{Barycentric Subdivision}

\begin{Def}\label{def barycentric subdivision}
Consider a graph $G$. The \textbf{barycentric subdivision} $\smallint G$ of $G$ is the category which has as
\begin{itemize}
    \item \textbf{objects}: all the vertices of $G$ (\textit{vertex objects}) and all symmetry classes of edges of $G$ (\textit{edge objects})
    \item \textbf{morphisms}: identity morphisms for each object, and two morphisms $e_x : e \rightarrow x$ and $e_y : e \rightarrow x$ for each edge object $e=xy$.
\end{itemize}
Note that composition in  $\smallint G$ is trivial, since one of the morphisms in each composable pair is an identity morphism.
\end{Def}

\begin{Ex}
The barycentric subdivision of the graph $G$ on the left yields the category $\smallint G$ on the right.
\begin{equation*}
\begin{tikzcd}
	x &&& x && z \\
	y & z && e & y & f
	\arrow["e"', no head, from=1-1, to=2-1]
	\arrow["f", no head, from=2-2, to=2-1]
	\arrow["{{e_x}}", from=2-4, to=1-4]
	\arrow["{{e_y}}"', from=2-4, to=2-5]
	\arrow["{{f_z}}"', from=2-6, to=1-6]
	\arrow["{{f_y}}", from=2-6, to=2-5]
\end{tikzcd}
\end{equation*}
\end{Ex}

\subsection{Structured Decompositions}

\begin{Def} \label{def structured decomposition}
Given a category $\cat{C}$ and a graph $G$, a $G$-\textbf{structured decomposition} is a functor $d: \smallint G \rightarrow \cat{C}$ such that the image of each morphism in $\smallint G$ is a monomorphism in $\cat{C}$. We call
\begin{itemize}
    \item images $d(v)$ of vertex objects $v \in \smallint G$ the \textbf{bags},
    \item images $d(e)$ of edge objects $e \in \smallint G$ the \textbf{adhesions},
    \item spans of the form $d(v_1) \hookleftarrow d(e) \hookrightarrow d(v_2)$ the \textbf{adhesion-spans} indexed by the edge object~$e$.
\end{itemize}
\end{Def}

\begin{Ex} \label{ex finset structured decomp}
For $\cat{C} = \ncat{FinSet}$ and $G$ a graph as below left, we have a structured decomposition below right.
\begin{equation*}
\begin{tikzcd}
	&& \bullet &&& {\{1,2,3\}} & {\{1,2,3, 4\}} & {\{3,4\}} \\
	& \bullet && \bullet & {\{1, 2\}} & {\{1,2,3,\pi, \sqrt{2}\}} & {\{3, \pi\}} & {\{3,4,10\}} \\
	\bullet && \bullet && {\{1,2, 3/2\}} && {\{3,\pi, 11\}}
	\arrow[no head, from=1-3, to=2-2]
	\arrow[no head, from=1-3, to=2-4]
	\arrow[hook', from=1-6, to=1-7]
	\arrow[hook, from=1-6, to=2-6]
	\arrow[dashed, no head, from=1-7, to=2-6]
	\arrow[dashed, no head, from=1-7, to=2-8]
	\arrow[hook, from=1-8, to=1-7]
	\arrow[hook', from=1-8, to=2-8]
	\arrow[no head, from=2-2, to=3-1]
	\arrow[no head, from=2-2, to=3-3]
	\arrow[hook', from=2-5, to=2-6]
	\arrow[hook, from=2-5, to=3-5]
	\arrow[dashed, no head, from=2-6, to=3-5]
	\arrow[dashed, no head, from=2-6, to=3-7]
	\arrow[hook, from=2-7, to=2-6]
	\arrow[hook', from=2-7, to=3-7]
\end{tikzcd}
\end{equation*}
\end{Ex}

\begin{Def}
Take a category $\cat{C}$, a graph $G$, an object $X \in \cat{C}$, and a structured decomposition $d : \smallint G \to \cat{C}$. If $X \cong \colim\:d$, we say that $d$ is a \textbf{structured decomposition of the object~$X$}.
\end{Def}

As a shorthand, we sometimes refer to structured decompositions of an object $X$ simply as decompositions of $X$.

\begin{Ex}
The structured decomposition of Example \ref{ex finset structured decomp} is a structured decomposition of the finite set $\{1,2,3,4,3/2, \pi, \sqrt{2}, 10, 11\}$.
\end{Ex}

\subsection{Grothendieck Construction}

\begin{Def} \label{def diagram category}
Given a category $\cat{C}$, let $\ncat{Diag}(\cat{C})$ denote the category whose objects are (small) diagrams $d : I \to \cat{C}$ and for two diagrams $d : I \to \cat{C}$ and $d' : I' \to \cat{C}$, a morphism $A : d \to d'$ in $\ncat{Diag}(\cat{C})$ consists of a functor $a : I \to I'$ and a natural transformation $\alpha : d \Rightarrow A d'$
\begin{equation*}
    \begin{tikzcd}
	I && {I'} \\
	& {\cat{C}}
	\arrow["a", from=1-1, to=1-3]
	\arrow[""{name=0, anchor=center, inner sep=0}, "d"', from=1-1, to=2-2]
	\arrow[""{name=1, anchor=center, inner sep=0}, "{d'}", from=1-3, to=2-2]
	\arrow["\alpha", shorten <=7pt, shorten >=7pt, Rightarrow, from=0, to=1]
\end{tikzcd}
\end{equation*}
We call this the \textbf{diagram category}\footnote{The diagram category construction is well understood, see \cite{tholen2021diagrams}, and for further applications see \cite{patterson2021diagrams}.} of $\cat{C}$. Given a functor $F : \cat{C} \to \cat{D}$, we obtain a functor $\ncat{Diag}(F) : \ncat{Diag}(\cat{C}) \to \ncat{Diag}(\cat{D})$ by postcomposing diagrams with $F$ and by sending a natural transformation as below left to the whiskered natural transformation as below right.
\begin{equation*}
\begin{tikzcd}
	I && {I'} & I && {I'} \\
	& {\cat{C}} &&& {\cat{D}}
	\arrow["a", from=1-1, to=1-3]
	\arrow[""{name=0, anchor=center, inner sep=0}, "d"', from=1-1, to=2-2]
	\arrow[""{name=1, anchor=center, inner sep=0}, "{d'}", from=1-3, to=2-2]
	\arrow["a", from=1-4, to=1-6]
	\arrow[""{name=2, anchor=center, inner sep=0}, "Fd"', from=1-4, to=2-5]
	\arrow[""{name=3, anchor=center, inner sep=0}, "{Fd'}", from=1-6, to=2-5]
	\arrow["\alpha", shorten <=7pt, shorten >=7pt, Rightarrow, from=0, to=1]
	\arrow["{F\alpha}", shorten <=7pt, shorten >=7pt, Rightarrow, from=2, to=3]
\end{tikzcd}
\end{equation*}
Thus we obtain a functor $\ncat{Diag} : \ncat{CAT} \to \ncat{CAT}$. We let $\mathfrak{D}(\cat{C})$ denote the full subcategory of $\ncat{Diag}(\cat{C})$ whose objects are structured decompositions in $\cat{C}$. If $\G$ is a set of graphs, then we let $\sd_{\G}(\cat{C})$ denote the full subcategory of the $\G$-structured decompositions. If $X \in \cat{C}$, then we let $\sd_\G(\cat{C}, X)$ denote the full subcategory on $\G$-structured decompositions of $X$.
\end{Def}

\begin{Rem}
In later sections, we will allow for more general notions of graphs, in which case our notion of graph above will be a special case. If $G$ is a multigraph (Definition \ref{def multigraphs}), then let $\smallint G$ denote the category defined just as above, but now with a span $x \leftarrow e \rightarrow y$ for each of the multiple possible edges between $x$ and $y$. This construction is functorial, and we use the notation $\smallint : \Gr_\mlt \to \ncat{Cat}$ for the resulting functor. We note that for a loop, we define its image in the barycentric subdivision to be a span $x \leftarrow e\rightarrow x$. By the full subcategory inclusions $\Gr, \Gr_\ell, \Gr_r \hookrightarrow \Gr_\mlt$, we can also index structured decompositions using any of the kinds of undirected graphs considered in Section \ref{section overview categories of graphs}. 

We can also index structured decompositions using digraphs (Definition \ref{def directed multigraphs}). If $G$ is a digraph, equivalently a functor $G: \ncat{dGrSch} \to \ncat{Set}$, then $\smallint G$ is precisely the \textbf{category of elements} or \textbf{Grothendieck construction}\footnote{See \cite[Section 2.4]{riehl2017category}.} of $G$.
We note that if $G$ is a digraph, and $U(G)$ denotes its underlying undirected multigraph, then there is an isomorphism of categories $\smallint G \cong \smallint U(G)$. Furthermore and generalizing the case of digraphs, one can clearly also index structured decompositions by objects of any category of presheaves or C-sets (e.g. Petri nets, simplicial sets, etc.). We will not be considering structured decompositions whose shapes are not graphs in the present paper; but we simply point out this possibility for completeness.
\end{Rem}

Let $\const : \cat{C} \to \sd(\cat{C})$ denote the constant diagram functor, which sends an object $X \in \cat{C}$ to the structured decomposition of $X$ given by the functor $X :  (\ncat{1} \cong \smallint \bullet) \to \cat{C}$ that sends the single object of the terminal category $\ncat{1}$ to $X$.

To conclude this section, let us characterize the category of structured decompositions in $\ncat{Set}$. Consider the (strict) functor $J \colon \Gr^\op \to \ncat{Cat}$ defined as:
\begin{align*}
    J \colon \Gr^\op &\to \ncat{Cat} \\
    G &\mapsto \ncat{Set}^{\smallint G} \\
    (G \xrightarrow{f} H) &\mapsto \bigl (\ncat{Set}^{\smallint H} \xrightarrow{\ncat{Set}^{\smallint f}} \ncat{Set}^{\smallint G}\bigr).
\end{align*}
Taking the Grothendieck construction \cite[Chapter 10]{johnson2021} of this functor, we obtain a Grothendieck fibration $\pi : \smallint J \to \Gr^\op$, where $\smallint J$ is the category whose objects are pairs $(G, d : \smallint G \to \ncat{Set})$ and whose morphisms are pairs $(f : G \to H, \eta : d \Rightarrow (d' \circ \smallint f))$, i.e. diagrams of the form
\begin{equation*}
    \begin{tikzcd}
	{\smallint G} && {\smallint H} \\
	& {\ncat{Set}}
	\arrow["{\smallint f}", from=1-1, to=1-3]
	\arrow[""{name=0, anchor=center, inner sep=0}, "d"', from=1-1, to=2-2]
	\arrow[""{name=1, anchor=center, inner sep=0}, "{d'}", from=1-3, to=2-2]
	\arrow["\eta", shorten <=8pt, shorten >=8pt, Rightarrow, from=0, to=1]
\end{tikzcd}
\end{equation*}
We recognize this as precisely the category of structured decompositions in $\ncat{Set}$.

\begin{Lemma}
There is an isomorphism of categories
\begin{equation*}
    \sd(\ncat{Set}) \cong \smallint J.
\end{equation*}
\end{Lemma}

\subsection{Width Categories}

In prior work, the first two of the authors introduced the formalism of spined categories~\cite{bumpus2023spined} to give a uniform account of several width notions in the graph theory literature, including treewidth, co-treewidth, and hypergraph treewidth. One can regard the definitions below as advancing this formalism by generalizing the spine category and requiring actual colimits rather than proxy pushouts. This allows us to capture a larger variety of examples from the literature.

\begin{Def} \label{def structured decomposition category}
A \textbf{structured decomposition category} or \textbf{sd-category} $\Gamma = (\cat{C}, \G)$ consists of a small category $\cat{C}$ equipped with a set $\G$ of graphs\footnote{We also allow for the indexing graphs to be multigraphs or digraphs to accommodate for examples like those in Section \ref{section caremsin's graph decompositions}. We consider graphs, loop graphs, and reflexive loop graphs as special kinds of multigraphs, see Section \ref{section overview categories of graphs}. For most examples of interest however $\G = \text{Trees}$.} containing the graph $\bullet$, called the \textbf{index set} such that 
\begin{enumerate}[left = \parindent + 1em]
\item[(\textbf{W1})] for every structured decomposition $d : \smallint G \to \cat{C}$, with $G \in \G$, the colimit of $d$ exists in $\cat{C}$.
\end{enumerate}
\end{Def}

We now need some additional structure on an sd-category with which we can use to obtain some notion of ``size'' of an object in the category.

\begin{Def} \label{def spined sd-category}
A \textbf{spined sd-category} $\Gamma = (\cat{C}, \G, \Omega)$ consists of an sd-category $(\cat{C}, \G)$, equipped with a filtered essentially full subcategory\footnote{In other words, if $X \in \Omega$ and $X \cong Y$, then $Y \in \Omega$.} $\Omega$ called the \textbf{spine}. By filtered we mean $\Omega$ is equipped with a filtration $\Omega_0 \subseteq \Omega_1 \subseteq \dots \subseteq \Omega$ of essentially full subcategories of the spine with
\begin{equation*}
    \Omega = \bigcup_{n \geq 0} \, \Omega_n,
\end{equation*}
such that
\begin{enumerate}[left = \parindent + 1em]
    \item[(\textbf{W2a})] for every object $X \in \cat{C}$ there is a monomorphism $X \hookrightarrow W$ for some $W \in \Omega$, and 
    \item[(\textbf{W2b})] for every $n \geq 1$ and $W \in (\Omega_n \setminus \Omega_{n-1})$, there exists no monomorphism $W \hookrightarrow W'$ into any $W' \in \Omega_{n-1}$.
\end{enumerate}

We write (\textbf{W2}) to mean the logical conjunction (\textbf{W2}) $=$ (\textbf{W2a}) $\wedge$ (\textbf{W2b}). We say an object $W \in \cat{C}$ is \textbf{complete} if $W \in \Omega$. If each $\Omega_n \setminus \Omega_{n-1}$ is a singleton $\{ W_n \}$, then we will often write $\Omega$ as $\{ W_n \}_{n \geq 0}$. We may also write $\Omega = \{ W_i \}_{i \in I}$ if $I$ is a convenient indexing set.
\end{Def}

\begin{Rem}
The main difference between spined sd-categories and the spined categories of \cite{bumpus2023spined} is that the latter's notion of spine only allows for $(\Omega_n \setminus \Omega_{n-1})$ to be a singleton, and our notion requires certain kinds of colimits rather than proxy pushouts.
\end{Rem}

\begin{Def} \label{def size function of a spined sd-category}
Given a spined sd-category $\Gamma = (\cat{C}, \G, \Omega)$, let $s : \text{Obj}(\cat{C}) \to \mathbb{N}$ be the function, which we call the \textbf{size function} of $\Gamma$, defined as follows. Given $X \in \cat{C}$, we say that the size $s(X)$ of $X$ is $n$ if $\Omega_n$ is the smallest full subcategory of $\Omega$ such that there exists a monomorphism $X \hookrightarrow W$ with $W \in \Omega_n$.
\end{Def}

\begin{Lemma} \label{lem size monotonic for spined sd-categories}
If $\Gamma = (\cat{C}, \G, \Omega)$ is a spined sd-category and $f : X \hookrightarrow Y$ is a monomorphism in $\cat{C}$, then $s(X) \leq s(Y)$.
\end{Lemma}

\begin{proof}
If $Y \hookrightarrow W$ is monic with $W \in \Omega$, then the composite $X \hookrightarrow Y \hookrightarrow W$ is monic and witnesses $s(X) \leq s(Y)$.
\end{proof}

\begin{Rem}
By Lemma \ref{lem size monotonic for spined sd-categories} say that the size $s$ of a spined sd-category $\Gamma$ is monotonic with respect to monomorphisms.
\end{Rem}

\begin{Def} \label{def width}
Given a spined sd-category $\Gamma = (\cat{C}, \G, \Omega)$, an object $X \in \cat{C}$, and a $\G$-structured decomposition $d: \smallint G \to \cat{C}$ of $X$, we define the \textbf{width} $w(d)$ of the decomposition to be the maximal size of all of its bags minus $1$,
\begin{equation*}
w(d) = \max_{x \in V(G)} s(d(x)) - 1.
\end{equation*}
We say that a $\G$-structured decomposition $d$ of $X$ is \textbf{minimum}, if $w(d) \leq w(d')$ for all other $\G$-structured decompositions $d'$ of $X$.

Given an object $X \in \cat{C}$ we say that its \textbf{$\Gamma$-width} -- written $\mathbf{w}_\Gamma(X)$ -- is the width of any one of its minimum decompositions, or  equivalently: 
\begin{equation*}
    \mathbf{w}_\Gamma(X) = \min_{d \in \mathfrak{D}_{\G}(\cat{C}, X)} w(d).
\end{equation*}
If the ambient sd-category $\Gamma$ is clear from context, then we will refer to the $\Gamma$-width of an object $X$ simply as its width.
\end{Def}

\begin{Ex}
If we take $\Gamma = (\Gr, \text{Trees}, \{ K^n \}_{n \geq 0})$, (see Definition \ref{def complete graphs, cycles, trees and paths}) and if $G \in \Gr$, then the $\Gamma$-width of $G$ is precisely its treewidth. We will prove this later: see Proposition \ref{prop sd-category for treewidth}.
\end{Ex}

\begin{Ex}
The poset $(\N, \leq)$ of the natural numbers forms a spined sd-category of the form $\Gamma = ((\N, \leq), \text{Graphs}, \{ n \}_{n \geq 0} )$. The colimit of a finite diagram here is the maximum, hence any structured decomposition $d$ of an object $n$ must include $n$ as a bag $d(x)$, and furthermore this must be an object of maximal size in the decomposition. Hence the $\Gamma$-width of any number $n$ is $n-1$.
\end{Ex}

\begin{Lemma}
Given a spined sd-category $\Gamma = (\cat{C}, \G, \Omega)$ and an object $X \in \cat{C}$,
\begin{equation*}
    \mathbf{w}_\Gamma(X) \leq s(X) - 1.
\end{equation*}
\end{Lemma}

\begin{proof}
Given any $X \in \cat{C}$, the $\G$-structured decomposition $d : \smallint \bullet \to \cat{C}$ that picks out $X$ has width $s(X) - 1$.
\end{proof}

\begin{Def} \label{def stable sd-category}
If $\Gamma = (\cat{C}, \G)$ is an sd-category, and $\cat{F}$ is a pullback-stable\footnote{We say that $\cat{F}$ is pullback-stable if for any $f : X \to Y$ in $\cat{F}$ and any map $g : Z \to Y$, the pullback $g^*(f) : Z \times_Y X \to Z$ exists and belongs to $\cat{F}$.} class of morphisms in $\cat{C}$, then we say that $\Gamma$ is $\cat{F}$-\textbf{stable} if it additionally satisfies
\begin{enumerate}[left = \parindent + 1em]
    \item [(\textbf{W3})] the pullback $f^*(d)$ of the colimit cocone of any $\G$-structured decomposition $d$ of an object $Y$ along a morphism $f : X \to Y$ with $f \in \cat{F}$ exists and is a colimit cocone of $X$.
\end{enumerate}
We say that $\Gamma$ is \textbf{stable} if it is $\text{Mor}(\cat{C})$-stable, and \textbf{monic-stable} if it is $\text{Mono}(\cat{C})$-stable.
\end{Def}

\begin{Rem} \label{rem universal colimits}
The condition (\textbf{W3}) is a special case of having \textbf{universal} or \textbf{pullback-stable colimits}. More precisely, given a category $\cat{C}$ with pullbacks, suppose that $\mathcal{I}$ is a class of (small) diagrams $d : I \to \cat{C}$ whose colimit exists in $\cat{C}$. We say that $\cat{C}$ has universal $\mathcal{I}$-colimits if pulling back a colimit cocone of a diagram $d \in \mathcal{I}$ over an object $Y$ along an arbitrary map $f : X \to Y$ is a colimit cocone over $X$.

We say that a category has universal colimits if it has universal $\mathcal{I}$-colimits for $\mathcal{I}$ the class of all small diagrams. Every locally cartesian closed category has universal colimits. The class of locally cartesian closed categories includes quasitoposes and toposes. Given a set $\G$ of graphs, we say that a category $\cat{C}$ has \textbf{universal $\G$-colimits}, if it has universal $\mathcal{I}$-colimits, where $\mathcal{I}$ is the class of all $\G$-structured decomposition diagrams $d : \smallint G \to \cat{C}$ with $G \in \G$. 
\end{Rem}

We repackage the above remark with the following lemma.

\begin{Lemma}
Suppose that $\mathcal{I}$ is a class of (small) diagrams and $\cat{C}$ is a category with pullbacks and universal $\mathcal{I}$-colimits. By choosing representatives for every pullback, one obtains a well-defined presheaf 
\begin{equation*}
    \text{Dgm}_\mathcal{I} : \cat{C}^\op \to \ncat{Set},
\end{equation*}
where for $X \in \cat{C}$, $\text{Dgm}_\mathcal{I}(X)$ is the set of diagrams $d : I \to \cat{C}$ in $\mathcal{I}$ equipped with a colimit cocone $\sigma : d \Rightarrow \Delta(X)$. Given a map $f : X \to Y$ in $\cat{C}$, $\text{Dgm}_\mathcal{I}(f)$ is the function that takes a colimit cocone $\sigma$ over $Y$ to its pullback $f^*(\sigma)$ over $X$.
\end{Lemma}

\begin{Ex}
Any finitely complete, cocomplete and locally cartesian closed category $\cat{C}$ is a stable spined sd-category if we set $\Omega_0 = \varnothing$, $\Omega_1 = \dots = \Omega = \cat{C}$ and $\G = \text{Graphs}$. In this case all objects then have width $0$.
\end{Ex}

We now single out the class of spined sd-categories that are monic-stable. These are the kinds of sd-categories that are well-behaved and appear most frequently in examples.

\begin{Def} \label{def width category}
We call a spined sd-category $\Gamma = (\cat{C}, \G, \Omega)$ a \textbf{width category} if it is monic-stable.
\end{Def}

\begin{Ex}
The main example of a width category for us is $(\Gr, \text{Trees}, \{ K^n \}_{n \geq 0})$. We prove that this is indeed a width category in Proposition \ref{prop sd-category for treewidth}.
\end{Ex}

However we note that the size and width of even a complete object can have arbitrarily large difference.

\begin{Ex}
The category of finite sets forms a width category $(\ncat{FinSet}, \text{Paths}, \{ \u{n} \}_{n \geq 0})$, where $\u{n} = \{1, 2, \dots, n \}$. The size of a finite set is its cardinality, but the width of $\u{1}$ is clearly $0$, and the width of every other finite set is $1$, because we can always obtain a $\text{Paths}$-structured decomposition of $\u{n}$ of the form
\begin{equation*}
    \begin{tikzcd}
	& {\{2\}} && {\{3\}} && {\{n-1\}} \\
	{\{1,2\}} && {\{2,3\}} && \dots && {\{n-1,n\}}
	\arrow[hook', from=1-2, to=2-1]
	\arrow[hook, from=1-2, to=2-3]
	\arrow[hook', from=1-4, to=2-3]
	\arrow[hook, from=1-4, to=2-5]
	\arrow[hook', from=1-6, to=2-5]
	\arrow[hook, from=1-6, to=2-7]
\end{tikzcd}
\end{equation*}
\end{Ex}

\begin{Prop} \label{prop Gamma width is monotone with respect to monomorphisms}
Given a width category $\Gamma = (\cat{C}, \G, \Omega)$, if $ f: X \hookrightarrow Y$ is monic, then
\begin{equation*}
    \mathbf{w}_\Gamma(X) \leq \mathbf{w}_\Gamma(Y).
\end{equation*}
\end{Prop}

\begin{proof}
Suppose that $d : \smallint G \to \cat{C}$ is a $\G$-structured decomposition of $Y$ of minimal width. So $w(d) = \mathbf{w}_\Gamma(Y)$. Then since $\Gamma$ is monic-stable, $f^*(d)$ is a $\G$-structured decomposition of $X$ and since monomorphisms are stable under pullback, we obtain monomorphisms $f^*(d)(x) \hookrightarrow d(x)$ between the bags for every $x \in V(G)$. So by Lemma \ref{lem size monotonic for spined sd-categories}, $s(f^*(d)(x)) \leq s(d(x))$ for all $x \in V(G)$. Thus $w(f^*(d)) \leq w(d) = \mathbf{w}_\Gamma(Y)$, and therefore $\mathbf{w}_\Gamma(X) \leq \mathbf{w}_\Gamma(Y)$.
\end{proof}

\begin{Rem}
    Observe that the definition of a spined sd-category (Definition~\ref{def spined sd-category}) forces every object to have a \textit{finite} size. This is actually not due to mathematical necessity, but purely a stylistic choice. The reader can verify that, by considering ordinal-indexed filtrations as the spine and by changing the definition of width (Definition~\ref{def width}) to simply be the maximum bag size (i.e. removing the minus one), one obtains a theory of spined sd-categories that does not impose any finiteness conditions on bags of decompositions. These considerations, although relevant to tree-width of infinite graphs, will not be further addressed in the present contribution.
\end{Rem}

\subsection{Chordal Completions}
Let us now introduce a competing notion of width for width categories, using Lemma \ref{lem tree width given by chordal completion} as our motivation.

\begin{Def} \label{def chordal width}
Given a spined sd-category $\Gamma = (\cat{C}, \G, \Omega)$, if $d : \smallint G \to \cat{C}$ is a $\mathcal{G}$-structured decomposition such that each of its bags $d(x)$ belong to $\Omega$, then we say that $d$ is a \textbf{$(\mathcal{G}, \Omega)$-structured decomposition}. We let $\sd_{(\G, \Omega)}(\cat{C})$ and $\sd_{(\G,\Omega)}(\cat{C}, X)$ denote the full subcategories of $(\G, \Omega)$-structured decompositions in $\cat{C}$ and of an object $X \in \cat{C}$ respectively.

A \textbf{chordal object} is an object $Y \in \cat{C}$ which is the colimit of a $(\G, \Omega)$-structured decomposition. The \textbf{$(\G,\Omega)$-width} $w_{(\G,\Omega)}(d)$ of a $(\G,\Omega)$-decomposition $d$ of $Y$ is $n - 1$ if $\Omega_n$ is the smallest full subcategory that the bags of $d$ are objects of. The $(\G,\Omega)$-width $\mathbf{w}_{(\G,\Omega)}(Y)$ of a chordal object $Y$ is the minimal $(\G,\Omega)$-width of all of its $(\G,\Omega)$-decompositions.

A \textbf{chordal completion} of an object $X$ is a monomorphism $i : X \hookrightarrow Y$ where $Y$ is chordal. We define the \textbf{chordal $\Gamma$-width} of an object $X$, denoted $\mathbf{cw}_\Gamma(X)$, to be the minimal $(\G,\Omega)$-width of every chordal completion of $X$. When $\Gamma$ is clear from context, we may refer to chordal $\Gamma$-width simply as chordal width.
\end{Def}

\begin{Rem} \label{rem G omega width}
We note that $(\G, \Omega)$-width $\mathbf{w}_{(\G, \Omega)}(Y)$ is only defined for chordal objects $Y$ in a spined sd-category. We will see in Section \ref{section bass-serre theory} that even though $(\G, \Omega)$-width is not defined for all objects, it can still be a useful width measure in categories where monomorphisms don't behave like they do in $\Gr$.
\end{Rem}

\begin{Rem}
Given a width category $\Gamma = (\cat{C}, \G, \Omega)$, an object $Y \in \cat{C}$ is chordal if and only if it is in the essential image of the left Kan extension
\begin{equation*}
    \begin{tikzcd}
	\Omega & {\cat{C}} \\
	{\mathfrak{D}_{(\G, \Omega)}(\cat{C})}
	\arrow["i", hook, from=1-1, to=1-2]
	\arrow["{\text{cst}}"', from=1-1, to=2-1]
	\arrow["{\text{Lan}_{\text{cst}}i}"', from=2-1, to=1-2]
\end{tikzcd}
\end{equation*}
\end{Rem}

Now given a width category $\Gamma = (\cat{C}, \G, \Omega)$, and $X \in \cat{C}$, suppose that $f : X \hookrightarrow Y$ is a chordal completion of an object $X$. Then by definition there is a $(\G,\Omega)$-structured decomposition $d : \smallint G \to \cat{C}$ of $Y$, and we can pull back this decomposition along $f$ to get a $\G$-structured decomposition $f^*(d)$ of $X$. Furthermore we obtain a morphism of structured decompositions
\begin{equation*}
    \begin{tikzcd}
	{\smallint G} && {\smallint G} \\
	& {\cat{C}}
	\arrow[Rightarrow, no head, from=1-1, to=1-3]
	\arrow[""{name=0, anchor=center, inner sep=0}, "{d}"', from=1-1, to=2-2]
	\arrow[""{name=1, anchor=center, inner sep=0}, "{f^*(d)}", from=1-3, to=2-2]
	\arrow["\eta", shorten <=7pt, shorten >=7pt, Rightarrow, from=0, to=1]
\end{tikzcd}
\end{equation*}
given by mapping each component along the pullback maps. 

Conversely, given a $\G$-structured decomposition $d$ of an object $X$, by (\textbf{W2}) there exists a monomorphism $d(x) \hookrightarrow W_x$ for each bag of the decomposition, and so we can consider the $\G$-structured decomposition $d'$ obtained by replacing each bag $d(x)$ by their corresponding object $W_x$ and the adhesion maps $d(e) \hookrightarrow d(x)$ by the composite morphisms $d(e) \hookrightarrow W_x$. We call this completing the $\G$-structured decomposition $d$. In more detail we have the following.

\begin{Def}
Let $\Gamma = (\cat{C}, \G, \Omega)$ be a width category with $X \in \cat{C}$ and let $d : \smallint G \to \cat{C}$ be a $\G$-structured decomposition of $X$. Suppose that $d' : \smallint G \to \cat{C}$ is a diagram of the same shape. We say that a natural transformation $\alpha : d \Rightarrow d'$ is a \textbf{completion} of $d$ if each bag in $d'$ is a complete object of $\cat{C}$, each component map $\alpha_{e} : d(e) \to d'(e)$ between the adhesions is an identity map, and the component maps $\alpha_x : d(x) \to d'(x)$ between the bags are monomorphisms.
\end{Def}

If we let $Y$ denote the colimit of $d'$, then we obtain an induced map $f : X \to Y$ on the colimits. If the maps $d(x) \hookrightarrow W_x$ are minimal in the sense that $d(x)$ does not have a monomorphism to a complete object of size strictly smaller than the size of $W_x$, then we say that the completion $d'$ is \textbf{minimal}.

We wish to be able to go between $\G$-structured decompositions and chordal completions of objects, as is the case with Lemma \ref{lem tree width given by chordal completion}. 

However, there is no reason why the map $f : X \to Y$ one obtains by completing a $\G$-structured decomposition will be a monomorphism. In order to get this correspondence to behave correctly, we ask for a stronger property from our width categories.

\begin{Def} \label{def complete width category}
We say that a width category $\Gamma = (\cat{C}, \G, \Omega)$ is \textbf{complete}\footnote{Not to be confused with having all small limits.} if 
\begin{enumerate}[left = \parindent + 1em]
    \item [(\textbf{W4})] for every completion $\alpha : d \Rightarrow d'$ of a $\G$-structured decomposition $d$, the induced map $f : X \rightarrow Y$ between their colimits is a monomorphism.
\end{enumerate}
\end{Def}

Adhesive categories are categories where pushouts are particularly well behaved~\cite{lack2004adhesive}. Many familiar categories of combinatorial data are adhesive, including all presheaf categories. The following result establishes a connection between such categories and complete width categories defined thereupon. 

\begin{Prop}\label{prop adhesive cats with tree-shaped decomps are complete}
    Let $\Gamma = (\cat{C}, \text{Trees}, \Omega)$ be a width category. If $\cat{C}$ is adhesive, then $\Gamma$ is complete.
\end{Prop}
\begin{proof}
    Let $\delta \colon \smallint T \to \cat{C}$ denote a completion of any given tree-shaped structured decomposition $d \colon \smallint T \to \cat{C}$. Our argument proceeds by induction on the edges of $T$. If $T$ has no edges, the claim is trivially true since the only completion arrow is also the colimit arrow. 
    
    For the inductive step, consider an edge $e = xy$ of $T$ which splits $T$ into two subtrees $T_1$ and $T_2$. By the induction hypothesis, there are monomorphisms 
    \[g_i \colon \colim \bigl(\smallint T_i \hookrightarrow \smallint T \xrightarrow{d} \cat{C} \bigr) \hookrightarrow \colim \bigl(\smallint T_i \hookrightarrow \smallint T \xrightarrow{\delta} \cat{C} \bigr) \text{ where }i \in \{1,2\}.\]
    Since colimits of tree-shaped decompositions can be computed recursively (one pushout at a time, for each adhesion span) it suffices to show that the dashed puhout arrow in the comutative cube shown below is monic. 

\[\begin{tikzcd}
	{\colim \bigl(\smallint T_1 \hookrightarrow \smallint T \xrightarrow{d} \cat{C} \bigr)} && {d(e)} \\
	& {\colim(d)} && {\colim \bigl(\smallint T_2 \hookrightarrow \smallint T \xrightarrow{d} \cat{C} \bigr)} \\
	{\colim \bigl(\smallint T_1 \hookrightarrow \smallint T \xrightarrow{\delta} \cat{C} \bigr)} && {d(e)} \\
	& {\colim(\delta)} && {\colim \bigl(\smallint T_2 \hookrightarrow \smallint T \xrightarrow{\delta} \cat{C} \bigr)}
	\arrow[from=1-1, to=2-2]
	\arrow["{g_x}", color={rgb,255:red,214;green,92;blue,92}, hook, from=1-1, to=3-1]
	\arrow[hook, from=1-3, to=1-1]
	\arrow[hook, from=1-3, to=2-4]
	\arrow["\cong"{description, pos=0.7}, color={rgb,255:red,214;green,92;blue,92}, no head, from=1-3, to=3-3]
	\arrow["\lrcorner"{anchor=center, pos=0.125, rotate=90}, draw=none, from=2-2, to=1-3]
	\arrow[dashed, from=2-2, to=4-2]
	\arrow[from=2-4, to=2-2]
	\arrow["{g_y}"', color={rgb,255:red,214;green,92;blue,92}, hook, from=2-4, to=4-4]
	\arrow[from=3-1, to=4-2]
	\arrow[hook, from=3-3, to=3-1]
	\arrow[hook, from=3-3, to=4-4]
	\arrow["\lrcorner"{anchor=center, pos=0.125, rotate=90}, draw=none, from=4-2, to=3-3]
	\arrow[from=4-4, to=4-2]
\end{tikzcd}\]

    To that end, observe that the top and bottom faces of the above cube are pushouts (by construction) and that the back faces must necessarily be pullbacks since the adhesion spans of the completion $\delta$ are constructed by composition. Thus, since $\cat{C}$ is adhesive, the previous observations amount to saying that the above cube in van Kampen~\cite[Definition 2.1]{lack2004adhesive} and hence the front faces are pushouts. This concludes the proof since, $\cat{C}$ being adhesive,  pushout squares of monomorphisms are also van Kampen (by definition) and hence preserve monormophisms~\cite[Lemma 2.3]{lack2004adhesive}.
\end{proof}

It now follows immediately that if $(\cat{C}, \G, \Omega)$ is a complete width category, and $d : \smallint G \to \cat{C}$ is a $\G$-structured decomposition of an object $X$, then completing $d$ results in a $(\G, \Omega)$-structured decomposition $d'$ and by taking colimits, we obtain a chordal completion $i : X \hookrightarrow Y$. Conversely, if we have a chordal completion $i : X \hookrightarrow Y$, then by definition $Y$ has a $(\G, \Omega)$-structured decomposition $d$, and pulling back along $i$ we obtain a decomposition $i^*(d)$ of $X$. 

However, unlike with Lemma \ref{lem chordal completion and tree decomp adjunction}, there is no canonical choice for how to complete a structured decomposition. However, every chordal completion does arise from a $\G$-structured decomposition. Indeed if $i : X \hookrightarrow Y$ is a chordal completion, then if we pull back any decomposition $d$ of $Y$, we get a decomposition, $i^*(d)$, and the colimit of the map $i^*(d) \Rightarrow d$ is isomorphic to the completion $i$. Furthermore, we have the following result, generalizing Lemma \ref{lem tree width given by chordal completion}.

\begin{Th} \label{th completions equiv to decomps in stable sd-categories}
Given a complete width category $\Gamma = (\cat{C}, \G, \Omega)$ with $X \in \cat{C}$, we have
\begin{equation*}
    \mathbf{w}_\Gamma(X) = \mathbf{cw}_\Gamma(X).
\end{equation*}
\end{Th}

\begin{proof}
Suppose that $d_0 : \smallint G \to \cat{C}$ is a minimal $\G$-structured decomposition of $X$, in the sense that it has the smallest possible size of all $\G$-structure decompositions of $X$. Then by the discussion above, we can obtain a minimal completion $d'_0 : \smallint G \to \cat{C}$ of $d_0$, whose colimit we denote $Y$. This provides a chordal completion $f : X \hookrightarrow Y$. Since $d'_0$ is a minimal completion of $d_0$, this means that $w(d_0) = w_{(\G,\Omega)}(d'_0)$. Since $\mathbf{w}_{(\G,\Omega)}(Y) \leq w_{(\G,\Omega)}(d'_0)$, we have that $\mathbf{w}_{(\G,\Omega)}(Y) \leq w(d_0)$.

But since $Y$ is chordal, then there exists a $(\G, \Omega)$-decomposition $d'_1$ that achieves its width $\mathbf{w}_{(\G,\Omega)}(Y)$. By pulling this decomposition back, we obtain a decomposition $d_1$ of $X$. By definition of the width of $d_1$, we have $w(d_1) \leq w(d'_1) = \mathbf{w}_{(\G,\Omega)}(Y)$. But $w(d_0) \leq w(d_1)$, since we assumed that $d_0$ is minimal. Hence $\mathbf{w}_\Gamma(X) = w(d_0) = \mathbf{w}_{(\G,\Omega)}(Y)$. Therefore $\mathbf{w}_\Gamma(X) \geq \mathbf{cw}_\Gamma(X)$.

Now suppose that $f : X \to Y'$ is a minimal chordal completion of $X$, so that $\mathbf{cw}_\Gamma(X) = \mathbf{w}_{(\G, \Omega)}(Y')$. Then if we pull back a minimal decomposition $d'_2$ of $Y'$, i.e. $w_{(\G,\Omega)}(d'_2) = \mathbf{w}_{(\G,\Omega)}(Y')$, we obtain a decomposition $d_2$ of $X$ such that $w(d_2) \leq \mathbf{w}_{(\G,\Omega)}(Y')$. Then we have
\begin{equation*}
    \mathbf{w}_\Gamma(X) \leq w(d_2) \leq \mathbf{w}_{(\G,\Omega)}(Y') = \mathbf{cw}_{\Gamma}(X).
\end{equation*}
Thus the $\Gamma$-width and chordal $\Gamma$-width of $X$ agree.
\end{proof}

Although it is beyond the scope of the current paper, Theorem~\ref{th completions equiv to decomps in stable sd-categories} is particularly relevant to algorithmic applications of structured decompositions. To see this, we will briefly recall why chordal completions are useful when dealing with graphs and their tree decompositions. In algorithmic applications~\cite{flum2006parameterized, cygan2015parameterized, courcelle2012book, GroheBook} one is not only interested in knowing whether a given graph or structure has bounded treewidth, say, but instead one is interested in algorithms that, given an input graph $G$, will compute a decomposition of $G$ of small width. The most famous such algorithm is due to Bodlaender and Kloks~\cite{BodlaenderKloks1996} (see Althaus and Ziegler for a simplified exposition~\cite{althaus2021optimal}) which builds on Perković and Reed's earlier algorithm~\cite{PerkovicReedTreeDecompAlg}. This algorithm is an FPT-time algorithm parameterized by treewidth which decides whether a given graph has treewidth at most $k$ (and outputs such a decomposition, if it exists). However, due to being fairly complicated to implement, in practice approximation or heuristic methods are often preferred. For these methods one usually relies on finding a chordal completion of the graph. These approaches usually build upon Bouchitté and Todinca's exact algorithm~\cite{BouchitteTodincaAlgorithm, KorhonenBergJarvisaloEvaluatingBouchitteTodinca} by imposing various vertex ordering schemes for computing a chordal completion~\cite{chordalityTestingTarjanYannakakis, RoseTarjanLuekerCompletions, CuthillMcKeeCompletions, berry2004maximum}. 

\subsection{sd-Functors}

We now define functors between sd-categories. Far from routine, pinning down the right notion requires considerable care: we use the data in the sd-category to define permitted structured decompositions, but we need functors not just to preserve the data, but to send decompositions to decompositions in the appropriate sense.

\begin{Def}\label{def sd functors}
Given spined sd-categories $\Gamma = (\cat{C}, \G, \Omega)$ and $\Gamma' = (\mathcal{C}', \G', \Omega')$, an \textbf{sd-functor} $F : \Gamma \to \Gamma'$ consists of
\begin{itemize}
\itemsep -.3em
    \item a functor $F: \mathcal{C} \rightarrow \mathcal{C}'$,
    \item a function $F_g: \mathcal{G} \rightarrow \mathcal{G'}$,
    \item for each $G \in \mathcal{G}$ a functor $F_G : \smallint F_g(G) \rightarrow \smallint G$
\end{itemize}
such that 
\begin{itemize}
    \item $F$ is size-nonincreasing, in the sense that for any object $X \in \cat{C}$, $s'(F(X)) \leq s(X)$, and
    \item given a structured decomposition $d: \smallint G \rightarrow \mathcal{C}$ of an object $X$ in $\mathcal{C}$, the composite
\begin{equation*}
\begin{tikzcd}
	{\smallint F_g(G)} & {\smallint G} & {\cat{C}} & {\cat{C}'}
	\arrow["{F_G}", from=1-1, to=1-2]
	\arrow["d", from=1-2, to=1-3]
	\arrow["F", from=1-3, to=1-4]
\end{tikzcd}
\end{equation*}
which we denote $F_*(d) : \int F_g(G) \to \cat{C}'$ is a structured decomposition of $F(X)$ in $\mathcal{C}'$.
\end{itemize}

\end{Def}

Intuitively, the intermediate functors $F_G$ of Definition~\ref{def sd functors} serve to track how each structured decomposition transforms as we change the structure graph. Their presence introduces some nuance into the definition of composition for sd-functors.

\begin{Def}\label{def composition of sd functors}
Given sd-functors $(\cat{C},\G,\Omega) \xrightarrow{A} (\cat{C}',\G',\Omega') \xrightarrow{B} (\cat{C}'',\G'',\Omega'')$ we define their composite $B\circ A : (\cat{C},\G,\Omega) \rightarrow (\cat{C}'',\G'',\Omega'')$ by the following data:
\begin{itemize}
\itemsep -.3em
    \item the functor $B \circ A : \mathcal{C} \rightarrow \mathcal{C}''$,
    \item the function $(B \circ A)_g = B_g \circ A_g : \G \rightarrow \G''$,
    \item for each $G \in \G$, the functor $(B \circ A)_G = \bigl( A_G \circ B_{A_g(G)}: \smallint (B \circ A)_g(G) \rightarrow \smallint G \bigr)$.
\end{itemize}
\end{Def}

\begin{Lemma} \label{lem composite sd-functors}
Given sd-functors $A$ and $B$ as in Definition \ref{def composition of sd functors}, their composite $B \circ A$ is also an sd-functor.
\end{Lemma}
 
\begin{proof}
Clearly $B \circ A$ is size-nonincreasing since both $A$ and $B$ are. Given a structured decomposition $d: \smallint G \rightarrow \mathcal{C}$ of an object $X$ in $\mathcal{C}$,
we claim that $(B \circ A) \circ d \circ (B \circ A)_G : \smallint (B\circ A)_g G \rightarrow \mathcal{C}''$ is a structured decomposition of $B(A(X))$ in $\mathcal{C}''$. Using the fact that $A$ is an sd-functor, we get that $A \circ d \circ A_G: \smallint A_g(G) \rightarrow \mathcal{C}'$ is a structured decomposition of $A(X)$ in $\mathcal{C}'$. Then, invoking the sd-functor condition for $B$ on the decomposition $A \circ d \circ A_G$, we obtain that $B \circ (A \circ d \circ A_G) \circ B_{A_g(G)}: \smallint B_g (A_g(G)) \rightarrow \mathcal{C}''$ is a structured decomposition of $B(A(X))$
and from there associativity of composition and the definition of $(B \circ A)_G$ lets us conclude that $$B \circ (A \circ d \circ A_G) \circ B_{A_g(G)} = (B \circ A) \circ d \circ (A_G \circ B_{A_g(G)}) = (B \circ A) \circ d \circ (B \circ A)_G$$
from which the claim immediately follows.
\end{proof}

\begin{Rem}
The identity sd-functor $F: (\cat{C}, \G, \Omega) \rightarrow (\cat{C}, \G, \Omega)$ for the composition defined above has the components
    \begin{itemize}
    \itemsep -.3em
        \item $F = \mathrm{id}_\cat{C}: \cat{C} \rightarrow \cat{C}$,
        \item $F_g = \mathrm{id}_\G : \G \rightarrow G$,
        \item for each $G \in \G$ the functor $F_G = \mathrm{id}_{\smallint G} : \smallint G \rightarrow \smallint G$.
    \end{itemize}
Using this along with Lemma \ref{lem composite sd-functors}, we can define the category $\ncat{SdCat}$ whose objects are sd-categories and whose morphisms are sd-functors.
\end{Rem}

\begin{Prop} \label{prop sd-functors give inequality on width}
Given spined sd-categories $\Gamma = (\cat{C}, \G, \Omega)$, $\Gamma' = (\cat{C}', \G', \Omega')$ and an sd-functor $F: \Gamma \to \Gamma'$, then for any $X \in \cat{C}$, we have
\begin{equation*}
    \mathbf{w}_{\Gamma'}(F(X)) \leq \mathbf{w}_\Gamma(X).
\end{equation*}
\end{Prop}

\begin{proof}
Suppose that $d : \smallint G \to \cat{C}$ is a minimal $\G$-structured decomposition of $X$, which implies that $w(d) = \mathbf{w}_\Gamma(X)$. Since $F$ is an sd-functor, $F_*(d) : \int F_g(G) \to \cat{C}'$ is a $\G'$-structured decomposition of $F(X)$. But since $F$ is size-nonincreasing, for every $x \in V(G)$ we have $s'(F(d(x))) \leq s(d(x))$. Hence $w(F_*(d)) \leq w(d) = \mathbf{w}_\Gamma(X)$. Thus $\mathbf{w}_{\Gamma'}(F(X)) \leq \mathbf{w}_\Gamma(X)$.
\end{proof}

The most common use case for sd-functors arises from the desire to relax structural constraints on decompositions (i.e. expand the set of allowed structure graphs) or to refine width measures by adding more permitted bounds.

\begin{Lemma} \label{lem subset of graphs inherits sd-category structure}
If $(\cat{C}, \G', \Omega)$ is a complete width category, width category or spined sd-category, and $\G \subseteq \G'$ is a subset, then $(\cat{C}, \G, \Omega)$ is a complete width category, width category or spined sd-category respectively, and furthermore letting $F = 1_{\cat{C}}$, $F_g : \G \hookrightarrow \G'$ be the inclusion map and for each $G \in \G$, letting $F_G : \int G \to \int G$ be the identity $F_G = 1_{\int G}$ defines an sd-functor $(\cat{C}, \G, \Omega) \to (\cat{C}, \G', \Omega)$.
\end{Lemma}

\begin{proof}
It is clear that $(\cat{C}, \G, \Omega)$ satisfies the conditions (\textbf{W1 - 4}) if $(\cat{C}, \G', \Omega)$ does. Since $\G$ does not affect the size of objects in $\cat{C}$, the identity functor is size-nonincreasing.
\end{proof}

\begin{Rem}
From now on if $\G \subseteq \G'$ and $\Omega = \Omega'$, then we will assume that $F_g : \G \hookrightarrow \G'$ is the inclusion and for each $G \in \G$, $F_G$ is the identity.
\end{Rem}

\begin{Def}
We say that an sd-functor $F: \Gamma \to \Gamma'$ is \textbf{width-preserving} if $\mathbf{w}_\Gamma(X) = \mathbf{w}_{\Gamma'}(F(X))$ for all $X \in \cat{C}$.
\end{Def}

\begin{Lemma} \label{lem sufficient condition for width-preserving sd-functor}
Given sd-categories $\Gamma = (\cat{C}, \G, \Omega)$, $\Gamma' = (\cat{C}', \G', \Omega')$, if $\G \subseteq \G'$, and $F: \Gamma \to \Gamma'$ is an sd-functor that is 
\begin{itemize}
\itemsep -.5em
    \item size-preserving, i.e. for all $X \in \cat{C}$, $s'(F(X)) = s(X)$, 
    \item the inclusion on graph classes, i.e. $F_g : \G \hookrightarrow \G'$, and for all $G \in \G$, $F_G = 1_{\int G}$, and
    \item takes a minimal $\G$-structured decomposition $d$ of an object $X \in \cat{C}$ to a minimal $\G'$-structured decomposition $F_*(d) = F(d)$ of $F(X) \in \cat{C}'$,
\end{itemize} 
then $F$ is width-preserving.
\end{Lemma}

\begin{proof}
If $X \in \cat{C}$ and $d : \smallint G \to \cat{C}$ is a minimal $\G$-structured decomposition of $X$, then $\mathbf{w}_\Gamma(X) = w(d)$. Since $F$ is size nonincreasing and $\G \subseteq \G'$, $F_*(d) = F(d)$ is a $\G'$-structured decomposition of $F(X)$, and furthermore $w(d) = w(F(d))$ because $F$ is size-preserving. Since $F$ preserves minimal decompositions, $F(d)$ is minimal, and therefore $\mathbf{w}_{\Gamma'}(F(X)) = w(F(d)) = w(d) = \mathbf{w}_\Gamma(X)$.
\end{proof}

\subsection{Tree Decompositions} \label{section graph tree decompositions}

The motivation for this work comes from the notion of treewidth in graph theory. To define treewidth, we first need the notion of a tree decomposition. 

\begin{Def}[{\cite{robertsonII}}] \label{def tree decomposition}
Given a graph $G$, an \textbf{RS-tree decomposition}\footnote{Here RS stands for Robertson-Seymour. We note that in \cite{robertsonII}, Robertson and Seymour allow for $G$ to be a multigraph, but here we will restrict to the more usual case of simple graphs.} of $G$ consists of a tree $T$ along with a function $X : V(T) \to P(V(G))$ such that
\begin{enumerate}
    \item $\bigcup_{t \in V(T)} X(t) = V(G)$,
    \item every edge of $G$ has both ends in some $X(t)$, and
    \item given three vertices $t_1, t_2, t_3 \in V(T)$, if $t_2$ lies on the unique path between $t_1$ and $t_3$, then $X(t_1) \cap X(t_3) \subseteq X(t_2)$.
\end{enumerate}
We call the $X(t)$ the \textbf{bags} of the tree decomposition, and the intersections $X(t) \cap X(t')$ the \textbf{adhesions} of the decomposition.
\end{Def}

\begin{Ex} \label{ex pic of tree-decomp}
Below left we have a graph $G$, and on the right an RS-tree decomposition of $G$, where the set labeling each node represents the induced subgraph of $G$ on those vertices. The adhesion between each node is given by the intersection of the subgraphs.
\definecolor{softBlack}{RGB}{45, 47, 49}
\definecolor{creamWhite}{RGB}{245,244,241}
\definecolor{softGray}{RGB}{220, 216, 214}
\definecolor{brick}{RGB}{232, 48, 48}
\begin{figure}[H]
    \centering
    \begin{tikzpicture}[-latex , auto , node distance =1.4 cm and 1.4 cm , on grid , semithick ,                                  state/.style={ circle ,top color =softGray , bottom color = creamWhite ,
                        draw, softBlack , text=softBlack , minimum width =0.2 cm}]
        \node[state] (A) {$a$};
        \node[state] (B) [right=of A] {{$b$}};
        \node[state] (C) [right=of B] {$c$};
        \node[state] (D) [below =of A] {{$d$}};
        \node[state] (F) [below =of C] {$f$};
        \node[state] (G) [below =of D] {$g$};
        \node[state] (H) [right =of G] {$h$};
        \node[state] (I) [below =of F] {$i$};

        \path (A) edge [-] (B);
        \path (B) edge [-] (C);
        \path (B) edge [-] (D);
        \path (D) edge [-] (A);
        \path (F) edge [-] (C);
        \path (G) edge [-] (D);
        \path (I) edge [-] (F);
        \path (G) edge [-] (H);
        \path (I) edge [-] (H);
        \path (H) edge [-] (D);
        \path (H) edge [-] (F);

        \node[state] (lu) [right =2.2cm of C ] {$\{a, b, d\}$};
        \node[state] (ru) [right =2.8cm of lu] {$\{b,c,f\}$};
        \node[state] (mu) [below right =2cm of lu] {$\{b,d,f\}$};
        \node[state] (md) [below =2cm of mu] {$\{d,f,h\}$};
        \node[state] (ld) [below left=2cm of md] {$\{d,g,h\}$};
        \node[state] (rd) [right =2.8cm of ld] {$\{f,h,i\}$};

        \path (lu) edge [-] (mu);
        \path (ru) edge [-] (mu);
        \path (mu) edge [-] (md);
        \path (ld) edge [-] (md);
        \path (rd) edge [-] (md);
    \end{tikzpicture}
    \caption{An example of a tree decomposition of a graph $G$.}
    \label{fig:tree_decomposition_drawing}
\end{figure}
\end{Ex}

Now it is easy to give a definition for the treewidth of a graph $G$.

\begin{Def} \label{def treewidth}
If $(X,T)$ is an RS-tree decomposition of a graph $G$, let 
\begin{equation*}
    w(X,T) = \max_{t \in T} |X(t)| - 1
\end{equation*}
denote the width of the tree decomposition, where $|X(t)|$ is the number of vertices in the bag $X(t)$. Then we define the \textbf{treewidth} $\tw(G)$ of $G$ to be the minimum width $w(X,T)$ over all RS-tree decompositions $(X,T)$ of $G$.
\end{Def}

\begin{Rem}
Note that there is always an RS-tree decomposition $(X,T)$ of a graph $G$ where $T = \bullet$ and $X(\bullet) = V(G)$. Hence $\text{tw}(G) \leq |V(G)| - 1$.
\end{Rem}

\begin{Ex} \label{ex trees have treewidth 1}
Every tree $T$ has $\tw(T) = 1$. To see this, note that one can obtain an RS-tree decomposition where there is one bag per edge of $T$, containing both endpoints of the edge. This decomposition has the smallest possible maximum bag size, and so achieves the treewidth.
\end{Ex}

\begin{Rem} \label{rem RS-tree decomps and sd-tree decomps}
In order to connect our theory to the classical graph theory literature, we must show that the RS-tree decompositions given above can be understood using structured decompositions.

In Appendix \ref{section appendix tree decompositions} we show that to every RS-tree decomposition $(X,T)$ of a graph\footnote{We will actually do this for hypergraphs, which will immediately imply the same result for graphs.} $G$ one can associate a structured decomposition $d : \smallint T \to \ncat{Gr}$ of $G$, and vice versa. These maps will not define a bijection, as the corresponding adhesions of $d$ can change. However, this ultimately makes little difference. In Corollary \ref{cor graph tree width equal to gamma width} we prove that the corresponding width-measure $\mathbf{w}_\Gamma(G)$ defined in terms of tree-shaped structured decompositions of $G$ (Definition \ref{def width}) is equal to $\mathbf{tw}(G)$.

Therefore, in what follows, we will be cavalier and refer to tree-shaped structured decompositions $d : \smallint T \to \ncat{Gr}$ of a graph $G$ simply as \textbf{tree decompositions} of $G$.
\end{Rem}

While Definition \ref{def treewidth} is rather elegant and simple, it turns out that another equivalent definition\footnote{In fact, there are numerous cryptomorphic definitions for treewidth \cite{bertele1972, halin1976s, robertsonII, courcelle2012book}.} for treewidth is very useful for algorithmic purposes.

\begin{Def}
 A \textbf{chord} of an $n$-cycle $\{ x_1, \dots, x_n \}$ is an edge connecting two vertices $x_i$ and $x_j$ such that $|i - j| > 1$. A graph $G$ is \textbf{chordal} if every $n$-cycle in $G$ has a chord for $n > 3$. By a \textbf{chordal completion} of a graph $G$, we mean a monomorphism $f: G \hookrightarrow H$ where $H$ is a chordal graph.
\end{Def}

The chordal graphs can also be characterized using the following theorem of Dirac.

\begin{Th}[{\cite{dirac1961rigid}}]
A graph $H$ is chordal if and only if
\begin{itemize}
\itemsep -.5em
    \item $H$ is a complete graph, i.e. $H \cong K^n$ for some $n \geq 0$, or
    \item there is a monic span $H_0 \hookleftarrow K^n \hookrightarrow H_1$ where $H_0$ and $H_1$ are chordal, and $H \cong H_0 +_{K^n} H_1$. 
\end{itemize}
\end{Th}

The connection between tree decompositions and chordal graphs is given by the following result.

\begin{Lemma}[{\cite[Proposition 12.3.11]{diestel2010graphtheory}}] \label{lem chordal iff has tree-decomp into complete graphs}
A graph $G$ is chordal if and only if it has a tree decomposition $d : \smallint T \to \ncat{Gr}$ where each bag $d(t)$ is a complete graph $K^n$.
\end{Lemma}

\begin{Def}
Given a graph $G$, let $\text{Tree}(G)$ denote the subcategory of $\sd_{\text{Trees}}(\ncat{Gr}, G)$ (Definition \ref{def diagram category}) whose morphisms are monomorphisms. Let $\text{Chord}(G)$ denote the category whose objects are pairs $(f,g)$ where $f : G \hookrightarrow H$ is a monomorphism to some chordal graph $H$ (a \textit{chordal completion} of $G$) and $d$ is a $(\text{Trees}, \{K^n\})$-structured decomposition of $H$. Define the morphisms to be monomorphisms $h : H \to H'$ making the following diagram commute
\begin{equation*}
\begin{tikzcd}[ampersand replacement=\&]
	\& G \\
	H \&\& {H'}
	\arrow["f"', hook, from=1-2, to=2-1]
	\arrow["{f'}", hook, from=1-2, to=2-3]
	\arrow["h", hook, from=2-1, to=2-3]
\end{tikzcd}
\end{equation*}
\end{Def}

Suppose that $(f: G \hookrightarrow H, d)$ is an object of $\text{Chord}(G)$. By Corollary \ref{cor can pullback tree decompositions} the pullback of $d$ along $f$ exists, and is a tree decomposition $f^*(d)$ of $G$. By choosing representatives for each pullback, this extends to a functor
$(-)^* : \text{Chord}(G) \to \text{Tree}(G)$.

Conversely, if $d: \smallint T \to \ncat{Gr}$ is a tree decomposition of $G$, then we can embed each bag $d(t)$, into $K^{|V(d(t))|}$ to form a new tree decomposition $d'$. We call this \textbf{completing} the tree decomposition. The colimit $H$ of $d'$ will then be a chordal graph, and we have a canonical chordal completion $f : G \hookrightarrow H$. This also extends to a functor $\text{Comp}: \text{Tree}(G) \to \text{Chord}(G)$.

\begin{Lemma} \label{lem chordal completion and tree decomp adjunction}
The functor $(-)^*$ is left adjoint to $\text{Comp}$. Furthermore $(-)^* \circ \text{Comp} \cong 1_{\text{Tree}(G)}$ and the functor $\text{Comp}$ is surjective on objects.
\end{Lemma}

It is now obvious how treewidth can be characterized using chordal completions. If $G$ is a graph, let $\omega(G)$ denote its \textbf{clique number}, i.e. the largest integer $n$ such that $K^n$ is a subgraph of $G$.

\begin{Lemma}[{\cite[Corollary 12.3.12]{diestel2010graphtheory}}] \label{lem tree width given by chordal completion}
Given a graph $G$, we have
\begin{equation*}
    \tw(G) = \text{min} \, \{ \omega(H) - 1 \; | \; G \hookrightarrow H, \text{ and $H$ is chordal} \}.
\end{equation*}
\end{Lemma}

\subsection{Nice Tree Decompositions}

The main algorithmic applications of pathwidth, treewidth and related graph parameters proceed by dynamic programming on the associated decompositions. Although the classic characterization of treewidth uses arbitrary tree decompositions, algorithms are more easily described by restricting attention to narrower decomposition structures. For instance, dynamic programming on trees benefits from restricting algorithms to binary trees: trees consisting of vertices with degree at most 3. Our formalism allows us to prove a general result showing that ordinary and binary tree decompositions yield the same width notion. 

\begin{Def}
Let $\mathrm{BTrees}$ denote the class of \textbf{binary trees}, i.e. trees whose vertices have degree at most 3.
\end{Def}

Suppose that $\Gamma = (\cat{C}, \text{Trees}, \Omega)$ is a spined sd-category. Then clearly $\Gamma' = (\cat{C}, \text{BTrees}, \Omega)$ is as well, and by Lemma \ref{lem subset of graphs inherits sd-category structure} we have an sd-functor $\iota: (\cat{C}, \mathrm{BTrees}, \Omega) \rightarrow (\cat{C}, \mathrm{Trees}, \Omega)$ which is the identity on the underlying category. We now show the existence of a width-preserving sd-functor $B: (\cat{C}, \mathrm{Trees}, \Omega) \rightarrow (\cat{C}, \mathrm{BTrees}, \Omega)$.

\begin{Prop}\label{prop binarization functor}
Given a spined sd-category $\Gamma = (\cat{C}, \text{Trees}, \Omega)$ and a $\text{Trees}$-structured decomposition $d : \smallint T \to \cat{C}$ of an object $X \in \cat{C}$, there exists a binary tree $U$ and a functor $U_\ast : \smallint U \rightarrow \smallint T$ such that the composite $(d \circ U_\ast) : \smallint U \rightarrow \mathcal{C}$ is a $\mathrm{BTrees}$-structured decomposition of $X$, and furthermore has the same width as $d$.
\end{Prop}

\begin{proof}
We prove this by simultaneous induction on the maximal degree $m$ of any vertex in $T$, and the number $n$ of vertices which actually have degree $m$. Assume that every tree with $m \leq M$ and $n < N$ satisfies the result. We show that a tree $T$ with maximal degree $M$ and $N$ vertices of maximal degree also satisfies the result. We can further assume that $M > 3$: otherwise $T$ is already a binary tree, so we can take $U = T$ and $U_* = 1_{\int T}$. 

Choose a degree $M$ vertex $v$ of $T$, and partition the adjacent vertices into two sets $L,R$ of respective sizes $\left\lceil \frac{M}{2} \right\rceil$ and $\left\lfloor \frac{M}{2} \right\rfloor$. Let $T^L$ (resp. $T^R$) denote the induced subgraph consisting of those nodes that are connected to $v$ via a path in $L$ ($R$). Since $T$ is a tree, $T^L$ and $T^R$ are disjoint subgraphs of $T$.

Consider the tree $T'$ formed by $T^L$, $T^R$ and 3 new vertices $\ell,v',r$, with edges $e_{\ell, x}$ connecting $\ell$ to each vertex $x$ in $L$, edges $e_{r,x}$ connecting $r$ to all vertices $x$ in $R$, and edges $e_{\ell, v'}$ and $e_{r, v'}$ connecting $v'$ to $\ell$ and $r$ respectively. In $T'$, $\ell$ has degree $|L|+1$, $r$ has degree $|R| + 1$ and $v'$ has degree $2$. Given that $M > 3$, we have $M  > |L| + 1 \geq |R| + 1$, so either the maximal degree of any vertex in $T'$ is less than $M$, or else the number of elements of degree $M$ in $T'$ is at most $N-1$, since $T^L$ and $T^R$ do not contain $v$, and $\ell,v',r$ each have degree strictly less than $M$.

Either way, the inductive hypothesis applies, and we obtain a binary tree $U'$ along with a functor $U'_\ast : \smallint U' \rightarrow \smallint T'$ so that for any $\text{Trees}$-structured decomposition $d': \smallint T' \rightarrow \mathcal{C}$ of an object $X$ in $\mathcal{C}$, the composite $(d' \circ U'_\ast) : \smallint U' \rightarrow \mathcal{C}$ is a $\mathrm{BTrees}$-structured decomposition of $X$ in $\Gamma'$ of the same width as $d'$.

Now let us construct a functor $T'_\ast : \smallint T' \rightarrow \smallint T$ as follows. Define $T'_\ast$ to
\begin{enumerate}
\itemsep -.3em
    \item act as the identity on the vertex and edge objects of $\smallint T^L$ and $\smallint T^R$,
    \item send the vertex objects of $\ell, v', r$ and the edge objects $e_{\ell, v'}, e_{r, v'}$ in $\int T'$ to the vertex object of $v$ in $\int T$,
    \item send the edge morphisms in $\int T'$ connecting $e_{\ell, v'}$ and $e_{r, v'}$ to the vertex object of $v'$, to the identity morphism on the vertex object of $v$,
    \item the edge morphisms connecting the edge objects $e_{\ell, x}$ between $\ell$ and each vertex $x$ of $L$ in $T'$ to the respective edge morphisms between $v$ and $L$ in $T$, and
    \item similarly for the edge morphisms between $r$ and the vertices of $R$ in $T'$.
\end{enumerate}

Now the diagram $d \circ T'_\ast : \smallint T' \rightarrow \mathcal{C}$ differs from $d$ only by the insertion of a trivial span, so any cocone for $d$ extends uniquely to a cocone for $(d \circ T'_\ast)$, and it clearly has the same width as $d$.

Now setting $U = U'$ and $U_\ast = T'_\ast \circ U'_\ast : \smallint U' \rightarrow \smallint T$ gives for any $\text{Trees}$-structured decomposition $d: \smallint T \rightarrow \mathcal{C}$ of an object $X$ in $\mathcal{C}$, the composite $d \circ U_\ast : \smallint U \rightarrow \mathcal{C}$ is a $\mathrm{BTrees}$-structured decomposition of $X$ with the same width as $d$.

Thus by setting $B : \cat{C} \to \cat{C}$ to be the identity, $B_g : \text{Trees} \to \text{BTrees}$ to be the map that assigns to a tree $T$ the binary tree $U'$ as constructed above, and for $T \in \text{Trees}$ we let $B_T : \int B_g(T) \to \int T$  be $U_*$, we obtain a width-preserving sd-functor $B : (\cat{C}, \text{Trees}, \Omega) \to (\cat{C}, \text{BTrees}, \Omega)$.
\end{proof}

Now consider a graph class $\G$ which satisfies closure under barycentric subdivision, so that whenever the class $\G$ contains some graph $G$, we have that its barycentric subdivision $\smallint G$, regarded as a graph, also belongs to $\G$. For example $\G = \text{Trees}$ satisfies this property. 

Let $\smallint \G$ denote the class of graphs in $\G$ that are in the image of the barycentric subdivision. For the remainder of this section, let $\iota : (\cat{C}, \smallint \G, \Omega) \to (\cat{C}, \G, \Omega)$ denote the identity on $\cat{C}$, which is an sd-functor by Lemma \ref{lem subset of graphs inherits sd-category structure}. By constructing an appropriate sd-functor from $(\cat{C}, \G, \Omega)$ to $(\cat{C}, \smallint \G, \Omega)$, we show that restricting to $(\cat{C}, \smallint \G, \Omega)$ does not change the associated $\Gamma$-width.

\begin{Def}
Define the \textbf{barycentric restriction functor} $S$ between $(\cat{C}, \G, \Omega)$ and $(\cat{C}, \smallint \G, \Omega)$ by setting the underlying functor $S: \cat{C} \rightarrow \cat{C}$ to be the identity, $S_g : \G \rightarrow \smallint G$ as barycentric subdivision so that $S_g(G) = \smallint G$ regarded as a graph, and for $G \in \mathcal{G}$ define $S_G : \smallint S_g G \rightarrow \smallint G$ in the following manner. 

In $\smallint S_g G$, objects come in two types: those corresponding to vertices of $S_gG$ (originating from vertex or edge objects of $\smallint G$) and those arising as the vertex resulting from the subdivision of an edge in $S_gG$. The functor $S_G$ sends each vertex object of the former kind to the corresponding vertex object in $\smallint G$, and both the edge object and the subdivision vertex object corresponding to each edge $e$ in $S_g G$ to the edge object $e$ in $\smallint G$, while sending the legs of every span connecting preimages of edge objects to the identity morphism. These functors $S_G$ are clearly final, so the sd-functor condition pertaining to preservation of colimits is satisfied.
\end{Def}

\begin{Ex}
Consider the graph $G$ consisting of a single edge $e$ and its end vertices $x,y$. The upper diagram depicts the category $\smallint G$, the diagram lower depicts the category $\smallint S_g G$, where now $e$ has become a vertex object, and there are two new edge objects $ex$ and $ey$. The dashed arrows depict the action of the functor $S_G: \smallint S_g G \rightarrow \smallint G$ on the objects of $\smallint S_g G$; the morphisms of $\smallint S_g G$ are labeled with their $S_G$-images.
\begin{equation*}
\begin{tikzcd}
	& x & e & y \\
	\\
	x & ex & e & ey & y
	\arrow[from=1-3, to=1-2]
	\arrow[from=1-3, to=1-4]
	\arrow[dashed, maps to, from=3-1, to=1-2]
	\arrow[dashed, maps to, from=3-2, to=1-3]
	\arrow[from=3-2, to=3-1]
	\arrow[from=3-2, to=3-3]
	\arrow[dashed, maps to, from=3-3, to=1-3]
	\arrow[dashed, maps to, from=3-4, to=1-3]
	\arrow[from=3-4, to=3-3]
	\arrow[from=3-4, to=3-5]
	\arrow[dashed, maps to, from=3-5, to=1-4]
\end{tikzcd}
\end{equation*}
\end{Ex}

\begin{Prop}\label{prop barycentric restriction preserves width}
The barycentric subdivision functor defines a width-preserving sd-functor $S : (\cat{C}, \G, \Omega) \to (\cat{C}, \smallint \G, \Omega)$.
\end{Prop}

\begin{proof}
Clearly $S$ is size non-increasing and the barycentric restriction functor associates to each $G$-structured decomposition of $X$ to an $\smallint G$-structured decomposition of $X$.
\end{proof}

Combining Proposition \ref{prop binarization functor} with Proposition \ref{prop barycentric restriction preserves width} gives the existence of so-called \textit{nice tree decompositions} in any sd-category of tree-structured decompositions. Nice tree decompositions of graphs are commonly used when defining dynamic programming algorithms on tree decompositions, see \cite[Chapter 10]{niedermeier2006invitation} for more information.

\begin{Cor}\label{cor nice decompositions}
Given an sd-category of the form $\Gamma = (\cat{C}, \text{Trees}, \Omega)$ and $X \in \cat{C}$, we can find a binary tree $G$ and a $G$-structured decomposition $d : \smallint G \rightarrow \cat{C}$ of $X$ such that each span $d(v_1) \hookleftarrow d(e) \hookrightarrow d(v_1)$ has one of the following forms:
\begin{enumerate}
    \item $d(v_1) \cong d(e) \cong d(v_2)$ and both monos of the span are isomorphisms (\textit{join nodes}),
    \item $d(v_1) \cong d(e)$, the left mono is an isomorphism but the right mono is not, or
    \item $d(e) \cong d(v_2)$, the right mono is an isomorphism and the left mono is not (\textit{introduce/forget nodes}).
\end{enumerate}
\end{Cor}

A similar, but simpler argument yields the existence of \textit{nice path decompositions} as well.

\section{Examples} \label{section examples}

In this section we detail several examples of our formalism that already exist in the literature. We are able to capture treewidth, pathwidth, complemented treewidth, tree independence number and hypergraph treewidth directly as examples of $\Gamma$-width, and we draw connections between sd-categories and $\Gamma$-width with Carmesin's graph decompositions, layered treewidth and $\mathcal{H}$-treewidth.

\subsection{Treewidth}

We obtain the main combinatorial invariant of a graph of interest in this paper, namely treewidth, as an example of $\Gamma$-width as follows.

\begin{Prop} \label{prop sd-category for treewidth}
The triple $\Gamma = (\ncat{Gr}, \text{Trees}, \{ K^n \}_{n \geq 0})$ is a complete width category and the $\Gamma$-width of a graph $G \in \Gr$ is precisely its treewidth $\tw(G) = \mathbf{w}_\Gamma(G)$.
\end{Prop}

\begin{proof}
Condition (\textbf{W1}) holds by Corollary \ref{cor colimits of tree decomps exist in Gr}, (\textbf{W2}) holds since there exist no monomorphisms $f : G \to H$ of graphs where $|V(G)| > |V(H)|$, and (\textbf{W3}) holds by Corollary \ref{cor can pullback tree decompositions}. Furthermore (\textbf{W4}) holds because completing a tree decomposition $d$ of a graph $G$ to a $(\text{Trees}, \{K^n \}_{n \geq 0})$-decomposition $d'$ can never reduce the number of vertices of the corresponding colimit $\colim \, d' = H$. This is because the bags of $d'$ can only be greater than or equal to the size of the bags of $d$, and the adhesions are the same. In effect we are quotienting larger graphs by the ``same relation.'' Hence the induced map $G \to H$ must be a monomorphism. The fact that $\mathbf{w}_\Gamma(G) = \tw(G)$ is Corollary \ref{cor graph tree width equal to gamma width}.
\end{proof}

In fact, we need not restrict ourselves to $\Gr$ to capture treewidth. Let $K^n_r$ denote the reflexive $n$-clique, i.e. the $n$-clique additionally having a single loop on each vertex. Let $K^n_\mlt$ denote the multi $n$-clique. This is the multigraph with $n$ vertices and between every pair of (not necessarily distinct) vertices there are $n$ edges.

\begin{Lemma} \label{lem different treewidths}
The triples $\Gamma = (\Gr_\ell, \text{Trees}, \{K^n_r \}_{n \geq 0})$ and $\Gamma' = (\Gr_\mlt, \text{Trees}, \{K^n_\mlt \}_{n \geq 0})$ both form width categories, and the inclusion functors
$\Gr \hookrightarrow \Gr_\ell \hookrightarrow \Gr_\mlt$ are sd-functors. Furthermore if $G$ is a simple graph, then 
\begin{equation*}
    \mathbf{w}_\Gamma(G) = \mathbf{w}_{\Gamma'}(G) = \tw(G).
\end{equation*}
\end{Lemma}

\begin{proof}
This follows from noticing that the inclusion functors above all preserve finite limits and pushouts of spans of monomorphisms. The size of a loop graph $G$ in $(\Gr_\ell, \text{Trees}, \{K^n \} \cup \{K^n_r\}_{n \geq 0})$ is its number of vertices, while the size of a multigraph $H$ in $(\Gr_\mlt, \text{Trees}, \{K^n \} \cup \{ K^n_\mlt \}_{n \geq 0})$ is the maximum of its vertices and edges. For a simple graph $G$, its number of vertices is always greater than or equal to its number of edges. Hence the size of any simple graph embedded in $\Gr_\mlt$ is still its number of vertices. Thus the inclusion functors clearly are size preserving. Now for a simple graph $G$ embedded in $\Gr_\ell$ or $\Gr_\mlt$, none of the multigraphs appearing in any of the bags or adhesions of any decomposition of $G$ can have loops or multiple edges. Hence every decomposition must factor through $\Gr$, and therefore the inclusion functors must take minimal decompositions of graphs to minimal decompositions, thus by Lemma \ref{lem sufficient condition for width-preserving sd-functor}, the inclusion functors are width-preserving.
\end{proof}

\subsection{Pathwidth}
Pathwidth is another common width-measure used for graphs \cite[Page 279]{diestel2010graphtheory}, \cite{robertsonI}.

\begin{Def}
A path-decomposition of a graph $G$ is a $\text{Paths}$-structured decomposition (Definition \ref{def complete graphs, cycles, trees and paths}) $d : \smallint P \to \Gr$. The corresponding notion of \textbf{pathwidth} $\textbf{pw}(G)$ of a graph $G$ is defined in the obvious way.
\end{Def}

Clearly every path is a tree, i.e. $\text{Paths} \subseteq \text{Trees}$.

\begin{Lemma}
The triple $\Gamma = (\Gr, \text{Paths}, \{K^n\}_{n \geq 0})$ forms a complete width category. Furthermore the $\Gamma$-width of a graph $G$ is precisely its pathwidth $\mathbf{pw}(G)$.
\end{Lemma}

\begin{proof}
This follows from Proposition \ref{prop sd-category for treewidth} and Lemma \ref{lem subset of graphs inherits sd-category structure}.
\end{proof}

The following result is obvious and well-known, but we include it anyway, as it follows immediately from our formalism.

\begin{Cor} \label{cor pathwidth greater than treewidth}
Given a graph $G$, we have $\textbf{tw}(G) \leq \textbf{pw}(G)$.
\end{Cor}

\begin{proof}
The identity functor $(\ncat{Gr}, \text{Paths}, \{ K^n \}_{n \geq 0}) \to (\ncat{Gr}, \text{Trees}, \{ K^n \}_{n \geq 0})$ is a sd-functor, so the result holds by Proposition \ref{prop sd-functors give inequality on width}.
\end{proof}

\subsection{Complemented Treewidth}

Another graph width measure of interest is the notion of \textbf{complemented treewidth}. Given a graph $G$, this is defined as
\begin{equation*}
   \overline{\tw}(G) := \tw(\overline{G}),
\end{equation*}
where $\overline{G}$ denotes the complement\footnote{Recall the complement of a graph $G$ is the graph $\overline{G}$ where $xy$ is an edge in $\overline{G}$ if and only $x$ and $y$ are not incident in $G$.} of $G$. Complemented treewidth has useful applications in algorithmics~\cite{sousa2021} since it allows for recursive algorithms on graph classes that are both dense and incomparable to bounded treewidth classes. Using the theory of spined categories \cite[Section 5.2]{bumpus2023spined}, one can already provide a category-theoretic characterization of complemented treewidth. Here we recast this characterization in the language of width categories.

\begin{Def} \label{def reflexive graph homomorphism}
Let $\overline{\Gr}$ denote the category whose objects are graphs, and whose morphisms $f : G \to H$ are vertex maps $f : V(G) \to V(H)$ with the property that if $f(x)f(y)$ is an edge in $H$, then $xy$ is an edge in $X$. We call this a \textbf{reflexive graph homomorphism}.
\end{Def}

Let $\overline{(-)} : \Gr \to \overline{\Gr}$ denote the \textbf{complementation functor} that sends a graph $G$ to its complement $\overline{G}$ and a graph homomorphism $f : G \to H$ to the corresponding reflexive graph homomorphism $\overline{f} : \overline{G} \to \overline{H}$.

\begin{Lemma} \label{lem monos/epis in complemented graphs}
A reflexive graph homomorphism $f : G \to H$ is a mono/epimorphism in $\overline{\Gr}$ if and only $V(f)$ is an in/surjective function.
\end{Lemma}

\begin{proof}
The vertex set functor $V : \overline{\Gr} \to \ncat{Set}$ is faithful, and furthermore has a left adjoint $\text{CoDisc}$ and a right adjoint $\text{Disc}$, so the result follows from the same reasoning as in the beginning of Section \ref{section hypergraphs}.
\end{proof}

\begin{Lemma}[{\cite[Proposition 5.10]{bumpus2023spined}}] \label{lem complementation functor is iso}
The complementation functor $\overline{(-)} : \Gr \to \overline{\Gr}$ is an isomorphism of categories.
\end{Lemma}

We abuse notation and let $\overline{(-)} : \overline{\Gr} \to \Gr$ denote the inverse of the isomorphism of Lemma \ref{lem complementation functor is iso}. If we let $\text{Disc}(n)$ denote the discrete graph on $n$ vertices, then it is obvious that $\overline{K^n} = \text{Disc}(n)$. Furthermore, since $\overline{(-)}$ is an isomorphism, we know that we can transfer the width category structure of $\Gr$ over to $\overline{\Gr}$.

\begin{Prop} \label{prop sd-category for co-treewidth}
Graphs with reflexive graph homomorphism form a width category $\Gamma = (\overline{\Gr}, \text{Trees}, \{ \text{Disc}(n)\}_{n \geq 0})$. Furthermore the $\Gamma$-width of a graph $G$ is precisely its complemented treewidth $\overline{\tw}(G)$.
\end{Prop}

\begin{proof}
Clearly $\Gamma$ satisfies (\textbf{W1}) and (\textbf{W3}) since it is isomorphic to the width category of Proposition \ref{prop sd-category for treewidth}. Since the size in $\Gamma$ is the number of vertices, by Lemma \ref{lem monos/epis in complemented graphs}, condition (\textbf{W2}) is satisfied. Now a functor $d : \smallint T \to \overline{\Gr}$ is a decomposition of a graph $G$ if and only if $\overline{d} = \overline{(-)} \circ d : \smallint T \to \Gr$ is a tree decomposition of $\overline{G}$, and $w(d) = w(\overline{d})$. Therefore $\mathbf{w}_\Gamma(G) = \tw(\overline{G}) = \overline{\tw}(G)$.
\end{proof}

\subsection{Tree Independence Number}

We can also capture another interesting graph invariant, called the \textbf{tree independence number} $\tw_\alpha(G)$ of a graph $G$. This invariant, introduced in \cite{dallard2024treeindependence}, is defined similarly to treewidth, but using a different measure of size for the bags. It is important in algorithmics as the \textsc{Maximum Weight Independent Packing} problem can be solved in polynomial time on a graph $G$ equipped with a tree decomposition of bounded independence number \cite[Theorem 7.2]{dallard2024treeindependence}.

\begin{Def}
Given a graph $G$, an \textbf{independent set} $I$ of $G$ is a subset $I \subseteq V(G)$ such that each pair $x,y \in I$ are not adjacent, i.e. $xy \notin E(G)$. The \textbf{independence number} $\alpha(G)$ of a graph $G$ is the cardinality of the largest independent set of $G$. The tree independence number $\tw_\alpha(d)$ of a tree decomposition $d : \smallint T \to \ncat{Gr}$ of a graph $G$ is $\max_{t \in V(T)} \alpha(d(t)) - 1$. The \textbf{tree independence number} $\tw_\alpha(G)$ of a graph $G$ is the minimal tree independence number of all of its tree decompositions.
\end{Def}

We wish to capture tree independence number using a width category. If we try to do this using $\Gr$ as the underlying category, we immediately run up against the fact that independence number is not monotone with respect to monomorphisms.

\begin{Ex}
Take $G = \text{Disc}(n)$ to be the discrete graph on $n > 1$ vertices and $H = K^n$ the complete graph on $n$ vertices. Then there is a canonical monomorphism $G \hookrightarrow H$ in $\Gr$, but $\alpha(G) = n$ and $\alpha(H) = 1$.
\end{Ex}

By Lemma \ref{lem size monotonic for spined sd-categories}, this implies that there does not exist a spine on $\Gr$ with which we can obtain tree independence number as $\Gamma$-width. However, the following result provides us with the correct category with which to capture tree independence number. Recall the notion of a reflexive graph homomorphism from Definition \ref{def reflexive graph homomorphism}.

\begin{Lemma} \label{lem reflexive homomorphisms push forward independence sets}
Let $f : G \to H$ be a reflexive graph homomorphism. If $I \subseteq V(G)$ is an independent set of $G$, then $f(I)$ is an independent set of $H$.
\end{Lemma}

\begin{proof}
Suppose that $x,y \in I$ and $f(x)f(y) \in E(H)$. Then since $f$ is reflexive, this implies $xy \in E(G)$, which is a contradiction as $I$ is an independent set. Thus $f(I)$ is an independent set of $H$.
\end{proof}

\begin{Lemma} \label{lem reflexive graph homomorphisms monotone wrt independence number}
If $ f: G \hookrightarrow H$ is a monomorphism in $\overline{\Gr}$, then $\alpha(G) \leq \alpha(H)$.
\end{Lemma}

\begin{proof}
This follows from Lemma \ref{lem monos/epis in complemented graphs} and Lemma \ref{lem reflexive homomorphisms push forward independence sets}.
\end{proof}

Let $\Omega^\alpha_k$ denote the set of graphs $G$ with $\alpha(G) \leq k$. Let $\Omega^\alpha$ denote the corresponding filtered full subcategory of $\overline{\Gr}$.

\begin{Prop} \label{prop sd-category for tree independence number}
Graphs with reflexive graph homomorphisms form a width category $\Gamma = (\overline{\Gr}, \text{Trees}, \Omega^\alpha)$. The $\Gamma$-width of a graph $G$ is precisely its tree independence number $\mathbf{tw}_\alpha(G)$.
\end{Prop}

\begin{proof}
The conditions (\textbf{W1}) and (\textbf{W3}) follow from the same reasoning as in Proposition \ref{prop sd-category for co-treewidth}. Condition (\textbf{W2}) holds by Lemma \ref{lem reflexive graph homomorphisms monotone wrt independence number}. The fact that $\mathbf{w}_\Gamma(G) = \mathbf{tw}_\alpha(G)$ is practically by definition.
\end{proof}

\begin{Cor}
Given a reflexive graph homomorphism $ f: G \to H$ which is furthermore injective on vertices we have
\begin{equation*}
    \mathbf{tw}_\alpha(G) \leq \mathbf{tw}_\alpha(H).
\end{equation*}
\end{Cor}

\begin{proof}
This follows from Proposition \ref{prop sd-category for tree independence number} and Proposition \ref{prop Gamma width is monotone with respect to monomorphisms}.
\end{proof}

Let $\overline{\Gr}_\alpha$ denote the width category from Proposition \ref{prop sd-category for tree independence number}, and let $\overline{\Gr}_\tw$ denote the width category from Proposition \ref{prop sd-category for co-treewidth}. 

\begin{Lemma}
The identity functor on $\overline{\Gr}$ is a sd-functor $1_{\overline{\Gr}} : \overline{\Gr}_\tw \to \overline{\Gr}_\alpha$.
\end{Lemma}

\begin{proof}
If $G$ is a graph, then clearly $\alpha(G) \leq |V(G)|$. Thus the result follows by Proposition \ref{prop sd-functors give inequality on width}.
\end{proof}

\begin{Cor} \label{cor tree independence number and co-treewidth}
Given a graph $G$,
\begin{equation*}
    \tw_\alpha(G) \leq \overline{\tw}(G).
\end{equation*}
\end{Cor}

There is a well-known duality between independent sets and cliques.

\begin{Lemma} \label{lem duality between cliques and independent sets}
Given a graph $G$, a subgraph $K \subseteq G$ is complete if and only if the corresponding subgraph of its complement $\overline{K} \subseteq \overline{G}$ is an independent set.
\end{Lemma}

Recall the notion of the clique number $\omega(G)$ of a graph $G$ from Lemma \ref{lem tree width given by chordal completion}. Let $\Omega^\omega_k$ denote the set of graphs with clique number less than or equal to $k \geq 0$, and let $\Omega^\omega$ denote the corresponding full subcategory of $\Gr$.

\begin{Lemma} \label{lem clique number monotone under monomorphisms}
If $f : G \hookrightarrow H$ is a monomorphism in $\Gr$, then $\omega(G) \leq \omega(H)$.
\end{Lemma}

\begin{proof}
If $K \subseteq G$ is a complete subgraph of $G$, then $f(K)$ is a complete subgraph of $H$ with the same number of vertices.
\end{proof}

\begin{Prop} \label{prop sd-category for tree clique number}
Graphs with clique number form a width category $\Gamma = (\Gr, \text{Trees}, \Omega^\omega)$.
\end{Prop}

\begin{proof}
Conditions (\textbf{W1}) and (\textbf{W3}) follow from the same reasoning as in Proposition \ref{prop sd-category for treewidth}. Condition (\textbf{W2}) follows from Lemma \ref{lem clique number monotone under monomorphisms}.
\end{proof}

We now define the following auxiliary graph parameter, which is intimately related with tree-independence number by Lemma \ref{lem sd-functor between tcn and tw alpha}.

\begin{Def}
Given a graph $G$, we call the $\Gamma$-width of $G$ with respect to the width category $\Gamma$ of Proposition \ref{prop sd-category for tree clique number} the \textbf{tree clique number} of $G$ and denote it 
\begin{equation*}
    \mathbf{tcn}(G) = \mathbf{w}_\Gamma(G).
\end{equation*}
\end{Def}

Let $\Gr_\omega$ denote the width category $\Gr_\omega = (\Gr, \text{Trees}, \Omega^\omega)$ and let $\overline{\Gr}_\alpha$ denote the width category $\overline{\Gr}_\alpha = (\overline{\Gr}, \text{Trees}, \Omega^\alpha)$. 

\begin{Lemma} \label{lem sd-functor between tcn and tw alpha}
The complementation functors $\overline{(-)}: \Gr_\omega \to \overline{\Gr}_\alpha$ and $\overline{(-)}: \overline{\Gr}_\alpha \to \Gr_\omega$ are size-preserving sd-functors.
\end{Lemma}

\begin{proof}
Since $\overline{(-)}$ is an isomorphism of the underlying categories by Lemma \ref{lem complementation functor is iso}, we need only to show that $\overline{(-)}$ is size nonincreasing. Let $G$ be a graph. Its size in $\Gr_\omega$ is equal to its clique number $\omega(G)$. Suppose that $K^n \subseteq G$ is a maximal clique, so that $\omega(G) = n$. Then $\overline{K^n} \subseteq \overline{G}$ is an $n$-independent set. If $I \subseteq V(\overline{G}) = V(G)$ is a $k$-independent set of $\overline{G}$ with $n \leq k$, then $\overline{I}$ is a $k$-clique in $G$ by Lemma \ref{lem duality between cliques and independent sets}. Hence $n = k$, and so $I$ is maximal. Thus $\omega(G) = \alpha(\overline{G})$, which immediately implies that $\omega(\overline{G}) = \alpha(G)$. 
\end{proof}

\begin{Cor} \label{cor tree independence number and vertex cover number}
Given a graph $G$ we have
\begin{equation*}
    \mathbf{tcn}(G) = \mathbf{tw}_\alpha(\overline{G}).
\end{equation*}
\end{Cor}

\begin{proof}
By Lemma \ref{lem sd-functor between tcn and tw alpha} and Proposition \ref{prop sd-functors give inequality on width}, we have
\begin{equation*}
    \mathbf{tw}_\alpha(\overline{G}) \leq \mathbf{tcn}(G)
\end{equation*}
and
\begin{equation*}
    \mathbf{tcn}(\overline{G}) \leq \mathbf{tw}_\alpha(G).
\end{equation*}
Hence $\mathbf{tcn}(G) \leq \mathbf{tw}_\alpha(\overline{G}) \leq \mathbf{tcn}(G)$.
\end{proof}

\subsection{Hypergraph Treewidth} \label{section hypergraph treewidth}

In this section, we construct sd-categories which capture the notions of hypergraph treewidth and generalized hypertreewidth.

\begin{Def} \label{def hypergraph}
A \textbf{hypergraph} $H$ consists of a finite set $V(H)$ of vertices and a subset $E(H) \subseteq P_{\neq \varnothing}(V(H))$ of nonempty subsets of $V(H)$ called hyperedges. If $e \in E(H)$ is a hyperedge with cardinality $n$, we call $e$ an $n$-edge, and let $P_n(V(H))$ denote the set of $n$-subsets of $V(H)$ and $E_n(H)$ the set of $n$-edges of $H$. A map $f : H \to H'$ of hypergraphs consists of a function $V(f) : V(H) \to V(H')$ such that if $e \in E_n(H)$, then $f(e) \in E_n(H')$. We call $f$ a \textbf{hypergraph homomorphism}. Let $\Hyp$ denote the category of hypergraphs and hypergraph homomorphisms.
\end{Def}

Definition \ref{def hypergraph} is equivalent to the definition of hypergraphs given in \cite[Section 2]{adler2007hypertree}, except that they do not allow for empty hypergraphs or hypergraphs with isolated vertices, namely those vertices that do not belong to any hyperedge. 

We can consider a graph $G$ as a special kind of hypergraph. Namely it is a hypergraph where each hyperedge $e \in E(G)$ has cardinality $2$. Then a function $V(f) : V(G) \to V(G')$ is a map of graphs if and only if it is a map of hypergraphs. In other words, we obtain a fully faithful inclusion functor $\iota : \Gr \hookrightarrow \Hyp$.

\begin{Rem} \label{rem multihypergraphs}
We mention another popular category of hypergraphs from the literature. Let $\Hyp_{\mlt}$ denote the category considered in \cite[Section 2]{Dorfler1980categoryhypergraphs}. An object $H$ of this category consists of finite sets $V(H)$, $E(H)$, and a function $\partial : E(H) \to P_{\neq \varnothing}(V(H))$. This means that multiple edges can contain the same vertices, and hence we call these \textbf{multihypergraphs}. If we require $\partial$ to be injective, then we get precisely our notion of hypergraph, which is sometimes called a \textbf{simple hypergraph}. We note that the corresponding notion of morphism in \cite{Dorfler1980categoryhypergraphs} is quite loose. It allows for example the unique morphism $f : H \to H'$ where $H$ has three vertices $x,y,z$ and one hyperedge $e = \{x, y, z\}$ and $H'$ has one vertex $w$ and one hyperedge $\ell = \{w \}$. We will return to this observation in Remark \ref{rem can't use Hyp mlt for gaifman}.
\end{Rem}

\begin{Def} \label{def gaifman graph}
Given a hypergraph $H$, we let $\Gaif(H)$ denote the \textbf{Gaifman graph}\footnote{Also called the \textbf{primal graph} or \textbf{clique graph} of $H$.} of $H$. This is the graph with $V(\Gaif(H)) = V(H)$ and where there is an edge $xy \in E(\Gaif(H))$ if and only if $x \neq y$ and there exists a hyperedge $e \in E(H)$ with $x,y \in e$.   
\end{Def}

The Gaifman graph construction extends to a functor $\Gaif : \Hyp \to \Gr$. Indeed if $f : H \to H'$ is a map of hypergraphs, then $V(\Gaif(f)) : V(\Gaif(H)) \to V(\Gaif(H'))$ is equal to $V(f) : V(H) \to V(H')$, so we need only to show that $\Gaif(f)$ is a graph homomorphism. Suppose that $xy \in E(\Gaif(H))$, so $x \neq y$. Then there exists a $n$-edge $e \in E_n(H)$ of cardinality $n > 1$ with $x,y \in e$. Since $f$ is a map of hypergraphs, this means that $f(e)$ is an $n$-edge, and $f(x) \neq f(y)$. Thus $f(x)f(y) \in \Gaif(H')$. Thus $\Gaif(f)$ is a graph homomorphism.

\begin{Rem} \label{rem can't use Hyp mlt for gaifman}
Had we chosen to work instead with $\Hyp_\mlt$, then we would not obtain a functor to $\Gr$, but instead to $\Gr_\ell$, as we would have to be able to collapse edges to loops, such as in Remark \ref{rem multihypergraphs}.
\end{Rem}

As with graphs in Section \ref{section graph tree decompositions}, there is an obvious extension of the notion of Robertson-Seymour tree decompositions for hypergraphs.

\begin{Def}[{\cite[Page 11]{heinz2013tree}}] \label{def hypergraph treewidth}
An RS-tree decomposition $(X, T)$ of a hypergraph $H$ consists of a tree $T$ along with a function $X : V(T) \to P(V(H))$ such that
\begin{enumerate}
    \item $\bigcup_{t \in V(T)} \, X(t) = V(G)$,
    \item for every $e \in E(H)$, there exists a bag $X(t)$ such that $e \subseteq X(t)$,
    \item given three vertices $t_1, t_2, t_3 \in V(T)$, if $t_2$ lies on the unique path between $t_1$ and $t_3$, then $X(t_1) \cap X(t_3) \subseteq X(t_2)$.
\end{enumerate}
The width $w(X,T)$ of a tree decomposition $(X,T)$ of $H$ is the maximal number of vertices of its bags minus $1$. The hypergraph treewidth $\tw_\Hyp(H)$ of $H$ is defined to be the minimal width of all of its tree decompositions.
\end{Def}

\begin{Rem}
We note that if $H$ is a graph, then Definition \ref{def hypergraph treewidth} is precisely the same as Definition \ref{def tree decomposition}.
\end{Rem}

In order to connect hypergraph treewidth with the $\Gamma$-width of Proposition \ref{prop hypergraph treewidth equal to gamma width}, we need to connect tree-shaped structured decompositions of hypergraphs with RS-tree decompositions. This takes a bit of work, so the construction and proofs are relegated to Section \ref{section appendix tree decompositions}. We therefore abuse notation and call both RS-tree decompositions and tree-shaped structured decomposition $d : \smallint T \to \Hyp$ of a hypergraph $H$ tree decompositions.

Let $K^n_\Hyp$ denote the hypergraph with $V(K^n_\Hyp) = \{1, \dots, n \}$ and where every nonempty subset of $\{1, \dots, n \}$ is a hyperedge. We call $K^n_\Hyp$ the \textbf{complete hypergraph} on $n$ vertices. A hypergraph $H$ has a monomorphism $H \hookrightarrow K^n_\Hyp$ if and only if $|V(H)| \leq n$.

\begin{Prop} \label{prop sd-category for hypergraph treewidth}
Hypergraphs form a width category $\Gamma = (\Hyp, \text{Trees}, \{K^n_\Hyp\}_{n \geq 0} )$. Furthermore, given a hypergraph $H$, its hypergraph treewidth and $\Gamma$-width are equal
\begin{equation*}
\mathbf{tw}_\Hyp(H) = \mathbf{w}_\Gamma(H).
\end{equation*}
\end{Prop}

\begin{proof}
It is clear that (\textbf{W1}) holds, (\textbf{W2}) holds by Lemma \ref{lem pushouts of monos in Hyp} and (\textbf{W3}) holds by Lemma \ref{lem hypergraphs have pullback stable pushouts of monos} and a similar argument for (\textbf{W4}) as in the proof of Proposition \ref{prop sd-category for treewidth}. The fact that $\mathbf{tw}_\Hyp(H) = \mathbf{w}_\Gamma(H)$ is Proposition \ref{prop hypergraph treewidth equal to gamma width}.
\end{proof}

It is well-known that the hypergraph treewidth of a hypergraph $H$ is equal to the treewidth of its Gaifman graph $\Gaif(H)$ \cite[Page 2]{gottlob2005hypertreedecomp}. We will now prove this using the sd-category framework.

For the next result, let $\Gr$ denote the sd-category from Proposition \ref{prop sd-category for treewidth}, with $\mathbf{w}_\Gr(G)$ its corresponding $\Gamma$-width, and let $\Hyp$ denote the sd-category of hypergraphs given above, with $\mathbf{w}_\Hyp(H)$ its corresponding $\Gamma$-width.

\begin{Lemma}
The inclusion functor $\iota : \Gr \to \Hyp$ and the Gaifman graph functor $\Gaif : \Hyp \to \Gr$ are both sd-functors.
\end{Lemma}

\begin{proof}
The inclusion functor preserves pushouts of spans of monomorphisms by construction, and is size-preserving. The Gaifman graph functor preserves pushouts of spans of monomorphisms by Lemma \ref{lem gaifman graph functor preserves pushouts of monos} and is also size-preserving.
\end{proof}

\begin{Cor}
Given a hypergraph $H$ and a graph $G$,
\begin{equation*}
    \mathbf{w}_\Hyp(\iota(G)) \leq \tw(G),
\end{equation*}
and
\begin{equation*}
\tw(\Gaif(H)) \leq \mathbf{w}_\Hyp(H).
\end{equation*}
\end{Cor}

Thus we see that if $G$ is a graph, then $\mathbf{w}_\Hyp(\iota(G)) = \tw(G)$. Let us now prove the reverse of the second inequality.

\begin{Lemma} \label{lem hypergraph treewidth less than gaifman treewidth}
If $H$ is a hypergraph, then
\begin{equation*}
    \mathbf{w}_\Hyp(H) \leq \tw(\Gaif(H)).
\end{equation*}
\end{Lemma}

\begin{proof}
Suppose that $d : \smallint T \to \Gr$ is a tree decomposition of $\Gaif(H)$. Then by Proposition \ref{prop hypergraph treewidth equal to gamma width}, we can construct an RS-tree decomposition $(X,T)$ with the same bags as $d$, and then by \cite[Lemma 12.3.4]{diestel2010graphtheory}, we know that every complete subgraph of $\Gaif(H)$ will belong to a single bag of $(X,T)$ and hence to $d$. Thus $d$ can be extended uniquely to a structured decomposition $d' : \smallint T \to \Hyp$ of $H$ with the same width as $d$. Hence the inequality follows.
\end{proof}

\begin{Cor} \label{cor hypergraph treewidth equal to treewidth of Gaifman graph}
Given a hypergraph $H$,
\begin{equation*}
    \mathbf{w}_\Hyp(H) = \tw(\Gaif(H)) = \tw_\Hyp(H).
\end{equation*}
\end{Cor}

\subsection{Carmesin's Graph Decompositions} \label{section caremsin's graph decompositions}

In \cite[Definition 9.3]{carmesin2022}, Carmesin introduced the notion of graph-decomposition, which we give in the language of structured decompositions. Let $\text{mltGraphs}$ denote the set of multigraphs.

\begin{Def}
Given a multigraph $H$, a \textbf{graph decomposition} of $H$ is a $\text{mltGraphs}$-structured decomposition $d: \smallint G \to \ncat{Gr}_\mlt$ of $H$.
\end{Def}

Now one might wonder why the definition was not given as a structured decomposition in $\Gr$. This is because $\Gr$ does not have colimits of all $\text{Graphs}$-structured decompositions. For example, consider the following structured decomposition $d: \smallint K^3 \to \ncat{Gr}$

\begin{equation*}
\begin{tikzcd}
	&& {K^2} \\
	& {\color{winered}K^1} & {\color{winered}K^1} \\
	{K^1} & {\color{winered}K^1} & {K^1}
	\arrow["f", hook, from=2-2, to=1-3]
	\arrow[hook', from=2-2, to=3-1]
	\arrow["g"', hook', from=2-3, to=1-3]
	\arrow[hook, from=2-3, to=3-3]
	\arrow[hook', from=3-2, to=3-1]
	\arrow[hook, from=3-2, to=3-3]
\end{tikzcd}
\end{equation*}
where the images of $f,g$ are disjoint and we colored the adhesions red for readability. The diagram $d$ does not have a colimit in $\Gr$. Such a colimit would have to coequalize the morphisms $f$ and $g$. But constructing a map $h : K^2 \to G$ so that $hf = hg$ requires sending both vertices of $K^2$ to the same vertex in $G$, and hence collapsing the edge of $K^2$. This cannot be done in $\Gr$.


Carmesin works around this issue by taking the colimit in the category $\Gr_\mlt$ of multigraphs, which is  finitely cocomplete (see Definition \ref{def multigraphs}). This also has the effect that the indexing graph itself must be allowed to be a multigraph, see \cite[Example 9.6, 9.7]{carmesin2022}.

Let $\text{mltGraphs}_{\leq n}$ denote the set of multigraphs with at most $n$ vertices.

\begin{Prop}
The triple $(\Gr_\mlt, \text{mltGraphs}, \{ \text{mltGraphs}_{\leq n} \}_{n \geq 0} )$ is a width category. Furthermore, the width of a $\text{mltGraphs}$-structured decomposition of a graph $G$ is precisely the same as the width of its corresponding graph-decomposition in the sense of Carmesin.
\end{Prop}

\begin{proof}
It is easy to see that (\textbf{W2}) holds, and (\textbf{W1}) holds because $\Gr_\mlt$ has all finite colimits. Now (\textbf{W3}) holds by Corollary \ref{cor finite colimits of multigraphs are pullback-stable}. Carmesin defines the width of a graph-decomposition $d$ of a graph $G$ as the maximum number of vertices of its bags, which in this case is precisely the size of the structured decomposition.
\end{proof}

Note that Carmesin does not actually ever define the graphwidth of a graph $G$. In other words, he considers width only as a propert of graph decompositionss, not as a property of graphs. This is explained by the following result.

\begin{Lemma}
Given $\Gamma = (\Gr_\mlt, \text{mltGraphs}, \{ \text{mltGraphs}_{\leq n} \}_{n \geq 0} )$, every multigraph $G$ with more than one vertex has $\Gamma$-width $1$.
\end{Lemma}

\begin{proof}
Every multigraph can be glued together all at once using colimits from three kinds of bags, the graph $\bullet$, the multigraph with one vertex and one loop, and the multigraphs with two vertices and $n$-many edge between them, all of which have size less than or equal to $2$.
\end{proof}

\subsection{Layered Treewidth}

Layered treewidth is a graph invariant introduced independently in \cite{shahrokhi2015layered} and \cite{dujmovic2017layered}. Unlike treewidth, layered treewidth is bounded on the class of planar graphs \cite[Theorem 12]{dujmovic2017layered}.

\begin{Def}[{\cite[Section 2]{Bose2022}}]
A \textbf{layering} $L$ for a graph $G$ consists of a partition $(L_0, L_1, \dots, L_n, \dots)$ of $V(G)$ such that if $xy \in E(G)$, and $x \in L_i$, $y \in L_j$, then $|i - j| \leq 1$.
\end{Def}

There is a convenient way to think about layerings using graph homomorphisms. Let $\Gr^\infty, \Gr_\ell^\infty, \Gr_r^\infty$ and $\Gr_\mlt^\infty$ denote the categories of graphs corresponding to those in Section \ref{section overview categories of graphs} but whose sets of vertices (and sets of edges in the case of multigraphs) can now be infinite. We can consider the faithful inclusion $i : \Gr \hookrightarrow \Gr_r^\infty$. Then a layering of a graph $G$ is equivalent to a morphism $G \to P^\infty_r$ in $\Gr_r^\infty$, where $P^\infty_r$ is the infinite reflexive graph with $V(P^\infty_r) = \N$ and where for $x,y \in V(P^\infty_r)$, there is an edge $xy$ if and only if $|x - y | \leq 1$. Let $\Gr_{\Lyr}$ denote the comma category $(i \downarrow P_r^\infty)$. We call the objects of $\Gr_\Lyr$ \textbf{layered graphs} and the morphisms \textbf{layered graph homomorphisms}. Note that every graph $G$ has at least one layering, by say putting all of the vertices of $G$ in $L_0$. Let $\pi : \Gr_\Lyr \to \Gr$ denote the projection functor.

Note that the inclusion $i : \Gr \hookrightarrow \Gr^\infty_r$ preserves pushouts along monomorphisms and finite limits. Hence by Lemma \ref{lem computing (co)limits in comma categories}, $\Gr_\Lyr$ has pushouts along monomorphisms and finite limits, and they are computed as in $\Gr$. 

Let $\mathbf{i} = (i_0, \dots, i_n, \dots)$ denote an infinite sequence of non-negative integers where only finitely many of the $i_k$ are non-zero. We say $\mathbf{i}$ is a \textbf{sequence of finite support} and let $\mathbf{I}$ denote the set of all sequences of finite support. Given a sequence $\mathbf{i} \in \mathbf{I}$, let $K^\mathbf{i}$ denote the following layered graph. Its $n$th layer $L_n$ is the complete graph $K^{i_n}$ (where recall that $K^0 = \varnothing$), and for every $x \in V(L_{n-1}) = V(K^{i_{n-1}})$ and $y \in V(L_n) = V(K^{i_n})$, there is an edge $xy \in K^\mathbf{i}$ for $n \geq 1$.

\begin{Ex}
If $\mathbf{i} = (1,2,2, 0, \dots, 0, \dots)$, then $K^\mathbf{i}$ is the graph
\begin{equation*}
    \begin{tikzcd}
	&& \bullet & \bullet \\
	{K^{\mathbf{i}}} & \bullet & \bullet & \bullet \\
	{P^\infty_r} & \bullet & \bullet & \bullet & \dots
	\arrow[no head, from=1-3, to=1-4]
	\arrow[no head, from=1-3, to=2-3]
	\arrow[no head, from=1-3, to=2-4]
	\arrow[no head, from=1-4, to=2-4]
	\arrow[no head, from=2-2, to=1-3]
	\arrow[no head, from=2-2, to=2-3]
	\arrow[no head, from=2-3, to=1-4]
	\arrow[no head, from=2-3, to=2-4]
	\arrow[no head, from=3-2, to=3-2, loop, in=300, out=240, distance=5mm]
	\arrow[no head, from=3-2, to=3-3]
	\arrow[no head, from=3-3, to=3-3, loop, in=300, out=240, distance=5mm]
	\arrow[no head, from=3-3, to=3-4]
	\arrow[no head, from=3-4, to=3-4, loop, in=300, out=240, distance=5mm]
\end{tikzcd}
\end{equation*}
\end{Ex}

Given a sequence of finite support $\mathbf{i}$, we call $K^\mathbf{i}$ a \textbf{layered complete graph}, and we let $\max(\mathbf{i})$ denote the maximal $i_k \in \mathbf{i}$. We let $\Omega_n$ denote the set of layered complete graphs $K^\mathbf{i}$ where $\max(\mathbf{i}) \leq n$. Thus there exists a monomorphism $G \hookrightarrow K^\mathbf{i}$ if and only if the maximal number of vertices of each layer of $G$ is less than or equal to $\max \mathbf{i}$.

\begin{Def}[{\cite[Section 2]{Bose2022}}] \label{def layered tree decomposition}
A \textbf{layered tree decomposition} $(L,d)$ of a graph $G$ consists of a layering $L$ of $G$ along with a tree decomposition $d$ of the underlying graph $G$. The width of a layered tree decomposition $(L,d)$ of a graph $G$ is the maximal number of vertices in the intersection of the bags of $d$ and the layers of $L$, i.e. the width of $(L,d)$ is $\max |L_i \cap d(x)|$.
\end{Def}

Let $L : G \to P_r^\infty$ be a layered graph. Note that if $d: \smallint T \to \Gr_\Lyr$ is a structured-decomposition of $G$ in $\Gr_\Lyr$, where $T$ is a tree, then we obtain a tree decomposition of $G$ in $\Gr$. Furthermore, using Lemma \ref{lem computing (co)limits in comma categories}, it is easy to see that to every tree decomposition $d : \smallint T \to \Gr$ of $G$ taken in $\Gr$, there is a unique diagram $d' : \smallint T \to \Gr_\Lyr$ with $\pi d' = d$, where each bag and adhesion $d'(x)$ of $d'$ are given the layering $d'(x) = d(x) \to P^\infty_r$ given by composing with the colimit cocone maps $d(x) \to G \to P^\infty_r$.

In other words, structured decompositions of the form $d : \smallint T \to \Gr_\Lyr$ where $T$ is a tree are equivalent to layered tree decompositions.

\begin{Prop} \label{prop sd-category for layered treewidth}
Layered graphs form a width category $\Gamma = (\Gr_\Lyr, \text{Trees}, \{K^\mathbf{i} \}_{\mathbf{i} \in \mathbf{I}} )$. Furthermore, if $d : \smallint T \to \Gr_\Lyr$ is a structured decomposition, then the width of $d$ is precisely the width of the layered tree decomposition $(L,d)$ as in Definition \ref{def layered tree decomposition}. 
\end{Prop}

\begin{proof}
It is easy to check that the conditions for $\Gamma$ to be a monic-stable sd-category hold thanks to Lemma \ref{lem computing (co)limits in comma categories}. The size of a layered graph $L : G \to P^\infty_r$ is clearly the largest number of vertices of its layers. Hence by definition, the width of a decomposition $d : \smallint T \to \Gr_\Lyr$ is the largest bag size $s(d(x))$ minus $1$, which is precisely the maximal number of vertices in the largest layer of each bag minus $1$.
\end{proof}

\begin{Def}[{\cite[Section 2]{Bose2022}}]
The \textbf{layered treewidth} of a graph $G$ is the minimal width of every layered tree decomposition $(L,d)$ of $G$.
\end{Def}

Hence the layered treewidth of a graph $G$ is the minimal $\Gamma$-width of each layered graph in the fiber $\pi^{-1}(G)$ where $\pi : \Gr_\Lyr \to \Gr$ is the projection functor.

\subsection{H-Treewidth}

Given a class of graphs $\mathcal{H}$, $\mathcal{H}$-treewidth is a graph invariant introduced in \cite{eiben2021measuring} for the purpose of obtaining ``hybrid'' width parameters.

\begin{Def}[{\cite[Definition 3.4]{jansen2022vertex}}]
Given a class of graphs $\mathcal{H}$, an \textbf{$\mathcal{H}$-tree decomposition} $(X, T, L)$ of a graph $G$ consists of a tree decomposition $(X, T)$ of $G$ along with a subset $L \subseteq V(G)$ such that 
\begin{enumerate}
    \item for each $v \in L$, there exists a unique leaf $t \in T$ with $v \in X(t)$, and
    \item for each $t \in T$, the induced subgraph $G[X(t) \cap L]$ belongs to $\mathcal{H}$.
\end{enumerate}
The width of an $\mathcal{H}$-tree decomposition $(X,T,L)$ is 
\begin{equation*}
    w(X,T,L) = \max(0, \max_{t \in T} |X(t) \setminus L| - 1)
\end{equation*}
The $\mathcal{H}$-treewidth $\tw_{\mathcal{H}}(G)$ of a graph $G$ is defined to be the minimal width of all of its $\mathcal{H}$-tree decompositions.
\end{Def}

\begin{Rem} \label{rem H tree decomps}
If $v \in L$, then by (1), it must belong to a unique bag $X(t)$ where $t \in T$ is a leaf vertex. This implies that $X(t) \cap L$ is either a nonempty subset of $X(t)$ and $t$ is a leaf node or $X(t) \cap L = \varnothing$. In other words, $L$ must be a subset of a union of bags indexed by leaf vertices of $T$, and each induced subgraph must belong to $\mathcal{H}$.
\end{Rem}

\begin{Lemma}
If $(X,T,L)$ is a $\mathcal{H}$-tree decomposition of a graph $G$, let $(X', T', L')$ denote the $\mathcal{H}$-tree decomposition of $G$ where to every leaf vertex $t$ we set $X(t) = X(t) \setminus L$, add a new vertex $t'$ and a new bag $X(t') = X(t) \cap L$.
\end{Lemma}

We can easily obtain a width measure on graphs which is comparable to $\mathcal{H}$-treewidth as follows. Let $\Omega^\mathcal{H}$ be the spine of $\Gr$ where $\Omega^\mathcal{H}_n = \{ K^n \}_{n \geq 0} \cup \mathcal{H}$ for $n \geq 0$.

\begin{Prop} \label{prop sd-category for H-treewidth}
The triple $\Gamma = (\Gr, \text{Trees}, \Omega^\mathcal{H})$ forms a width category.
\end{Prop}

\begin{proof}
Conditions (\textbf{W1}) and (\textbf{W3}) follow from Proposition \ref{prop sd-category for treewidth}. Condition (\textbf{W2}) holds because every graph has a monomorphism to some complete graph, and because there are no monos $K^n \hookrightarrow K^m$ if $m < n$, and each graph in $\mathcal{H}$ belongs to each $\Omega^\mathcal{H}_n$, so (\textbf{W2b}) holds there vacuously.
\end{proof}

\begin{Rem}
Note that the size of each $H \in \mathcal{H}$ with respect to $\Gamma$ is $0$.
\end{Rem}

\begin{Prop}
Given a graph $G$, we have
\begin{equation*}
    \mathbf{w}_\Gamma(G) \leq \tw_\mathcal{H}(G).
\end{equation*}
\end{Prop}

\begin{proof}
Given a $\mathcal{H}$-tree decomposition $(X,T,L)$, let us construct a new $\mathcal{H}$-tree decomposition $(X',T', L')$ with the same width as $(X,T,L)$. For each leaf vertex $t \in T$, we set $X'(t) = X(t) \setminus L$. We add a new vertex $t'$ and a single edge $tt'$ to $T$ to form $T'$, and set $X'(t') = X(t) \cap L$. Immediately it is clear that $w(X,T,L) = w(X',T',L)$, and it is also clear that $(X', T')$ is still a tree decomposition of $G$.

Now let $d : \smallint T \to \ncat{Gr}$ denote the structured decomposition of $G$ corresponding to $(X',T')$. Then by construction $w(d) \leq w(X', T', L')$, as $d$ may contain bags which belong to $\mathcal{H}$ that are not leaf bags. So if $(X,T,L)$ is a minimal $\mathcal{H}$-tree decomposition of $G$, then
\begin{equation*}
    \mathbf{w}_\Gamma(G) \leq w(d) \leq w(X',T',L) = w(X,T,L) = \tw_\mathcal{H}(G).
\end{equation*}
\end{proof}

It is easy to see that the reverse inequality does not hold in general. For example, let $G$ be the following graph
\begin{equation*}
    \begin{tikzcd}
	& \bullet && \bullet \\
	\bullet && \bullet && \bullet
	\arrow[no head, from=1-2, to=1-4]
	\arrow[no head, from=1-4, to=2-5]
	\arrow[no head, from=2-1, to=1-2]
	\arrow[no head, from=2-1, to=2-3]
	\arrow[no head, from=2-3, to=1-2]
	\arrow[no head, from=2-3, to=1-4]
	\arrow[no head, from=2-3, to=2-5]
\end{tikzcd}
\end{equation*}
and let $\mathcal{H} = \{K^3 \}$. Then there is an obvious tree decomposition where each bag is a $K^3$. So $\mathbf{w}_\Gamma(G) = 0$, but there exists no $\mathcal{H}$-tree decomposition of $G$ with width less than $2$.

\subsection{Bass-Serre Theory} \label{section bass-serre theory}

A well-studied precursor of our structured decompositions originates outside of combinatorics, in the field of geometric group theory. Specifically, \textit{graphs of groups}, the central objects of Bass-Serre theory, correspond exactly to structured decomposition of the form $d: \smallint G \rightarrow \ncat{Grp}$ in the category of groups and homomorphisms.

By general algebraic considerations, any diagram
\begin{tikzcd}
A & C \arrow[r, "b", hook] \arrow[l, "a"', hook'] & B
\end{tikzcd}
in the category $\ncat{Grp}$ has a colimit. One can construct the colimit as the quotient of the free product $A \ast B$ by the normal closure of the set $\left\{ a(x)b(x)^{-1} \mid x \in C \right\}$. As customary in the group theory literature, we call the colimit a \textit{free product with amalgamation} of $A, B$ over $C$ and (suppressing the monomorphisms) denote it $A \ast_C B$.

Fix a group $H$ and a graph $D$. Bass-Serre theory (introduced in \textit{Arbres, amalgames, $SL_2$} \cite{serre1977arbres}, later translated into English under the title \textit{Trees}~\cite{serre1980trees}) exhibits a correspondence between:
\begin{itemize}
    \item actions of $H$ on barycentric subdivisions of a tree\footnote{Here and only here, we allow infinite trees!} $T$ with quotient graph $D$, and
    \item graphs of groups on $D$ with \textit{fundamental group} $H$,
where a \textbf{graph of groups} is precisely a (not necessarily tree-shaped) structured decomposition $d : \smallint G \to \ncat{Grp}$.
\end{itemize}

As long as the quotient graph itself is a tree, the colimit of the graph of groups $d: \smallint D \rightarrow \mathrm{Grp}$ coincides with the fundamental group. Thus, by Bass-Serre theory, tree-structured decompositions of a group $G$ correspond to well-behaved actions of $G$ on trees.

Using this correspondence, one can equip the category of groups with the structure of a spined sd-category, and characterize the resulting width notion. However it turns out that $\Gamma$-width is less interesting than the $(\G, \Omega)$-width of Definition \ref{def chordal width}.

This also serves as an example of a spined sd-category which is not stable, but is nonetheless well-behaved enough to admit a notion of monotonicity~(Corollary~\ref{cor bass serre monotonicity}) for $(\G, \Omega)$-width. We will give a simple example of an action with tree quotient in Example~\ref{ex bass serre dihedral action}.

\begin{Def}
For $n\in\mathbb{N}$, let $\Omega_n$ denote the class of all groups with order at most $n$. The \textbf{sd-category of groups} is the sd-category $\Gamma_{\text{BS}} = (\ncat{Grp}, \mathrm{Trees}, \Omega)$ where $\Omega = \left\{\Omega_n \mid n \geq 0 \right \}$.
\end{Def}

Clearly, $\Gamma_{\text{BS}}$ is a spined sd-category. Let $w_{\text{BS}}$ denote the $(\G, \Omega)$-width of Definition \ref{def chordal width} for $\Gamma_{\text{BS}}$. We say that an object $X$ in an sd-category $(\cat{C}, \G, \Omega)$ has infinite $(\G, \Omega)$-width if there exists no $(\G, \Omega)$-structured decomposition of $X$, i.e. $X$ is not chordal.

Now since $\ncat{Grp}$ contains both finite and infinite groups, let us begin by characterizing the behavior of $w_{\text{BS}}$ on finite groups. To aid with this, we introduce the notion of inseparable object.

\begin{Def} \label{def inseparable object}
Consider a category $\cat{C}$ and an object $K$ of C. We call $K$ \textit{inseparable} if for any pushout $A +_C B$ and monomorphism $\iota: K \hookrightarrow A +_B C$ in $\cat{C}$, one can find a monomorphism $\iota_A : K \hookrightarrow A$ or a monomorphism $\iota_B: K \hookrightarrow B$.
\end{Def}

Keep in mind that Definition~\ref{def inseparable object} does not imply that the monomorphism $\iota : K \rightarrow A +_B C$ factors through either $\iota_A$ or $\iota_B$.

\begin{Lemma}\label{lemma finite groups are inseparable}
In the category $\ncat{Grp}$ of groups, every finite group is inseparable.
\begin{proof}
Immediate from \cite{serre1980trees}, Theorem 8.
\end{proof}
\end{Lemma}

\begin{Cor}\label{cor finite group width}
For a finite group $G$, we have $w_{\text{BS}}(G) = |G|$.
\begin{proof}
Since the one-vertex graph is a permitted structure graph, we have $w_{\text{BS}}(G) \leq |G|$. Now assume that $G$ admits a $(\text{Trees}, \Omega)$-structured decomposition $d: \smallint T \rightarrow \ncat{Grp}$ of width $n < |G|$. Then we have some vertex $v \in V(T)$ so that $|d(v)| < n$. But $G$ is a subobject of itself (i.e. the colimit of $d$), so Lemma~\ref{lemma finite groups are inseparable} guarantees that $G$ is a subobject of $d(v)$. But then $|G| \leq n < |G|$, a contradiction. Thus, $w_{\text{BS}}(G) = |G|$.
\end{proof}
\end{Cor}

Before characterizing the width of infinite groups, we present the simplest nontrivial example of a group acting on the barycentric subdivision of a tree with a tree quotient.

\begin{Ex} \label{ex bass serre dihedral action}
Consider the infinite path $P_\mathbb{Z}$ with vertex set the integers ($V(P_\mathbb{Z}) = \mathbb{Z}$), and edges connecting $x,y \in \mathbb{Z}$ precisely when $|x - y| = 1$. We depict the resulting tree (all black vertices) and its barycentric subdivision (black and white vertices) below:
\begin{equation*}\begin{tikzcd}
	{~} & {\bullet_{-3}} & {\bullet_{-2}} & {\bullet_{-1}} & {\bullet_0} & {\bullet_1} & {\bullet_2} & {\bullet_3} & {~}
	\arrow[dotted, no head, from=1-2, to=1-1]
	\arrow[no head, from=1-3, to=1-2]
	\arrow[no head, from=1-4, to=1-3]
	\arrow[no head, from=1-5, to=1-4]
	\arrow[no head, from=1-5, to=1-6]
	\arrow[no head, from=1-6, to=1-7]
	\arrow[no head, from=1-7, to=1-8]
	\arrow[dotted, no head, from=1-8, to=1-9]
\end{tikzcd}\end{equation*}
\begin{equation*}\begin{tikzcd}
	{~} & {\circ_{\frac{-3}{2}}} & {\bullet_{-1}} & {\circ_{\frac{-1}{2}}} & {\bullet_0} & {\circ_{\frac{1}{2}}} & {\bullet_1} & {\circ_{\frac{3}{2}}} & {~}
	\arrow[dotted, no head, from=1-2, to=1-1]
	\arrow[no head, from=1-3, to=1-2]
	\arrow[no head, from=1-4, to=1-3]
	\arrow[no head, from=1-5, to=1-4]
	\arrow[no head, from=1-5, to=1-6]
	\arrow[no head, from=1-6, to=1-7]
	\arrow[no head, from=1-7, to=1-8]
	\arrow[dotted, no head, from=1-8, to=1-9]
\end{tikzcd}\end{equation*}
The infinite dihedral group $D_\infty$, presented as $\left\langle a,b \mid a^2=b^2 = 1 \right\rangle$ acts on the integers by
\begin{align*}
    & a(x) = 0 -x \\
    & b(x) = 1 - x
\end{align*}
We can extend this to an action of $D_\infty$ on the barycentric subdivision of the tree $P_\mathbb{Z}$. The generator $a$ then acts by turning the tree 180 degrees around the zero vertex, whereas the generator $b$ acts by turning the tree 180 degrees around the midpoint $\frac{1}{2}$ between $0$ and $1$ ($b)$. Computing the quotient yields the graph $K^2$,
\begin{equation*}
\begin{tikzcd}
	\bullet & \circ
	\arrow[no head, from=1-1, to=1-2]
\end{tikzcd}
\end{equation*}
and $K^2$ is a tree. By Bass-Serre theory, we can thus write $D_\infty$ as the colimit of some $K^2$-structured decomposition. The theory even lets us determine the appropriate bags and adhesions by considering vertex and edge stabilizers, but here we can read off such a decomposition purely from the presentation of $D_\infty$ given above: $D_\infty$ arises as the colimit of the span
\begin{equation*}
\begin{tikzcd}
	{\mathbb{Z}/2\mathbb{Z}} & {\{1\}} & {\mathbb{Z}/2\mathbb{Z}}
	\arrow[hook', from=1-2, to=1-1]
	\arrow[hook, from=1-2, to=1-3]
\end{tikzcd}
\end{equation*}
i.e. the diagram $d : \smallint K^2 \rightarrow \ncat{Grp}$ with bags $d(\bullet) = d(\circ) = \mathbb{Z}/2\mathbb{Z}$ and adhesion $d(-) = \{1\}$.
\end{Ex}

\begin{Rem}
Example~\ref{ex bass serre dihedral action} shows that $w_{\text{BS}}(D_\infty) \leq 2$. Hence the width notion $w_{\text{BS}}$ is nontrivial. It also shows that infinite groups are generally not inseparable in $\ncat{Grp}$. For example, the morphism that maps $1$ to the product of the two generators $ab$ is a monomorphism from $\mathbb{Z}$ to $D_\infty$, even though no monomorphism of signature $\mathbb{Z} \hookrightarrow \mathbb{Z}/2\mathbb{Z}$ exists. This also shows the failure of monic-stability in the sd-category $\Gamma_{\text{BS}}$.
\end{Rem}

We now characterize the width of infinite groups. Specifically, we will show that for every bounded width group $G$, the order of the largest finite subgroup $W < G$ coincides with $w_{\text{BS}}(G)$.

\begin{Prop} \label{prop bass serre lower bound}
Consider a group $G$ such that $w_{\text{BS}}(G)$ is finite and let $W < G$ be a finite subgroup. Then every decomposition of $G$ has width bounded below by $|W|$.
\begin{proof}
Consider a tree decomposition $d : \smallint T \rightarrow \ncat{Grp}$ of width $n$. By repeated applications of Lemma~\ref{lemma finite groups are inseparable}, $W$ is a subobject of some bag $d(v)$. But then $|W| \leq d(v) \leq n$.
\end{proof}
\end{Prop}

We now show that the order of a maximal finite subgroup also upper bounds the width of the group. Recall that a group is \textit{virtually free} if it has a free subgroup of finite index. Since the free group on no generators (the trivial group $\{1\}$) has finite index in any finite group, in what follows we regard finite groups as virtually free.

\begin{Lemma}\label{lemma subgroup of virtually free group}
Subgroups of virtually free groups are virtually free.
\begin{proof}
Immediate from the Nielsen-Schreier theorem (subgroups of free groups are free).
\end{proof}
\end{Lemma}

\begin{Lemma}\label{lemma finite width group is virtually free}
Groups of finite width are virtually free.
\begin{proof}
See \cite{serre1980trees}, Proposition 11.
\end{proof}
\end{Lemma}

\begin{Rem}
Lemma~\ref{lemma finite width group is virtually free} gives an easy example of a group with unbounded width: $w_{\text{BS}}(\mathbb{Z}^2) = \infty$.
\end{Rem}

\begin{Prop} \label{prop bass serre upper bound}
Consider a finite-width group $G$ and let the nontrivial finite subgroup $W < G$ have maximal order. Then $G$ has a decomposition of width at most $|W|$.
\begin{proof}
Since $G$ has finite width, all finite subgroups of $G$ have order at most $v$. A result of Linnel~\cite{linnel1983accessibility} shows that one can write all such groups as free products with amalgamation where the edge groups are all finite, and the vertex groups are subgroups of $G$ with Cayley graphs that have at most one end. By Lemma~\ref{lemma finite width group is virtually free} the vertex groups are themselves virtually free. The Cayley graph of a virtually free group has either zero, two or infinitely many ends. Thus, in this case the vertex groups have no ends, and are finite as required. This yields a decomposition of $G$ whose bags are finite subgroups of $F$. Since $W$ has maximal order among the finite subgroups, this decomposition has width at most $|W|$ as required.
\end{proof}
\end{Prop}

From Propositions~\ref{prop bass serre lower bound}~and~\ref{prop bass serre upper bound} it follows that, as long as a group $\Gamma$ has a nontrivial finite subgroup of maximal order, its width is this maximal order.

\begin{Cor}\label{cor bass serre monotonicity}
Consider a group $G$. Then either $w_{\text{BS}}(G) = \infty$, or else $w_{\text{BS}}(G)$ is the order of a maximal-order finite subgroup of $G$. Consequently, when $H < G$ and both $w_{\text{BS}}(G) < \infty$ and $w_{\text{BS}}(H) < \infty$ hold, then $w_{\text{BS}}(H) \leq w_{\text{BS}}(G)$ holds as well.
\end{Cor}

The failure of pullback-stability and monic-stability in $\Gamma_{\text{BS}}$ sheds light on a general phenomenon: in $\ncat{Gr}$, subobjects of well-behaved objects usually remain similarly well-behaved. This is very far from the case in algebra: for example, a subgroup of a finitely presented group will in general fail to have a finite presentation. Moreover, via the Bass-Serre correspondence, the general theory developed in Section~\ref{section sd-categories} lets us conclude facts purely about group actions as corollaries. We present one example of this below.

\begin{Cor}\label{cor bass serre binary tree actions}
If a group $G$ acts on the barycentric subdivision of a tree with a tree-shaped quotient $T$, then $G$ acts on some other tree with a binary tree quotient $T'$.
\begin{proof}
Use Bass-Serre theory to obtain a $(\ncat{Gr}, \mathrm{Trees}, \Omega)$-decomposition of $G$. Apply Proposition \ref{prop binarization functor} to obtain a $(\ncat{Gr}, \mathrm{BTrees}, \Omega)$-decomposition of $G$. Using the Bass-Serre correspondence in reverse yields the desired action.
\end{proof}
\end{Cor}

\subsection{Hybrid Dynamical Systems}

The thesis \cite{ames2006categorical} sets up a framework with which to discuss hybrid systems categorically, where here a hybrid system is a dynamical system that has both discrete and continuous components. 

\begin{Def}[{\cite[Definition 1.6]{ames2006categorical}}]
A \textbf{D-category} is a small category $\cat{C}$ such that
\begin{enumerate}
    \item for every object $X \in \cat{C}$, $X$ is either the domain or codomain of a non-identity morphism, but not both, and
    \item if $X \in \cat{C}$ is the domain for a non-identity morphism $f : X \to Y$, then there exists exactly one other morphism $g : X \to Z$ in $\cat{C}$ with domain $X$.
\end{enumerate}
\end{Def}

Clearly if $G$ is a (possibly infinite) digraph, then $ \smallint G$ is a $D$-category. Ames defines a corresponding notion of functor of $D$-categories, resulting in a category $\ncat{DCat}$ of $D$-categories, and the Grothendieck construction extends to a functor $\smallint : \dGr^\infty \to \ncat{DCat}$.

\begin{Th}[{\cite[Theorem 1.1]{ames2006categorical}}]
The functor
\begin{equation*}
    \smallint : \dGr^\infty \to \ncat{DCat},
\end{equation*}
is an isomorphism of categories.
\end{Th}

\begin{Def}[{\cite[Definition 1.9]{ames2006categorical}}]
A \textbf{hybrid object} in a category $\cat{C}$ is a functor $d : I \to \cat{C}$ where $I$ is a $D$-category.
\end{Def}

So hybrid objects for Ames are precisely our structured decompositions. The main object of study for this thesis is the notion of a hybrid system, which turns out \cite[Proposition 2.1]{ames2006categorical} to be equivalent to certain hybrid objects in a category of smooth manifolds.

\section{Future Work} \label{section future work}

In this brief section we consider possible avenues for future research.

\subsection*{Other Shapes}
Recall that a \(G\)-shaped structured decomposition valued in some category \(\cat{C}\) is just a diagram \(d \colon \smallint G \to \cat{C}\) whose shape is given by the Grothendieck construction applied to a \textit{graph} \(G\). In principle there is nothing stopping us from letting \(G\) be another kind of presheaf and indeed, one could for instance study structured decompositions whose shapes are given not by graphs, but simplicial complexes or other combinatorial objects.

\subsection*{Obstructions to Low Width}
In graph-theory, one natural direction for further work is to characterize what obstacles prevent general objects from admitting structured decompositions of low width. When it comes to graphs and tree-shaped structured decompositions, there is a well-established theory of obstacles to having low tree-width, including notions such as \textit{brambles}~\cite{seymour-thomas-BRAMBLES}, \textit{$k$-blocks}~\cite{carmesin-k-blocks}, \textit{tangles}~\cite{robertsonX} and \textit{abstract separation systems}~\cite{diestel2019profiles}. It is a fascinating, but highly non-trivial research direction to lift these ideas to the more general, category-theoretic setting of structured decompositions. 

\subsection*{Co-Decompositions}
Zooming out, the focus of this paper has been to investigate which combinatorial invariants can
be expressed as colimits. In particular, in a structured decomposition category we think of objects as being
decomposed if they arise as colimits of a structured decompositions. An obvious question for future work is to dualize these ideas: rather than building objects as colimits, we might want to build objects as \textit{limits}. Thus one has the notion of a \textbf{strctured \textit{co}-decomposition}, namely diagrams of the form $d \colon (\smallint G)^{\op} \to \cat{C}$ which are thought to decompose an object $X \in \cat{C}$ if $\lim d \cong X$. It is not at all clear what graph classes can be captured via structured co-decompositions and indeed even the case of path-shaped structured co-decompositions of graphs seems to yield interesting graph classes. Preliminary investigation suggests that graphs of bounded tree-co-width should admit a highly symmetric structure, but much further research is still needed to understand the structure of such graphs. We believe that this is a fascinating direction for future work harboring many applications in graph theory and algorithmics.

\subsection*{Submodular Width Measures}

The astute reader will notice that nearly all of our examples come from tree shaped decompositions. Other important width measures for graphs exist that are not asymptotically equivalent to tree width, such as clique-width, mim-width and sim-width \cite{brettell2023comparingwidth}. Each of these width measures can be obtained using a submodular function and branch decompositions. Studying the relationship between these submodular width measures and width categories is future work.

\subsection*{Generalized Hypertreewidth}

In Section \ref{section hypergraph treewidth}, we showed how to capture hypergraph treewidth using width categories. Hypergraph treewidth is already a useful invariant of hypergraphs in algorithmics. For example the Conjunctive Query Containment Problem is polynomial time solvable on hypergraphs of bounded Gaifman treewidth \cite{kolaitis1998conjunctive}. However, as discussed in \cite[Section 5.2]{hlinveny2008width}, there exists a notion of hypergraph acyclicity, called $\alpha$-acyclicity \cite{fagin1983degrees}, which is not well behaved with respect to hypergraph treewidth. Namely there exist classes of $\alpha$-acyclic hypergraphs that have unbounded hypergraph treewidth. The asymptotically equivalent width measures of hypertreewidth and generalized hypertreewidth, introduced in \cite{gottlob1999hypertree}, do not suffer from this defect. We note that generalized hypertreewidth $\mathbf{ghw}(H)$ of a hypergraph $H$ can be defined using tree decompositions where the measure of size for the bags of the decomposition is given by the edge cover number. This is problematic for our formalism, as edge cover number is not monotonic with respect to monomorphisms of hypergraphs, and hence there exists no spine on $\Hyp$ which we can use to capture generalized hypertreewidth. So far we have been unable to find another category of hypergraphs where edge cover number is monotonic with respect to the monomorphisms. This suggests expanding the definition of spine and possibly also considering pushing forward structured decompositions (as in \cite{bumpus2024pushingforward}) rather than pulling them back. This avenue of research is future work.

\appendix

\section{Categories of Graphs} \label{section categories of graphs}

In this appendix, we detail some facts about categories of graphs. First we give an overview of several categories of graphs we are considering, and then examine the categories $\Gr$ and $\Gr_\mlt$ more deeply.

\subsection{Overview} \label{section overview categories of graphs}

There are many different notions of graphs, with different corresponding categories. These various categories, their properties and relationships are at this point well understood \cite{bumby1986categorical, brown2008graphs, plessas2011categories, schmidt2019, campion2021graphs}. Thus we will only quickly review the relevant categories for our purposes. 

\begin{Def} \label{def directed multigraphs}
Let $\ncat{dGrSch}$\footnote{which stands for directed graph schema.} denote the category with two objects and two non-identity morphisms
\begin{equation*}
\begin{tikzcd}
	E & V
	\arrow["s", shift left, from=1-1, to=1-2]
	\arrow["t"', shift right, from=1-1, to=1-2]
\end{tikzcd}
\end{equation*}
Then let $\ncat{dGr} = \mathsf{Fun}(\ncat{dGrSch}, \ncat{FinSet})$. We call the objects of this category \textbf{directed multigraphs}. By this we mean that the edges of the graph are directed, there can be multiple edges and loops at vertices. Morphisms are directed graph homomorphisms.
\end{Def}

\begin{Def} \label{def symmetric multigraphs}
Let $\ncat{sGrSch}$ denote the category given by
\begin{equation*}
    \begin{tikzcd}
	E & V
	\arrow["i", from=1-1, to=1-1, loop, in=150, out=210, distance=5mm]
	\arrow["s", shift left, from=1-1, to=1-2]
	\arrow["t"', shift right, from=1-1, to=1-2]
\end{tikzcd}
\end{equation*}
and where $t = si$, $s = ti$ and $i i = 1_{E}$. Let $\ncat{sGr} = \mathsf{Fun}(\ncat{sGrSch}, \ncat{FinSet})$. We call the objects of this category \textbf{symmetric multigraphs}.
\end{Def}

The objects of this category can be thought of as symmetric directed graphs, where the involution $i$ sends a directed edge to the same edge with the opposite direction. This does introduce a mild pathology however, as loops on vertices can be their own involution, and such loops cannot be mapped to arbitrary loops in other graphs\footnote{In \cite{brown2008graphs}, edges on a single vertex whose involution does not map to themselves are called bands, and those that do are called loops. In fact, it is known that in any conventional category of graphs that has exponential objects, this phenomenon must appear \cite[Introduction]{schmidt2019}.}.

We similarly have a reflexive version of $\ncat{sGr}$.

\begin{Def}
Let $\ncat{sGr}_r\ncat{Sch}$ denote the category given by
\begin{equation*}
\begin{tikzcd}
	E & V
	\arrow["i", from=1-1, to=1-1, loop, in=150, out=210, distance=5mm]
	\arrow["s", shift left=2, from=1-1, to=1-2]
	\arrow["t"', shift right=2, from=1-1, to=1-2]
	\arrow["u"{description}, from=1-2, to=1-1]
\end{tikzcd}
\end{equation*}
and where $t = si$, $s = ti$, $su = tu = 1_V$, $i u = u$ and $i i = 1_{E}$. Let $\ncat{sGr}_r = \mathsf{Fun}(\ncat{sGr}_r\ncat{Sch}, \ncat{FinSet})$. We call the objects of this category \textbf{reflexive symmetric multigraphs}.
\end{Def}

Although $\ncat{dGr}$, $\ncat{sGr}$ and $\ncat{sGr}_r$ are categorically nice, being toposes, we are primarily interested in undirected graphs. Let us introduce the categories of graphs under consideration.

We have already seen the category $\ncat{Gr}$ of graphs (Definition \ref{def graphs}). We will explore it more deeply in the next section. Another practical choice for a category of graphs is the category of loop graphs. 

\begin{Def} \label{def loop graphs}
By a \textbf{loop graph}, we mean a finite set equipped with a binary, symmetric relation. A morphism $f : G \to H$ of loop graphs consists of a function $V(f) : V(G) \to V(H)$ which preserves the relation. We let $\ncat{Gr}_\ell$ denote the category of loop graphs. We can visualize these objects as simple graphs with loops allowed.
\end{Def}

We call a loop graph \textbf{irreflexive} if it has no loops and \textbf{reflexive} if every vertex has a loop. Thus we can obtain $\Gr$ as the full subcategory of $\Gr_\ell$ on the irreflexive loop graphs, and we let $\Gr_r$ denote the full subcategory of reflexive loop graphs. The category $\ncat{Gr}_r$ can also be thought of as the category whose objects are simple graphs, but where morphisms are allowed to collapse edges.

We can characterize the categories $\Gr_\ell$ and $\Gr_r$ as reflective subcategories of $\ncat{sGr}$ and $\ncat{sGr}_r$ respectively.

\begin{Def}
We say that a small category $\cat{C}$ is \textbf{terminally concrete} if it has a terminal object $*$, and such that the functor
\begin{equation*}
    \cat{C}(*, -) : \cat{C} \to \ncat{Set}
\end{equation*}
is faithful. If $\cat{C}$ is terminally concrete, then a presheaf $X : \cat{C}^\op \to \ncat{Set}$ is \textbf{concrete} if for every $U \in \cat{C}$, the canonical map
\begin{equation*}
    X(U) \to \ncat{Set}(\cat{C}(*,U), X(*))
\end{equation*}
is injective.
\end{Def}

\begin{Prop}[{\cite[Lemma 47]{baez2011convenient}}]
Given a terminally concrete category $\cat{C}$, the category of concrete presheaves on $\cat{C}$ is a reflective subcategory of presheaves on $\cat{C}$
\begin{equation*}
\begin{tikzcd}
	{\ncat{ConPre}(\cat{C})} & {\ncat{Pre}(\cat{C})}
	\arrow[""{name=0, anchor=center, inner sep=0}, "i"', shift right=2, hook, from=1-1, to=1-2]
	\arrow[""{name=1, anchor=center, inner sep=0}, "{\text{Con}}"', shift right, from=1-2, to=1-1]
	\arrow["\dashv"{anchor=center, rotate=-90}, draw=none, from=1, to=0]
\end{tikzcd}
\end{equation*}
Furthermore this still holds if we consider categories of presheaves of finite sets.
\end{Prop}

The category $\Gr_\ell$ is precisely the category of concrete presheaves of finite sets on $\ncat{sGrSch}^\op$, and $\Gr_r$ is precisely the category of concrete finite presheaves of finite sets on $\ncat{sGr}_r\ncat{Sch}^\op$. Hence both $\Gr_\ell$ and $\Gr_r$ are quasitoposes. These are very nice categories, which are finitely (co)complete, locally cartesian closed, and have a weak subobject classifier \cite[Definition 8]{baez2011convenient}.

\begin{Def} \label{def multigraphs}
If $S$ is a finite set, let $S^{(2)}$ denote\footnote{We took this notation from \cite[Remark 2.2]{huh2020logarithmic}. We note that $S^{(2)}$ is isomorphic to the set of cardinality two multi-subsets of $S$, i.e. subsets of cardinality two that allow for repeated elements. Our category of multigraphs is precisely the category of graphs that Eur and Huh consider in their paper.} the set $S^2$ quotiented by the relation $(x,y) \sim (y,x)$. A \textbf{multigraph} $G$ consists of two finite sets $V(G)$, $E(G)$ and a function $\partial : E(G) \to V(G)^{(2)}$. We can visualize these as unordered graphs that can have multiple loops and parallel edges. Morphisms $ f: G \to H$ of multigraphs are then functions $V(f) : V(G) \to V(H)$ and $E(f) : E(G) \to E(H)$ making the following diagram commute
\begin{equation*}
\begin{tikzcd}
	{E(G)} & {E(H)} \\
	{V(G)^{(2)}} & {V(H)^{(2)}}
	\arrow["{{E(f)}}", from=1-1, to=1-2]
	\arrow["{{\partial_G}}"', from=1-1, to=2-1]
	\arrow["{{\partial_H}}", from=1-2, to=2-2]
	\arrow["{V(f)^{(2)}}", from=2-1, to=2-2]
\end{tikzcd}
\end{equation*}
\end{Def}

We note that the category of multigraphs has all finite limits and colimits but is \textit{not} cartesian closed \cite[Proposition 2.3.1]{plessas2011categories}.

\begin{Rem}
We see that there is a chain of fully faithful functors $\Gr \hookrightarrow \Gr_\ell \hookrightarrow \Gr_\mlt$, and furthermore note that each functor preserves finite limits and pushouts of spans of monomorphisms. There is a left adjoint $L_\mlt : \Gr_\mlt \to \Gr_\ell$ to the inclusion, given by collapsing multiple edges.

We can also embed $\Gr_r \hookrightarrow \Gr_\ell$ by just thinking of reflexive graphs as graphs where every vertex has a loop. This has a left adjoint $L_\ell : \Gr_\ell \to \Gr_r$ that just adds a loop to every vertex of a loop graph if it doesn't already have one.
\end{Rem}

\begin{Lemma}[{\cite[Proposition 2.3.13, 2.3.14]{plessas2011categories}}] \label{lem monos and epis of graphs}
A map $f : G \to H$ of (reflexive/loop) graphs is a mono/epimorphism in ($\Gr_r$/$\Gr_\ell$) $\Gr$ if and only if $V(f)$ is an in/surjective function. A map $f : G \to H$ of multigraphs is a mono/epimorphism in $\Gr_\mlt$ if and only if both $E(f)$ and $V(f)$ are in/surjective functions.
\end{Lemma}

\begin{Rem}\label{rem graphs notation comment}
We took the term loop graph from \cite{nlab:graph} and the term multigraph from \cite[Section 1.10]{diestel2010graphtheory}. In \cite{nlab:graph} what we call multigraphs are called pseudographs.
\end{Rem}

We summarize the relationships between our notation and that of \cite{plessas2011categories}, along with the properties of the categories.

\begin{table}[H]
    \centering
    \begin{tabular}{|c|c|c|c|}
    \hline
    \textsc{Graph name} & \textsc{Category} & \textsc{Plessas} & \textsc{Categorical properties}\\
    \hline
        simple graphs & $\ncat{Gr}$ & $\ncat{SiLlStGraphs}$ & products, pullbacks, pushouts of monos \\
        \hline
        loop graphs & $\ncat{Gr}_\ell$ & $\ncat{SiStGraphs}$ & quasitopos \\
        \hline
        reflexive loop graphs & $\Gr_r$ & $\ncat{SiGraphs}$ & quasitopos \\
        \hline
        multigraphs & $\ncat{Gr}_\mlt$ & $\ncat{StGraphs}$ & finite (co)limits, regular \\
    \hline
    \end{tabular}
    \caption{Different notation for categories of graphs.}
    \label{tab:my_label}
\end{table}

\subsection{Simple Graphs} \label{section gr}
The category that is perhaps closest to what combinatorialists and computer scientists think of for simple graphs we call $\ncat{Gr}$ (Definition \ref{def graphs}). This is the category whose objects are finite sets equipped with a binary, symmetric, irreflexive relation, and whose morphisms are functions that preserve the relation. We summarize some of the properties of this category below. For relations to other categories of graphs see Section \ref{section overview categories of graphs}

Let us work out some of the structure of $\ncat{Gr}$. First we note that there is a functor $V : \ncat{Gr} \to \ncat{FinSet}$ that gives the set of vertices of a graph. This functor has a left adjoint $\text{Disc} : \ncat{Set} \to \ncat{Gr}$ that sends a finite set $S$ to the discrete graph $\text{Disc}(S)$. Therefore, $V$ must preserve whatever limits $\ncat{Gr}$ has. So we can conclude that $\ncat{Gr}$ does not have a terminal object. If it did, call it $*$, then $V(*) \cong \mathbf{1}$, where $\mathbf{1}$ denotes a singleton set. There is only one graph with a single vertex, which we denote by $\bullet$. But any graph with an edge does not have a map to $\bullet$. Therefore $\ncat{Gr}$ does not have a terminal object. However it does have binary products. Given $G, H \in \ncat{Gr}$, $G \times H$ is the graph with $V(G \times H) = V(G) \times V(H)$ and where $(g,h)(g',h') \in E(G \times H)$ if and only if $gg' \in E(G)$ and $hh' \in E(H)$. The category $\ncat{Gr}$ also has pullbacks. If $f : G \to K$ and $g : H \to K$ are maps of graphs, then $G \times_K H$ is the graph with $V(G \times_K H) = V(G) \times_{V(K)} V(H)$ and where $(g,h)(g',h') \in E(G \times_K H)$ if and only if $gg' \in E(G)$ and $hh' \in E(H)$. Therefore it also has equalizers. Note that $\ncat{Gr}$ is not cartesian closed \cite[Proposition 2.3.6]{plessas2011categories}.

The category $\ncat{Gr}$ clearly has an initial object given by the empty graph $\varnothing$, and it has finite coproducts, where $G + H$ is the graph with $V(G + H) = V(G) + V(H)$ and $E(G + H) = E(G) + E(H)$. The category $\Gr$ does not have all pushouts in general, see \cite[Page 7]{bumpus2023spined}. However it does have pushouts along spans of monomorphisms.

\begin{Def} \label{def pushouts of spans of monos for graphs}
Let $i : G \to H$ and $j : G \to K$ be monomorphisms in $\Gr$. Let $H+_G K$ be the graph defined as follows. Let $V(H +_G K)$ be the pushout $V(H) +_{V(G)} V(K)$, i.e. the set $V(H) + V(K)/{\sim}$, where $\sim$ is the smallest equivalence relation such that for $h \in V(H)$ and $k \in V(K)$, $h \sim k$ if there exists a $g \in V(G)$ such that $i(g) = j(g)$. Since $i$ and $j$ are monomorphisms, $V(i)$ and $V(j)$ are injective. Given $h \in V(H)$, $k \in V(K)$, $[h] = [k]$ in $V( H +_G K)$ if and only if there exists a finite zig-zag of elements
\begin{equation*}
    \begin{tikzcd}
	& {g_0} && {g_1} && {g_2} && {g_n} \\
	h && {k_0} && {h_1} && \dots && k
	\arrow["i"', maps to, from=1-2, to=2-1]
	\arrow["j", maps to, from=1-2, to=2-3]
	\arrow["j"', maps to, from=1-4, to=2-3]
	\arrow["i", maps to, from=1-4, to=2-5]
	\arrow["i"', maps to, from=1-6, to=2-5]
	\arrow["j", maps to, from=1-6, to=2-7]
	\arrow["i"', maps to, from=1-8, to=2-7]
	\arrow["j", maps to, from=1-8, to=2-9]
\end{tikzcd}
\end{equation*}
i.e. $i(g_0) = h$, $j(g_0) = k_0$, $j(g_1) = k_0$, $\dots$, $j(g_n) = k$. But since $V(i)$ and $V(j)$ are injective, this means that $g_0 = g_1 = \dots = g_n$. Thus if $[h]=[k]$, then there exists a two-arrow zig-zag of elements 
\begin{equation*}
\begin{tikzcd}
	& g \\
	h && k
	\arrow["i"', maps to, from=1-2, to=2-1]
	\arrow["j", maps to, from=1-2, to=2-3]
\end{tikzcd}    
\end{equation*}
in other words, there exists a unique $g \in V(G)$ such that $i(g) = h$ and $j(g) = k$. Thus every element $x \in V(H +_G K)$ can either be represented uniquely by an $h \in V(H)$ such that $i^{-1}(h) = \varnothing$, a $k \in V(K)$ such that $j^{-1}(k) = \varnothing$, or there exists a unique $g \in V(G)$ such that $x = [i(g)] = [j(g)]$. This also implies that for $h, h' \in V(H)$, if $[h] = [h']$, then $h = h'$, and similarly for $V(K)$.

Now let us define $E(H +_G K)$. For $x,y \in V(H +_G K)$ there is an edge $xy \in E(H +_G K)$ if and only if $x \neq y$ and there exists an edge $hh' \in E(H)$ such that $x = [h]$, $y=[h']$ or $kk' \in E(K)$ such that $x=[k]$, $y=[k']$. Thus $H+_G K$ is a well-defined graph. 
   
\end{Def}

Let $r : H \to H+_G K$ and $s : K \to H+_G K$ be defined simply as the quotient map on vertices. They are easily seen to be graph morphisms.

\begin{Lemma} \label{lem pushouts of monos exist in Gr}
Given a span of monomorphisms in $\Gr$
\begin{equation} \label{eq pushouts of monos in Gr}
   \begin{tikzcd}
	G & K \\
	H & {H +_G K}
	\arrow["j", hook', from=1-1, to=1-2]
	\arrow["i"', hook, from=1-1, to=2-1]
	\arrow["s", hook', from=1-2, to=2-2]
	\arrow["r"', hook, from=2-1, to=2-2]
	\arrow["\lrcorner"{anchor=center, pos=0.125, rotate=180}, draw=none, from=2-2, to=1-1]
\end{tikzcd} 
\end{equation}
their pushout exists, and furthermore the maps $r$ and $s$ are monomorphisms as well.
\end{Lemma}

\begin{proof}
Suppose that $p, q: G \to Q$ are graph homomorphisms such that $pa = qb$. By the universal property of pushouts in $\ncat{Set}$, there exists a unique $\ell: V(H +_G K) \to V(Q)$ such that $\ell r = p$ and $ \ell s = q$, given by $\ell(x) = p(h)$ if $x = [h]$ or $\ell(x) = q(k)$ if $x = [k]$. If $xy \in E(H+_G K)$, then we want to show that $\ell(x)\ell(y) \in E(Q)$. But if $xy \in E(H+_G K)$, then there must be either an edge $hh' \in E(H)$ or $kk' \in E(K)$ that quotient to $xy$. Thus $p(h)p(h') \in E(Q)$ or $q(k)q(k') \in E(Q)$, but then $\ell(x)\ell(y) = \ell[h]\ell[h'] = p(h)p(h')$ or $\ell(x)\ell(y) = \ell[k]\ell[k'] = q(k)q(k')$. Thus $\ell$ is a graph homomorphism such that $\ell r = p$ and $\ell s = q$. Hence $H+_G K$ is a pushout in $\Gr$.

Now suppose that $r(h) = r(h')$ for $h,h' \in V(H)$. Then $[h] = [h']$ in $V(H)+_{V(G)} V(K)$. But by the discussion above, this implies that $h = h'$. Similarly, $s$ is a monomorphism.
\end{proof}

Given a graph $G$, a \textbf{subgraph} of $G$ is a graph $H$ with $V(H) \subseteq V(G)$. We say that $H$ is an \textbf{induced subgraph} if for all $x,y \in V(H)$, there is an edge $xy \in E(H)$ if and only $xy \in E(G)$. We write $H \subseteq G$ to mean that $H$ is a subgraph of $G$. We will sometimes also call a monomorphism $i : H \hookrightarrow G$ a subgraph, though we really mean the isomorphism equivalence class of $i$. If $G$ is a graph, then subgraphs of $G$ are in bijection with subobjects of $G$ in $\Gr$. 

Given a graph $G$ and subgraphs $H, K \subseteq G$, let $H \cap K$ denote the subgraph of $G$ with $V(H \cap K) = V(H) \cap V(K)$ and $E(H \cap K) = E(H) \cap E(K)$. We call $H \cap K$ the \textbf{intersection} of $H \subseteq G$ and $K \subseteq G$. There are monomorphisms $i_H : H \cap K \hookrightarrow H$ and $i_K : H \cap K \hookrightarrow K$.

Similarly we let $H \cup K$ denote the subgraph of $G$ with $V(H \cup K) = V(H) \cup V(K)$ and where there is an edge $xy \in E(H \cup K)$ if and only if $x,y \in V(H)$ and $xy \in E(H)$ or $x,y \in V(K)$ and $xy \in E(K)$. We note that there are monomorphisms $j_H : H \hookrightarrow H \cup K$ and $j_K : K \hookrightarrow H \cup K$.

\begin{Lemma}
Given a graph $G$ and subgraphs $H, K \subseteq G$, then the following square commutes and is a pullback
\begin{equation*}
    \begin{tikzcd}
	{H \cap K} & K \\
	H & G
	\arrow["{i_K}", hook', from=1-1, to=1-2]
	\arrow["{i_H}"', hook, from=1-1, to=2-1]
	\arrow["\lrcorner"{anchor=center, pos=0.125}, draw=none, from=1-1, to=2-2]
	\arrow[hook', from=1-2, to=2-2]
	\arrow[hook, from=2-1, to=2-2]
\end{tikzcd}
\end{equation*}
furthermore, the following square commutes and is a pushout
\begin{equation*}
\begin{tikzcd}
	{H \cap K} & K \\
	H & {H \cup K}
	\arrow["{i_K}", hook', from=1-1, to=1-2]
	\arrow["{i_H}"', hook, from=1-1, to=2-1]
	\arrow["{j_K}", hook', from=1-2, to=2-2]
	\arrow["{j_H}"', hook, from=2-1, to=2-2]
	\arrow["\lrcorner"{anchor=center, pos=0.125, rotate=180}, draw=none, from=2-2, to=1-1]
\end{tikzcd}    
\end{equation*}
\end{Lemma}
    
\begin{Lemma} \label{lem can replace pushouts by unions of subgraphs}
Suppose that the following commutative square
\begin{equation*}
\begin{tikzcd}
	G & K \\
	H & P
	\arrow["g", hook', from=1-1, to=1-2]
	\arrow["f"', hook, from=1-1, to=2-1]
	\arrow["s", hook', from=1-2, to=2-2]
	\arrow["r"', hook, from=2-1, to=2-2]
	\arrow["\lrcorner"{anchor=center, pos=0.125, rotate=180}, draw=none, from=2-2, to=1-1]
\end{tikzcd}    
\end{equation*}
is a pushout of monomorphisms in $\Gr$. Then the pullback diagram
\begin{equation*}
\begin{tikzcd}
	{H \cap K} & K \\
	H & P
	\arrow["{i_K}", hook', from=1-1, to=1-2]
	\arrow["{i_H}"', hook, from=1-1, to=2-1]
	\arrow["\lrcorner"{anchor=center, pos=0.125}, draw=none, from=1-1, to=2-2]
	\arrow["s", hook', from=1-2, to=2-2]
	\arrow["r"', hook, from=2-1, to=2-2]
\end{tikzcd}    
\end{equation*}
is also a pushout, $V(G) = V(H \cap K)$, and $G \subseteq H \cap K$.
\end{Lemma}

\begin{proof}
Suppose we have a commutative diagram of the form
\begin{equation} \label{eq pushout and intersection diagram in Gr}
\begin{tikzcd}
	&& K \\
	G & {H \cap K} && P & Q \\
	&& H
	\arrow["s"', hook', from=1-3, to=2-4]
	\arrow["b", curve={height=-6pt}, from=1-3, to=2-5]
	\arrow["g", curve={height=-12pt}, hook', from=2-1, to=1-3]
	\arrow["\ell", from=2-1, to=2-2]
	\arrow["f"', curve={height=12pt}, hook, from=2-1, to=3-3]
	\arrow["{{i_K}}", hook', from=2-2, to=1-3]
	\arrow["{{i_H}}"', hook, from=2-2, to=3-3]
	\arrow["h", dashed, from=2-4, to=2-5]
	\arrow["r", hook, from=3-3, to=2-4]
	\arrow["a"', curve={height=6pt}, from=3-3, to=2-5]
\end{tikzcd}
\end{equation}
we want to construct the dashed map $h$. First we note that the unique map $\ell : G \to H \cap K$ obtained from the universal property of the pullback must be a monomorphism since, for instance, $i_K \ell = g$ is a monomorphism. Hence, as subgraphs, $G \subseteq H \cap K$. Now since $ai_H = b i_K$, we have $af = bg$. Thus there exists a unique map $h : P \to Q$ such that $hr = a$ and $hs = b$. But this implies that $P$ is a pushout of $i_K$ and $i_H$, which is what we wanted to show.
Now clearly $V(G) \subseteq V(H \cap K)$. But by construction $V(P) = V(H) +_{V(G)} V(K)$. So the only way that a point $x \in V(P)$ can be considered as a vertex of the subgraph $G$ is if there exists a $z \in V(G)$ such that $x = [f(z)] = [g(z)]$. So if $x \in V(H) \cap V(K)$, then $x \in V(G)$.
\end{proof}

We furthermore note that pushouts along monomorphisms are pullback-stable. This observation is important for the proof of Proposition \ref{prop sd-category for treewidth}.

\begin{Lemma} \label{lem Gr has universal pushouts along monomorphisms}
Suppose that we have a diagram in $\ncat{Gr}$ of the form
\begin{equation*}
\begin{tikzcd}
	{G_0} && {G_1} \\
	& G && {G_2} \\
	{H_0} && {H_1} \\
	& H && {H_2}
	\arrow["r"', from=1-1, to=2-2]
	\arrow["{{f_0}}"', from=1-1, to=3-1]
	\arrow["{{g_0}}"', hook, from=1-3, to=1-1]
	\arrow["{{g_2}}", hook', from=1-3, to=2-4]
	\arrow["\lrcorner"{anchor=center, pos=0.125, rotate=-90}, draw=none, from=1-3, to=3-1]
	\arrow["{{f_1}}"'{pos=0.8}, from=1-3, to=3-3]
	\arrow["f"'{pos=0.2}, from=2-2, to=4-2]
	\arrow["s"{pos=0.4}, from=2-4, to=2-2]
	\arrow[""{name=0, anchor=center, inner sep=0}, "{{f_2}}"', from=2-4, to=4-4]
	\arrow["k"', from=3-1, to=4-2]
	\arrow["{{h_0}}"{pos=0.3}, hook, from=3-3, to=3-1]
	\arrow["{{h_2}}"'{pos=0.3}, hook', from=3-3, to=4-4]
	\arrow["\ell", from=4-4, to=4-2]
	\arrow["\lrcorner"{anchor=center, pos=0.125}, draw=none, from=1-3, to=0]
\end{tikzcd}
\end{equation*}
where $h_0, h_2$ are monomorphisms and the four vertical faces are all pullbacks. If $H$ is a pushout, then so is $G$.
\end{Lemma}

\begin{proof}
We first note that the functor $V : \Gr \to \ncat{Set}$ preserves finite limits and pushouts along spans of monomorphisms. Hence applying $V$ to the diagram, we know that $V(G) \cong V(G_0) +_{V(G_1)} V(G_2)$ since $\ncat{Set}$ has universal colimits. Now suppose we have a commutative diagram
\begin{equation*}
    \begin{tikzcd}
	{G_1} & {G_2} \\
	{G_0} & G \\
	&& Q
	\arrow["{g_2}", hook', from=1-1, to=1-2]
	\arrow["{g_0}"', hook, from=1-1, to=2-1]
	\arrow["s"', hook', from=1-2, to=2-2]
	\arrow["q", curve={height=-6pt}, from=1-2, to=3-3]
	\arrow["r", hook, from=2-1, to=2-2]
	\arrow["p"', curve={height=6pt}, from=2-1, to=3-3]
	\arrow["h", dashed, from=2-2, to=3-3]
\end{tikzcd}
\end{equation*}
we wish to define a unique dashed map $h$. Applying $V$, we obtain a unique map $V(h) : V(G) \to V(Q)$ by the universal property of the pushout in $\ncat{Set}$. We want to show that this is a graph homomorphism. Suppose that $xy \in E(G)$. Then $f(xy) \in E(H)$, and hence since $H$ is a pushout, there exists some $vv' \in E(H_0)$ or $ww' \in E(H_2)$ such that $f(xy) = [v][v']$ or $f(xy) = [w][w']$. But $G_0$ and $G_2$ are pushouts, hence there exists $x_0y_0 \in E(G_0)$ or $x_1 y_1 \in E(G_1)$ such that $f_0(x_0y_0) = vv'$ or $f_1(x_1y_1) = ww'$. Since $p$ and $q$ are graph homomorphisms, this implies that $h(xy) = p(x_0y_0) \in E(Q)$ or $h(xy) = q(x_1 y_1) \in E(Q)$. Hence $h$ is a graph homomorphism, and thus $G$ is a pushout.
\end{proof}

\begin{Rem}
Lemma \ref{lem Gr has universal pushouts along monomorphisms} is equivalent to saying that $\ncat{Gr}$ has universal pushouts alongs spans of monomorphisms. This is closely related to, but is more general than, the concept of adhesivity, see \cite{lack2004adhesive}\footnote{But notice there that adhesive categories in particular have universal pushouts along monomorphisms. The category $\ncat{Gr}$ does not even have all pushouts along monomorphisms. Both morphisms in the span must be monos in order for the pushout to exist.}.  
\end{Rem}

Now let us show that colimits of tree decompositions (Definition \ref{def tree decomposition}) always exist in $\ncat{Gr}$.

\begin{Def}
We say a category $\cat{C}$ is \textbf{connected} if it is non-empty, and for every pair of objects $c, c' \in \cat{C}$ there exists of a finite zig-zag of the form
\begin{equation*}
    \begin{tikzcd}
	c & {c_0} & {c_1} & \dots & {c_n } & c'
	\arrow[from=1-1, to=1-2]
	\arrow[from=1-3, to=1-2]
	\arrow[from=1-3, to=1-4]
	\arrow[from=1-5, to=1-4]
	\arrow[from=1-5, to=1-6]
\end{tikzcd}
\end{equation*}
or of the form where the direction of any of the arrows is reversed, in $\cat{C}$. 
\end{Def}

We now wish to define what it means for a category to be simply connected. To do this, we need to discuss the Gabriel-Zisman localization of a small category at the set of all of its morphisms. This construction first appeared in \cite{gabriel1967calculus}, but is now well-known. For an in-depth discussion see \cite{simpson2005explaining}. We will merely sketch its construction here.

\begin{Def}
Given a small category $\cat{C}$, let $\widetilde{\cat{C}}$ denote the category with the same objcts as $\cat{C}$ and whose morphisms consist of equivalence classes of finite zig-zags consisting of morphisms in $\cat{C}$ and $\cat{C}^\op$ under the smallest equivalence relation such that
\begin{itemize}
\itemsep -.5em
    \item adjacent arrows pointing in the same direction may be composed,
    \item adjacent pairs of the forms
    \begin{equation*}
        c \xleftarrow{f} d \xrightarrow{f^\op} c, \qquad  d \xrightarrow{f^\op} c \xleftarrow{f} d
    \end{equation*}
    are equivalent to identities, and
    \item identity arrows pointing either forwards or backwards may be removed.
\end{itemize}
Note that $\widetilde{\cat{C}}$ is a groupoid. We say that a small, connected category $\cat{C}$ is \textbf{simply connected} if $\widetilde{\cat{C}}$ has precisely one morphism between any two objects.
\end{Def}

In other words, a category is simply connected if it is connected and there are no ``non-trivial'' zig-zags from each object to itself. 

\begin{Ex}[{\cite[Page 733]{pare1990simply}}]
If a category $\cat{C}$ is finitely cofiltered, i.e. every finite diagram in $\cat{C}$ has a cone, then it is simply connected.
\end{Ex}

By a beautiful result of Par\'e \cite[Theorem 1]{pare1990simply}, limits of simply connected diagrams can be computed purely by pullbacks. However his proof is quite involved and much more powerful than what we need for our purposes, so we record our own simplified version of this result. We have taken the following proof from \cite{nlabconnectedlimit}.

\begin{Prop} \label{prop simply connected colimit computed by pushout}
Let $\cat{C}$ be a category with pushouts. If $d : I \to \cat{C}$ is a diagram where $I$ is a finite, simply connected category, then the colimit of $d$ exists and is given by an iterated pushout.
\end{Prop}

\begin{proof}
Since $I$ is non-empty, let us fix an object $x_0 \in I$. For any $y \in I$ we define $\ell(y)$ to be the number of morphisms of any minimal zig-zag between $x_0$ and $y$. Since $I$ is finite and simply connected, this is well-defined. Now choose a linear order of the objects of $I$ as $x_0, x_1, \dots, x_n$ such that $\ell(x_i) \leq \ell(x_{i+1})$ for $0 \leq i \leq n-1$.

Now let us inductively define objects $P_i$ in $\cat{C}$ equipped with maps $p_{ij} : d(x_j) \to P_i$ for all $0 \leq j \leq i \leq n$ as follows.

For the base case, let $P_0 = d(x_0)$ and let $p_{00} : d(x_0) \to d(x_0)$ be the identity map.

Now for the inductive step, suppose that we have objects $P_i$ for $0 \leq i \leq k$ and maps $p_{ij} : d(x_j) \to P_i$ for all $0 \leq j \leq i \leq k$. Now choose a zig-zag of minimal length between $x_0$ and $x_{k+1}$. It will look like $$x_0 \leftrightarrow y_1 \leftrightarrow y_2 \leftrightarrow \dots \leftrightarrow y_{\ell(x_{k+1}) - 1} \leftrightarrow x_{k+1}.$$
But since we have ordered all of the objects of $I$, we know that each $y_r = x_j$ for some $j \leq k$ and so we have a morphism $p_{ij} : y_r \to P_i$ for every $j \leq i$ and $r \leq \ell(x_{k+1})-1$. Let us write $y_\ell = y_{\ell(x_{k+1})-1} = x_j$ for the final object. Suppose that the final map in the zig-zag is a map directed as $f : y_\ell \leftarrow x_{k+1}$. Then let $P_{k+1} = P_k$ and let $p_{(k+1)(k+1)} : d(x_{k+1}) \to P_{k+1}$ be the composite map
\begin{equation*}
   d(x_{k+1}) \xrightarrow{d(f)} d(y_\ell) \xrightarrow{p_{kj}} P_{k+1} = P_k 
\end{equation*}
and set each other $p_{(k+1)r} :d(x_r) \to P_{k+1}$ equal to $p_{kr}$, which we assumed has already been defined for all $r \leq k$. 

If the final map in the zig-zag is instead a map directed as $f : y_\ell \to x_{k+1}$, then let $P_{k+1}$ be the pushout
\begin{equation*}
\begin{tikzcd}
	{d(y_\ell)} & {P_k} \\
	{d(x_{k+1})} & {P_{k+1}}
	\arrow["{{p_{kj}}}", from=1-1, to=1-2]
	\arrow["{{d(f)}}"', from=1-1, to=2-1]
	\arrow[from=1-2, to=2-2]
	\arrow["{{p_{(k+1)(k+1)}}}"', from=2-1, to=2-2]
	\arrow["\lrcorner"{anchor=center, pos=0.125, rotate=180}, draw=none, from=2-2, to=1-1]
\end{tikzcd}
\end{equation*}
and this defines the map $p_{(k+1)(k+1)}$. The other maps $p_{(k+1)i}$ are defined with composition of $d(f)$ with $p_{ki}$.

Thus we have defined an object $P_n$ in $\cat{C}$ along with a cocone $p$ over $d$ given by the maps $p_{ni} :d(x_i) \to  P_n$ for $0 \leq i \leq n$, and let $p_k = p_{kk}$. Let us show that $p$ is a colimit cocone over $d$. Suppose there was an object $Q$ and a cocone $q :d \Rightarrow \Delta_Q$ with components $q_i : d(x_i) \to Q$. Let us prove that it must factor uniquely through the $p$ by induction. In the base case, $q_0$ clearly factors uniquely through $p_0 : P_0 \to d(x_0)$ since it is the identity map. Now suppose that each $q$ factors uniquely through $p$ restricted to $x_0, \dots, x_k$. Then in the inductive step, either $P_{k+1}$ is the same as $P_k$ and therefore $q_{k+1}$ factors uniquely through $p_{k+1}$ because $q_k$ factors uniquely through $p_k$ by assumption, or $P_{k+1}$ is a pushout, and therefore $q_{k+1}$ factors uniquely through $p_{k+1}$ by the universal property of the pullback. Therefore $p :d \Rightarrow \Delta_{P_n}$ is a colimit cocone.
\end{proof}

\begin{Ex} \label{ex colimit of a zig-zag}
Computing the colimit of a diagram $d : \smallint T \to \Gr$ where $T$ is a tree with three vertices $x,y,z$ and edges $e = xy$, $f = yz$
\begin{equation}
\begin{tikzcd}
	{d(x)} & {d(e)} & {d(y)} & {d(f)} & {d(z)} \\
	{d(x)} & {d(x)} & {P_2} & {P_2} & {P_4} \\
	{P_0} & {P_1} & {P_2} & {P_3} & {\text{colim} \, d \cong G}
	\arrow[Rightarrow, no head, from=1-1, to=2-1]
	\arrow[hook', from=1-2, to=1-1]
	\arrow[hook, from=1-2, to=1-3]
	\arrow[hook, from=1-2, to=2-2]
	\arrow[hook, from=1-3, to=2-3]
	\arrow[hook', from=1-4, to=1-3]
	\arrow[hook, from=1-4, to=1-5]
	\arrow[hook, from=1-4, to=2-4]
	\arrow[hook, from=1-5, to=2-5]
	\arrow[Rightarrow, no head, from=2-1, to=2-2]
	\arrow[Rightarrow, no head, from=2-1, to=3-1]
	\arrow[hook, from=2-2, to=2-3]
	\arrow[Rightarrow, no head, from=2-2, to=3-2]
	\arrow["\lrcorner"{anchor=center, pos=0.125, rotate=180}, draw=none, from=2-3, to=1-2]
	\arrow[Rightarrow, no head, from=2-3, to=2-4]
	\arrow[Rightarrow, no head, from=2-3, to=3-3]
	\arrow[hook, from=2-4, to=2-5]
	\arrow[Rightarrow, no head, from=2-4, to=3-4]
	\arrow["\lrcorner"{anchor=center, pos=0.125, rotate=180}, draw=none, from=2-5, to=1-4]
	\arrow[Rightarrow, no head, from=2-5, to=3-5]
	\arrow[hook, from=3-1, to=3-2]
	\arrow[hook, from=3-2, to=3-3]
	\arrow[hook, from=3-3, to=3-4]
	\arrow[hook, from=3-4, to=3-5]
\end{tikzcd}
\end{equation}
\end{Ex}

\begin{Cor} \label{cor colimits of tree decomps exist in Gr}
If $d : \smallint T \to \ncat{Gr}$ is a tree-shaped structured decomposition, then the colimit of $d$ exists. Furthermore, the colimit cocone maps are monomorphisms.
\end{Cor}

\begin{proof}
This follows from Lemma \ref{lem pushouts of monos exist in Gr} and Proposition \ref{prop simply connected colimit computed by pushout}.
\end{proof}

\begin{Cor} \label{cor can pullback tree decompositions}
Let $d : \smallint T \to \ncat{Gr}$ be a tree-shaped structured decomposition of a graph $H$, and suppose there is a map $f : G \to H$ of graphs. Let $f^*(d) : \smallint T \to \ncat{Gr}$ denote the diagram obtained by pulling back the colimit cocone maps along $f$. Then $f^*(d)$ is a structured decomposition of $G$.
\end{Cor}

\begin{proof}
This follows from Lemma \ref{lem Gr has universal pushouts along monomorphisms} and Proposition \ref{prop simply connected colimit computed by pushout}.
\end{proof}

In the language of Remark \ref{rem universal colimits}, Corollary \ref{cor can pullback tree decompositions} proves that $\Gr$ has universal $\text{Trees}$-colimits.

\subsection{Multigraphs}

This section is used primarily in Section \ref{section caremsin's graph decompositions}. We first collect a few facts about $\Gr_\mlt$, the category of multigraphs, and then show that $\Gr_\mlt$ has pullback-stable or universal finite colimits.

\begin{Lemma}[{\cite[Theorem 3.6.1, 3.6.2]{plessas2011categories}}] \label{lem monos and epis are very nice in cat of multigraphs}
Given a morphism $f: G \to H$ in $\Gr_\mlt$, the following are equivalent:
\begin{itemize}
\itemsep -.5em
    \item $f$ is a mono/epimorphism,
    \item $f$ is a regular mono/epimorphism,
    \item $f$ is an effective mono/epimorphism,
    \item $f$ is an extremal mono/epimorphism.
\end{itemize}
\end{Lemma}

Since $\Gr_\mlt$ is finitely complete, thanks to Lemma \ref{lem monos and epis are very nice in cat of multigraphs}, we need only to show that epimorphisms of multigraphs are stable under pullback in order to show that $\Gr_\mlt$ is regular.

\begin{Lemma}
Epimorphisms in $\Gr_\mlt$ are stable under pullback. 
\end{Lemma}

\begin{proof}
Given maps $f : K \to H$ and $g : G \to H$, the pullback $K \times_H G$ is the multigraph with $V(K \times_H G) = V(K) \times_{V(H)} V(G)$ and where there is an edge $e$ between $(k,g)$ and $(k', g')$ for every pair of edges $e'$ between $k$ and $k'$ and $e''$ between $g$ and $g'$ such that $f(e') = g(e'')$. Therefore, if $g$ is an epimorphism of multigraphs, then by Lemma \ref{lem monos and epis of graphs}, $E(g)$ and $V(g)$ are surjective. So if we let $f^*(g) : K \times_H G \to K$ denote the induced map from the pullback, then $E(f^*(g))$ and $V(f^*(g))$ will be surjective as well, so again by Lemma \ref{lem monos and epis of graphs}, $f^*(g)$ is an epimorphism of multigraphs.
\end{proof}

\begin{Cor}
The category $\Gr_\mlt$ is regular.
\end{Cor}

We recall (Definition \ref{def multigraphs}) that a multigraph $G$ consists of finite sets $V(G)$, $E(G)$ and a function
\begin{equation*}
    \partial : E(G) \to V(G)^{(2)}.
\end{equation*}
We note that $(-)^{(2)}$ itself forms a functor $(-)^{(2)} : \ncat{FinSet} \to \ncat{FinSet}$, and so we can write $\Gr_\mlt$ as the comma category $(1_{\ncat{FinSet}}, (-)^{(2)})$. It is useful to know the following result on computing (co)limits in comma categories.

\begin{Lemma}[{\cite[Theorem 5.2.3]{rydeheard1988computational}}] \label{lem computing (co)limits in comma categories}
Let $F: \cat{C} \to \cat{E}$ and $G: \cat{D} \to \cat{E}$ be functors, and let $(F \downarrow G)$ denote the corresponding comma category, with projection functors $\pi : (F \downarrow G) \to \cat{C}$, $\pi' : (F \downarrow G) \to \cat{D}$. If $d : I \to (F \downarrow G)$ is a diagram such that the colimits $U = \colim \, \pi d$ and $V = \colim \, \pi' d$ exist in $\cat{C}$ and $\cat{D}$ respectively, and $F$ preserves this colimit, then the colimit of $d$ exists in $(F \downarrow G)$ and is given by the unique map
\begin{equation*}
    \colim \, F\pi d \cong F(U) \to G(V)
\end{equation*}
in $\cat{E}$ that commutes with the colimit cocone maps. Conversely if $G$ preserves the limit of $\pi' d$ then the limit of $d$ exists in $(F \downarrow G)$ and is given by the corresponding unique map
\begin{equation*}
    F(U) \to G(V) \cong \lim G\pi' d.
\end{equation*}
\end{Lemma}

Thanks to Lemma \ref{lem computing (co)limits in comma categories}, we know how to compute colimits in $\Gr_\mlt$.

\begin{Lemma}
Let $d : I \to \Gr_\mlt$ be a finite diagram. Then the colimit $\colim \, d$ exists and is given by the multigraph
\begin{equation*}
    \partial : \colim \, E(d(i)) \to \left( \colim \, V(d(i)) \right)^{(2)},
\end{equation*}
where $\partial$ is the unique induced map from the colimit cocone maps.
\end{Lemma}

The category $\Gr_\mlt$ has both an initial object $\varnothing$ given by the empty graph, and a terminal object $*$ given by the multigraph with one vertex and one loop. Given multigraphs $G$ and $H$, their product $G \times H$ is the multigraph with $V(G \times H) = V(G) \times V(H)$, $E(G \times H) = E(G) \times E(H)$, and if $e \in E(G)$ with $\partial_G(e) = \{x,x'\}$ and $e' \in E(H)$ with $\partial_H(e') = \{y, y'\}$, then $\partial_{G \times H}(e,e') = \{(x,y), (x', y') \}$. Pullbacks are similarly easy to describe. If $f : G \to H$ and $g : K \to H$ are maps of multigraphs, then $G \times_H K$ is the multigraph with $V(G \times_H K) = V(G) \times_{V(H)} V(K)$, $E(G \times_{H} K) = E(G) \times_{E(H)} E(K)$ and where if $e \in E(G)$ with $\partial_G(e) = \{x,x'\}$ and $e' \in E(K)$ with $\partial_K(e') = \{y, y'\}$, then $\partial_{G \times_H K}(e,e') = \{(x,y), (x', y') \}$. If $f : G \to H$ is a map of multigraphs, and $H' \hookrightarrow H$ is a monomorphism, i.e. a subgraph inclusion, then the pullback $f^{-1}(H')$ is given by the multigraph with vertex set $V(f)^{-1}(V(H'))$ and edge set $E(f)^{-1}(E(H'))$.

\begin{Lemma}
Coproducts in $\Gr_\mlt$ are stable under pullbacks.
\end{Lemma}

\begin{proof}
If $f : G \to H$ is a map of multigraphs, where $H \cong H_0 + H_1$, then we have $G \cong f^{-1}(H_0) + f^{-1}(H_1)$.
\end{proof}

\begin{Lemma}
Coequalizers in $\Gr_\mlt$ are stable under pullbacks.
\end{Lemma}

\begin{proof}
This follows from coequalizers being stable under pullback in $\ncat{Set}$. Indeed, suppose that $f, g: G \to H$ are maps of multigraphs, with colimit $X = \text{coeq}(f,g)$. If $h : Y \to X$ is a map of multigraphs, then we obtain a pair of pullback diagrams
\begin{equation*}
    \begin{tikzcd}
	{Y \times_X G} & G \\
	{Y\times_X H} & H \\
	Y & X
	\arrow[from=1-1, to=1-2]
	\arrow["{h^*(f)}"', shift right, from=1-1, to=2-1]
	\arrow["{h^*(g)}", shift left, from=1-1, to=2-1]
	\arrow["\lrcorner"{anchor=center, pos=0.125}, draw=none, from=1-1, to=2-2]
	\arrow["g", shift left, from=1-2, to=2-2]
	\arrow["f"', shift right, from=1-2, to=2-2]
	\arrow[from=2-1, to=2-2]
	\arrow["{q_Y}"', from=2-1, to=3-1]
	\arrow["\lrcorner"{anchor=center, pos=0.125}, draw=none, from=2-1, to=3-2]
	\arrow["{q_X}", from=2-2, to=3-2]
	\arrow["h", from=3-1, to=3-2]
\end{tikzcd}
\end{equation*}
So we have $V(q_Y) : V(Y \times_X H) \to V(Y)$ is the coequalizer of $V(h^*(f))$ and $V(h^*(g))$, and the same holds for the edges. Its then easy to see this causes $Y$ to be a coequalizer of $h^*(f)$ and $h^*(g)$.
\end{proof}

\begin{Cor} \label{cor finite colimits of multigraphs are pullback-stable}
Finite colimits in $\Gr_\mlt$ are stable under pullbacks. 
\end{Cor}

\subsection{Hypergraphs} \label{section hypergraphs}

This section is used primarily in Section \ref{section hypergraph treewidth}. We prove several facts about the category $\Hyp$ of hypergraphs (Definition \ref{def hypergraph}).

The vertex functor $V = V_\Hyp : \Hyp \to \ncat{FinSet}$ has both a left adjoint $\text{Disc}: \ncat{FinSet} \to \Hyp$ given by the discrete hypergraph on a set, and right adjoint $\text{CoDisc}: \ncat{FinSet} \to \Hyp$ given by taking every nonempty subset of a set $S$ to be a hyperedge edge. Furthermore $V$ is a faithful functor, making $\Hyp$ into a concrete category. Hence if $V(f)$ is inj/surjective, then $f: H \to H'$ is a mono/epimorphism. Since $V$ has a left and right adjoint, a map $f : H \to H'$ in $\Hyp$ is a mono/epimorphism if and only if $V(f)$ is inj/surjective.

Given a hypergraph $H$, a \textbf{subhypergraph} $H' \subseteq H$ consists of subsets $V(H') \subseteq V(H)$ and $E(H') \subseteq E(H)$ such that if $e \in E(H')$ and $v \in e$, then $v \in V(H')$. Subhypergraphs of a hypergraph $H$ are in bijection with subobjects of $H$ in $\Hyp$. Given a subset of vertices $S \subseteq V(H)$ of a hypergraph $H$, the \textbf{induced subhypergraph} $H[S]$ is the subhypergraph of $H$ containing all of the vertices of $S$ and all the hyperedges $e$ in $H$ such that $e \subseteq S$.

The intersection $H \cap K$ of subhypergraphs of $H \subseteq H'$ and $K \subseteq H'$ of a hypergraph $H'$ is the subhypergraph $H \cap K \subseteq H'$ with $V(H \cap K) = V(H) \cap V(K)$, and $E(H \cap K) = E(H) \cap E(K)$. There are monomorphisms $i_H : H \cap K \hookrightarrow H$ and $i_K : H \cap K \hookrightarrow K$.

The empty hypergraph $\varnothing$ is clearly an initial object in $\Hyp$, but there is no terminal object.

\begin{Def} \label{def pushout of span of monos for hypergraphs}
Suppose that we have a span of monomorphisms $i : G \to H$ and $j : G \to K$ of hypergraphs. Let $H +_G K$ denote the hypergraph defined as follows. Let $V(H +_G K) = V(H) +_{V(G)} V(K)$, and if $e = \{x_1, \dots, x_n \}$ is a nonempty subset of $V(H+_G K)$, then it is a hyperedge of $H+_G K$ if and only if there exists a hyperedge $e' = \{v_1, \dots, v_n\}$ in $H$ or $K$ such that $x_i = [v_i]$ for all $1 \leq i \leq n$. By the same reasoning as in Definition \ref{def pushouts of spans of monos for graphs}, since $i$ and $j$ are monomorphisms, a vertex $x \in H+_G K$ belongs either to $H$ or $K$ or if it belongs to both of them, then there is a unique vertex $x' \in G$ such that $i(x') = j(x') = x$. It is not hard to see that $H+_G K$ is a well-defined hypergraph.
\end{Def}

Let $r : H \to H+_G K$ and $s : K \to H+_G K$ be defined simply as the quotient map on vertices. They are easily seen to be hypergraph homomorphisms.

The proofs of the following four results are practically identical to the proofs of Lemma \ref{lem pushouts of monos exist in Gr}, Lemma \ref{lem Gr has universal pushouts along monomorphisms}, Corollary \ref{cor colimits of tree decomps exist in Gr} and Corollary \ref{cor can pullback tree decompositions}.

\begin{Lemma} \label{lem pushouts of monos in Hyp}
Given a span of monomorphisms in $\Hyp$
\begin{equation*}
  \begin{tikzcd}
	G & K \\
	H & {P = H +_G K}
	\arrow["j", hook', from=1-1, to=1-2]
	\arrow["i"', hook, from=1-1, to=2-1]
	\arrow["s", hook', from=1-2, to=2-2]
	\arrow["r"', hook, from=2-1, to=2-2]
	\arrow["\lrcorner"{anchor=center, pos=0.125, rotate=180}, draw=none, from=2-2, to=1-1]
\end{tikzcd}     
\end{equation*}
their pushout exists, is given by $P$, and furthermore the maps $r$ and $s$ are monomorphisms.
\end{Lemma}

\begin{Lemma} \label{lem hypergraphs have pullback stable pushouts of monos}
Suppose that we have a diagram in $\ncat{Hyp}$ of the form
\begin{equation*}
\begin{tikzcd}
	{G_0} && {G_1} \\
	& G && {G_2} \\
	{H_0} && {H_1} \\
	& H && {H_2}
	\arrow["r"', from=1-1, to=2-2]
	\arrow["{{f_0}}"', from=1-1, to=3-1]
	\arrow["{{g_0}}"', hook, from=1-3, to=1-1]
	\arrow["{{g_2}}", hook', from=1-3, to=2-4]
	\arrow["\lrcorner"{anchor=center, pos=0.125, rotate=-90}, draw=none, from=1-3, to=3-1]
	\arrow["{{f_1}}"'{pos=0.8}, from=1-3, to=3-3]
	\arrow["f"'{pos=0.2}, from=2-2, to=4-2]
	\arrow["s"{pos=0.4}, from=2-4, to=2-2]
	\arrow[""{name=0, anchor=center, inner sep=0}, "{{f_2}}"', from=2-4, to=4-4]
	\arrow["k"', from=3-1, to=4-2]
	\arrow["{{h_0}}"{pos=0.3}, hook, from=3-3, to=3-1]
	\arrow["{{h_2}}"'{pos=0.3}, hook', from=3-3, to=4-4]
	\arrow["\ell", from=4-4, to=4-2]
	\arrow["\lrcorner"{anchor=center, pos=0.125}, draw=none, from=1-3, to=0]
\end{tikzcd}
\end{equation*}
where $h_0, h_2$ are monomorphisms and the four vertical faces are all pullbacks. If $H$ is a pushout, then so is $G$.
\end{Lemma}

\begin{Cor} \label{cor colimits of hypergraph tree decomps exist in Hyp}
If $d : \smallint T \to \Hyp$ is a tree-shaped structured decomposition, then the colimit of $d$ exists. Furthermore, the colimit cocone maps are monomorphisms.
\end{Cor}

\begin{Cor} \label{cor can pullback hypergraph tree decompositions}
Let $d : \smallint T \to \Hyp$ be a tree-shaped structured decomposition of a hypergraph $H'$, and suppose there is a map $f : H \to H'$ of hypergraphs. Let $f^*(d) : \smallint T \to \Hyp$ denote the diagram obtained by pulling back the colimit cocone maps along $f$. Then $f^*(d)$ is a structured decomposition of $H$.
\end{Cor}

By construction the inclusion $\iota : \Gr \hookrightarrow \Hyp$ preserves finite limits, coproducts and pushouts of spans of monomorphisms.

\begin{Lemma} \label{lem gaifman graph functor preserves pushouts of monos}
The Gaifman graph functor $\Gaif : \Hyp \to \Gr$ preserves pushouts of spans of monomorphisms.
\end{Lemma}

\begin{proof}
Suppose we have a pushout of a span of monomorphisms in $\Hyp$
\begin{equation*}
    \begin{tikzcd}
	G & K \\
	H & P
	\arrow["j", hook', from=1-1, to=1-2]
	\arrow["i"', hook, from=1-1, to=2-1]
	\arrow["s", hook', from=1-2, to=2-2]
	\arrow["r"', hook, from=2-1, to=2-2]
	\arrow["\lrcorner"{anchor=center, pos=0.125, rotate=180}, draw=none, from=2-2, to=1-1]
\end{tikzcd}
\end{equation*}
and suppose that $p : \Gaif(H) \to Q$ and $q : \Gaif(K) \to Q$ are graph morphisms such that $p \Gaif(i) = q \Gaif(j)$. We want to obtain a unique graph morphism $h$ such that $h \Gaif(r) = p$ and $h \Gaif(s) = q$.
\begin{equation*}
\begin{tikzcd}
	{\Gaif(G)} & {\Gaif(K)} \\
	{\Gaif(H)} & {\Gaif(P)} \\
	&& Q
	\arrow["{{\Gaif(j)}}", hook', from=1-1, to=1-2]
	\arrow["{{\Gaif(i)}}"', hook, from=1-1, to=2-1]
	\arrow["{{\Gaif(s)}}", hook', from=1-2, to=2-2]
	\arrow["q", curve={height=-18pt}, from=1-2, to=3-3]
	\arrow["{{\Gaif(r)}}"', hook, from=2-1, to=2-2]
	\arrow["p"', curve={height=18pt}, from=2-1, to=3-3]
	\arrow["\lrcorner"{anchor=center, pos=0.125, rotate=180}, draw=none, from=2-2, to=1-1]
	\arrow["h", dashed, from=2-2, to=3-3]
\end{tikzcd}
\end{equation*}
Applying $V_\Hyp$ to the above commutative diagram, we obtain a pushout diagram in $\ncat{Set}$, since $V_\Hyp \circ \Gaif = V_\Gr$. So we get a unique map $V(h) : V(\Gaif(P)) = V(P) \to V(Q)$. We need only to show that it is a graph map. Suppose that $xy \in E(\Gaif(P))$. Then there exists an $n$-edge $e \in E(P)$ such that $x,y \in e$. Thus there is either an $n$-edge $e' \in E(H)$ or $e'' \in E(K)$ such that $e = [e']$ or $e = [e'']$. Thus there is either a $vw \in E(\Gaif(H))$ or $ab \in E(\Gaif(K))$ such that $v,w \in e'$ or $a,b \in e''$ and $x = [v]$, $y = [w]$ or $x = [a], y = [b]$. Hence $h(xy) = p(v)p(w)$ or $q(a)q(b)$, which in both cases are edges. Hence $h$ is a graph map.
\end{proof}

\subsubsection{RS and sd-Tree Decompositions} \label{section appendix tree decompositions}

Our goal for this section is to prove Proposition \ref{prop hypergraph treewidth equal to gamma width}, which says that if $\Gamma = (\Hyp, \text{Trees}, \{K^n_\Hyp \}_{n \geq 0})$ is the width category from Proposition \ref{prop sd-category for hypergraph treewidth} and $H \in \Hyp$, then $\mathbf{tw}(H) = \mathbf{w}_\Gamma(H)$. To do this, we need to get a better handle on the connection between structured decompositions and tree decompositions.

Recall Definition \ref{def hypergraph treewidth} for the notion of a (Robertson-Seymour) or RS-tree decomposition $(X, T)$ of a hypergraph $H$. We want to compare this kind of decomposition with the notion of a structured decomposition $d : \smallint T \to \Hyp$ of $H$ where $T$ is a tree. Let us call the first notion an RS-tree decomposition and the latter an sd-tree decomposition of $H$.

We shall define functions
\begin{equation} \label{eq maps between RS and sd-tree decomps}
\begin{tikzcd}
	{\{\text{sd-tree decompositions of }H\}} && {\{\text{RS-tree decompositions of }H\}}
	\arrow["\Psi", curve={height=-18pt}, from=1-1, to=1-3]
	\arrow["\Phi", curve={height=-18pt}, from=1-3, to=1-1]
\end{tikzcd}   
\end{equation} 
as follows.

First let us define $\Phi$. Suppose that $(X,T)$ is an RS-tree decomposition of a hypergraph $H$. Then consider the diagram $\Phi(X,T) = d_{(X,T)} : \smallint T \to \Hyp$ given by mapping each vertex $t \in T$ to the induced subhypergraph $H[X(t)]$, and to each edge $tt'$ of $T$ the adhesion $H[X(t) \cap X(t')]$, and the corresponding maps being the inclusions.

\begin{Lemma} \label{lem RS-tree hypergraph decomp to sd-tree decomp}
If $(X,T)$ is an RS-tree decomposition of a hypergraph $H$, then $\Phi(X,T) = d_{(X,T)}$ is an sd-tree decomposition of $H$.
\end{Lemma}

\begin{proof}
The inclusions $\iota_t: H[X(t)] \hookrightarrow H$ define a cocone $\iota : d_{(X,T)} \Rightarrow \Delta(H)$. Let us show that this is a colimit cocone. If $\sigma : d_{(X,T)} \Rightarrow \Delta(H')$ is another cocone, we want to construct a unique map $h : H \to H'$ such that $\Delta(h) \iota = \sigma$. On vertices, let us define $V(h) : V(H) \to V(H')$ by setting $V(h)(x) = \sigma_t(x)$, for $x \in H$. Let us show that this is well-defined. First, since $(X,T)$ is an RS-tree decomposition, every vertex $x \in H$ belongs to some $H[X(t)]$, by condition (1). If $x$ belongs to $H[X(t)]$ and $H[X(t')]$, then by condition (3), it belongs to every bag $X(t_0)$ with $t_0$ lying on the unique path between $t$ and $t'$. Since $\sigma$ is a cocone, this means that for every such $t_0$, $\sigma_t(x) = \sigma_{t_0}(x) = \sigma_{t'}(x)$. Thus $V(h)$ is well defined.

Now if $e \in E_n(H)$, then there exists some bag $X(t)$ such that $e \in E_n(H[X(t)])$ by condition (2). Thus $h(e) \in E_n(H')$, since each $\sigma_t$ is a hypergraph homomorphism. Thus $h : H \to H'$ is a well-defined hypergraph homomorphism. Thus we need now only show that it the unique such map such that $\Delta(h) \iota = \sigma$.

If $h' : H \to H'$ is a hypergraph homomorphism such that $\Delta(h') \iota = \sigma$, then for every $x \in H$, by condition (1), there exists a $t \in T$ such that $x$ is in the image of $\iota_t : H[X(t)] \to H$. Thus $h'(x) = \sigma_t(x)$, and therefore $h$ is unique, so $\iota$ is a colimit cocone.
\end{proof}

\begin{Rem} \label{rem sd-tree to RS-tree graph decomp}
When $H$ is a graph, then an sd-tree decomposition $d : \smallint T \to \Hyp$ of $H$ factors uniquely through $\Gr$
\begin{equation*}
   \begin{tikzcd}[ampersand replacement=\&]
	{\smallint T} \&\& \Hyp \\
	\& \Gr
	\arrow["d", from=1-1, to=1-3]
	\arrow["{\tilde{d}}"', from=1-1, to=2-2]
	\arrow["\iota"', hook, from=2-2, to=1-3]
\end{tikzcd} 
\end{equation*}
simply because if $H'$ is a hypergraph that is not a graph, i.e. it contains an edge of cardinality not equal to $2$, then there exists no hypergraph homomorphism $H' \to H$. So the map $\Phi$, when restricted to $H$ a graph, lands in the set of sd-tree decompositions of $H$ in $\Gr$.
\end{Rem}

Now let us define the map $\Psi$. Suppose that $d : \smallint T \to \Hyp$ is an sd-tree decomposition of a hypergraph $H$. Let $\Psi(d) = (X_d, T)$ be defined as follows. Define $X_d : V(T) \to P(V(H))$ as $X_d(t) = V(d(t))$ for each $t \in T$.

\begin{Lemma} \label{lem sd-tree decomp to RS-tree hypergraph decomp}
If $d : \smallint T \to \Hyp$ is an sd-tree decomposition of a hypergraph $H$, then $\Psi(d) = (X_d, T)$ is an RS-tree decomposition of $H$.
\end{Lemma}

\begin{proof}
By Lemma \ref{lem pushouts of monos in Hyp} and Lemma \ref{prop simply connected colimit computed by pushout}, the colimit cocone maps $\lambda_t : d(t) \to H$ are all monomorphisms. Hence we can think of the bags and adhesions of $d$ as subhypergraphs of $H$. It is clear that conditions (1) and (2) of Definition \ref{def tree decomposition} hold for $\Phi(d) = (X_d, T)$, as every vertex and hyperedge of $H$ must belong to some bag in order for them to exist in the colimit. Now we need only to show that condition (3) holds.

Let $x, y, z \in T$, we want to show that if the unique path between $x$ and $z$ contains $y$ then $X_d(x) \cap X_d(z) \subseteq X_d(y)$. Let $T_0$ denote the subgraph of $T$ consisting of the unique path from $x$ to $z$, and let $d_0$ denote the restriction of $d$ to $T_0$. The colimit of $d_0$ will still be given by an iterated pushout, with $x$ and $z$ the endpoints of a zig-zag diagram in $\Hyp$ by Proposition \ref{prop simply connected colimit computed by pushout}.

Suppose that $v \in d(x) \cap d(z)$. If $v$ belongs to both $d(x)$ and $d(z)$, we want to show that it belongs to each intermediate bag of the zig-zag $d_0$. Let us prove this by induction on the length of the path $T_0$. For the base case, suppose that $T_0$ consists of the three vertices $x$, $y$, $z$ and two edges $e = xy$ and $f = yz$. Then the diagram in Example \ref{ex colimit of a zig-zag} is precisely the colimit of $d_0$, interpreted in $\Hyp$. We see that if $v \in d(x)$, then $v \in P_2$. But in order for $v$ to belong to both $P_2$ and $d(z)$, by the set-theoretic construction of the pushout (Definition \ref{def pushout of span of monos for hypergraphs}), it must be the case that $v \in d(f)$. But $d(f) \subseteq d(y)$, hence $v \in d(y)$.

For the inductive step, suppose now that $T_0$ is a path of length $n$ between $x$ and $z$ passing through $y$, and suppose further that for every path of length less than $n$ the intersection of the bags of the endpoints must belong to every intermediate bag. Then let $d(w)$ denote the bag immediately to the left of $d(z)$, and let $g = wz$ denote the rightmost edge in $T_0$. Let $T_1$ denote the unique path between $x$ and $w$, and let $d_1$ denote the restriction of $d$ to $T_1$. Suppose that $P$ is constructed as an iterated pushout of the zig-zag $d_1$ as in Example \ref{ex colimit of a zig-zag}. By assumption we know that $d(x) \cap d(w) \subseteq d(y)$, and we have the pushout
\begin{equation*}
\begin{tikzcd}
	\dots & {d(w)} & {d(g)} & {d(z)} \\
	\dots & \dots & P & H
	\arrow[hook', from=1-1, to=1-2]
	\arrow[hook, from=1-1, to=2-1]
	\arrow[hook, from=1-2, to=2-2]
	\arrow[hook, from=1-3, to=1-2]
	\arrow[hook', from=1-3, to=1-4]
	\arrow[hook, from=1-3, to=2-3]
	\arrow[hook', from=1-4, to=2-4]
	\arrow[hook, from=2-1, to=2-2]
	\arrow[hook, from=2-2, to=2-3]
	\arrow[hook, from=2-3, to=2-4]
	\arrow["\lrcorner"{anchor=center, pos=0.125, rotate=180}, draw=none, from=2-4, to=1-3]
\end{tikzcd}
\end{equation*}

But if $v \in d(x) \cap d(z)$, then $v \in P$. But then as in the base step, we know that since $v \in P \cap d(z)$, then $v \in d(g)$. But $d(g) \subseteq d(w)$, so $v \in d(w)$. Thus $v \in d(x) \cap d(w)$, and so $v \in d(y)$. Therefore $\Phi(d) = (X_d, T)$ satisfies condition (3), and is therefore an RS-tree decomposition of $H$.
\end{proof}

\begin{Rem}\label{rem RS-tree decomp to sd-tree decomp for graphs}
If $H$ is a graph and $d : \smallint T \to \Hyp$ an sd-tree decomposition, then as in Remark \ref{rem sd-tree to RS-tree graph decomp}, $d$ must factor uniquely through $\Gr$. Hence the map $\Psi$ restricted to sd-tree deceompositions in $\Gr$ lands in the set of RS-tree decompositions of $H$ as in Definition \ref{def tree decomposition}.
\end{Rem}

Thus we have constructed the functions of (\ref{eq maps between RS and sd-tree decomps}). It is clear that if $(X,T)$ is an RS-tree decomposition, then $(\Psi \circ \Phi)(X, T) = (X,T)$. Thus $\Psi$ is surjective.

We also note that since $\Phi$ and $\Psi$ do not change the bags of the decompositions, if $(X,T)$ is an RS-tree decomposition of a hypergraph $H$, then $w(X,T) = w(\Phi(X,T))$, where the former width is from Definition \ref{def hypergraph treewidth} and the latter is from Definition \ref{def width}. Conversely if $d$ is an sd-tree decomposition of $H$, then $w(d) = w(\Psi(d))$.

\begin{Prop} \label{prop hypergraph treewidth equal to gamma width}
If $\Gamma = (\Hyp, \text{Trees}, \{K^n_\Hyp \}_{n \geq 0})$ is the sd-category from Proposition \ref{prop sd-category for hypergraph treewidth}, and $H$ is a hypergraph, then
\begin{equation*}
    \mathbf{w}_\Gamma(H) = \mathbf{tw}_\Hyp(H).
\end{equation*}
\end{Prop}

\begin{proof}
Suppose that $d : \smallint T \to \ncat{Hyp}$ is a minimal sd-tree decomposition of $H$. Then $\Psi(d)$ is a minimal RS-tree decomposition of $H$. Indeed, if $(X,T)$ were any other RS-tree decomposition, then
\begin{equation*}
 w(\Psi(d)) = w(d) \leq w(\Phi(X,T)) = w(X,T).
\end{equation*}
But $w(d) = w(\Psi(d))$, and since both decompositions are minimal $\mathbf{w}_\Gamma(H) = \mathbf{tw}_\Hyp(H)$.
\end{proof}

\begin{Cor} \label{cor graph tree width equal to gamma width}
If $\Gamma = (\Gr, \text{Trees}, \{K^n \}_{n \geq 0})$ is the sd-category from Proposition \ref{prop sd-category for treewidth}, and $G$ is a graph, then
\begin{equation*}
    \mathbf{w}_\Gamma(G) = \mathbf{tw}(G).
\end{equation*}
\end{Cor}

\begin{proof}
This follows from Remarks \ref{rem sd-tree to RS-tree graph decomp} and \ref{rem RS-tree decomp to sd-tree decomp for graphs} along with Proposition \ref{prop hypergraph treewidth equal to gamma width}.
\end{proof}

\addcontentsline{toc}{section}{References}
\printbibliography

\end{document}